\documentclass[11pt,a4paper,oneside]{report}
\usepackage[T1]{fontenc}
\usepackage[utf8]{inputenc}
\usepackage{lmodern,amsmath,amssymb,amsthm,mathtools}
\usepackage[margin=28mm,headheight=14pt]{geometry}
\usepackage{enumitem,booktabs,graphicx,microtype,longtable}
\usepackage{xcolor}
\usepackage[colorlinks=true,linkcolor=blue!55!black,citecolor=blue!55!black,urlcolor=blue!55!black]{hyperref}
\hypersetup{pdftitle={Michelson and Logarithmic Operator Contrasts: Spectral Geometry, Tracial Extensions, and Random Products},pdfauthor={Andrej Novikov and Irina Nikolaeva}}
\setlist[enumerate,1]{label=\textup{(\alph*)},leftmargin=2.2em,itemsep=2pt,topsep=3pt}
\newtheorem{theorem}{Theorem}[section]
\newtheorem{proposition}[theorem]{Proposition}
\newtheorem{lemma}[theorem]{Lemma}
\newtheorem{corollary}[theorem]{Corollary}
\newtheorem{conjecture}[theorem]{Conjecture}
\theoremstyle{definition}
\newtheorem{definition}[theorem]{Definition}
\newtheorem{example}[theorem]{Example}
\theoremstyle{remark}
\newtheorem{remark}[theorem]{Remark}
\newcommand{\OLC}{\operatorname{OLC}}
\newcommand{\OMC}{\operatorname{OMC}}
\newcommand{\Tr}{\operatorname{Tr}}
\newcommand{\tr}{\operatorname{tr}}
\newcommand{\supp}{\operatorname{supp}}
\newcommand{\Var}{\operatorname{Var}}
\newcommand{\GL}{\mathrm{GL}}
\newcommand{\R}{\mathbb{R}}
\newcommand{\C}{\mathbb{C}}
\newcommand{\N}{\mathbb{N}}
\newcommand{\T}{\mathbb{T}}
\newcommand{\E}{\mathbb{E}}
\newcommand{\one}{\mathbf{1}}
\newcommand{\cM}{\mathcal{M}}
\newcommand{\cH}{\mathcal{H}}
\newcommand{\tnorm}[1]{|\!|\!| #1 |\!|\!|}

\newcommand{\dd}{\mathrm{d}}
\newcommand{\boxx}{\boxtimes}
\newcommand{\rk}{\operatorname{rk}}
\title{Michelson and Logarithmic Operator Contrasts\\[6pt]
\large Spectral Geometry, Tracial Extensions, and Random Products}
\author{Andrej Novikov\thanks{Novikov Laboratories LLC, Kazan, Russia. Email: \texttt{a.hobukob@gmail.com}.}
\and Irina Nikolaeva\thanks{Lean Technologies Laboratory LLC, Kazan, Russia.}}
\date{Consolidated preprint, 11 November 2026}
\begin{document}
\hypersetup{pageanchor=false}
\maketitle
\clearpage
\pagenumbering{roman}
\hypersetup{pageanchor=true}
\chapter*{Abstract}
The Michelson contrast of a positive invertible operator is a bounded coordinate for its condition number. Its logarithmic lift,
$\OLC(A)=\frac12\log(\|A\|\|A^{-1}\|)$, connects spectral spread with projective Thompson geometry. This preprint collects the Novikov--Nikolaeva series and the associated single-author extensions by Novikov into one treatment. We develop the equality criterion for product subadditivity, its dimension and Jordan-rank thresholds, quantum-state examples, and contrasts associated with unitarily invariant norms. The latter have a sharp additive defect, while the full Schatten family records a logarithmic spectral profile up to homometry. Tracial extensions give characterisations of traces, dimension-sensitive defects, and spectral-edge asymptotics. Unbounded operators are described through critical exponents and finite sections. Tensor, matrix and free independence lead to distinct growth laws for products. Finally, symmetric cones provide a setting for the equality criterion, Lorentz boosts, and Jordan products. Classical ingredients are distinguished from the refinements developed in the series; conjectural asymptotic statements remain explicitly labelled.
\par\medskip
\noindent\textbf{2020 MSC:} 47A30, 15A60, 46L52; 47A10, 46L10, 17C20, 60B20, 81P45.\\
\textbf{Keywords:} Michelson contrast; condition number; Thompson metric; Schatten norm; homometric spectra; von Neumann algebra; free multiplicative convolution; symmetric cone.
\chapter*{Scope and organisation}
\addcontentsline{toc}{chapter}{Scope and organisation}
This is a consolidated research preprint, rather than a new claim of priority for each ingredient. The operator Michelson contrast for positive operators originates in the work of Abed, Nikolaeva and Novikov \cite{ANN24}; that joint work is credited as a prerequisite. The later separate three-author manuscripts on sharp weighted trace inequalities, quantitative centrality and monotonicity defects are not reproduced here.

Chapter~\ref{ch:core} consolidates the logarithmic-contrast manuscript of Nikolaeva and Novikov \cite{NN26}, including its quantum-information additions. Chapter~\ref{ch:norms} presents the finite-matrix continuation, using the later revised text of the Novikov--Nikolaeva manuscript. Chapters~\ref{ch:tracial}, \ref{ch:random} and \ref{ch:cones} incorporate the associated manuscripts originally authored by Novikov alone. Chapter~\ref{ch:infinite} incorporates the joint infinite-contrast manuscript. These provenance distinctions are retained explicitly; collecting the series does not change the authorship history of the component works.

The notation $\OMC,\OLC$ is used throughout; the original quantum-information text used $\mathrm{QMC},\mathrm{QLC}$ for the same functionals. The chapter structure permits local hypotheses to be stated precisely: an arbitrary invertible operator, a positive element, a finite tracial algebra, and a Euclidean Jordan algebra are different settings. In particular, a statement about $|A|$ alone does not transfer a commutativity assertion to $A$.

\paragraph{Reading guide.}
Chapters~\ref{ch:core}--\ref{ch:norms} establish the matrix framework. Chapter~\ref{ch:tracial} extends it to finite traces, Chapter~\ref{ch:infinite} to unbounded spectra, and Chapter~\ref{ch:random} to products. Chapter~\ref{ch:cones} returns to geometry and generalises the matrix equality phenomenon by Jordan rank.

\paragraph{Preparation of this version.}
The consolidation, English-language editing, cross-reference reconciliation and selected algebraic checks were assisted by ChatGPT. The present document is a working preprint and has not been independently refereed as a unified text. A separate author-facing preparation note records the source versions and corrections made during consolidation.
\tableofcontents
\clearpage
\pagenumbering{arabic}
\chapter{Operator Michelson contrast and spectral geometry}\label{ch:core}
\section{The logarithmic coordinate}
For an invertible bounded operator $A$ on a complex Hilbert space put
\[
\kappa(A)=\|A\|\|A^{-1}\|,\qquad
\OMC(A)=\frac{\kappa(A)-1}{\kappa(A)+1},
\]\[
\OLC(A)=\operatorname{artanh}\OMC(A)=\frac12\log\kappa(A).
\]
For positive $A$, with $m=\inf\sigma(A)>0$ and $M=\sup\sigma(A)$, this gives
$\OMC(A)=(M-m)/(M+m)$, the operator Michelson contrast of \cite{ANN24}.
The identity $\operatorname{artanh}((M-m)/(M+m))=\frac12\log(M/m)$ is elementary; condition-number inequalities and the scalar change of coordinate are background, not new inequalities.
The logarithmic formulation developed in \cite{NN26} makes the equality and geometric questions explicit.

Polar decomposition gives $\kappa(A)=\kappa(|A|)$, so these functionals extend from the positive cone to all invertibles. They are invariant under multiplication by a nonzero scalar, left or right multiplication by a unitary, adjoint and inversion. They are nonnegative and vanish exactly on scalar multiples of unitaries. On nonzero noninvertible operators we use the extended conventions $\OMC=1$ and $\OLC=+\infty$; the zero operator is excluded from this convention.

\begin{remark}
Passing to $|A|$ extends a functional, not the Thompson metric, and does not preserve commutation questions. For example, the Pauli matrices $\sigma_x,\sigma_y$ have absolute value $I$ and satisfy $\sigma_x\sigma_y=-\sigma_y\sigma_x$.
\end{remark}

\section{Projective Thompson geometry and continuity}
For positive invertible $A,B$, the Thompson metric and its quotient by positive scalar multiplication are
\[
d_T(A,B)=\|\log(A^{-1/2}BA^{-1/2})\|,\qquad
\delta_T([A],[B])=\inf_{c>0}d_T(A,cB).
\]
This quotient is one half of Hilbert's projective metric. The underlying order and metric results are classical; see \cite{Thompson63,CPR94}.

\begin{proposition}\label{core:metric}
For $A>0$ invertible,
\[
\delta_T([I],[A])=\OLC(A).
\]
Consequently $|\OLC(A)-\OLC(B)|\le\delta_T([A],[B])$.
\end{proposition}
\begin{proof}
Writing $m=\inf\sigma(A)$ and $M=\sup\sigma(A)$,
$d_T(I,cA)=\max\{|\log(cm)|,|\log(cM)|\}$. The minimum occurs at
$c=(mM)^{-1/2}$ and equals $\frac12\log(M/m)$. The second statement is the reverse triangle inequality in the metric on rays.
\end{proof}

\begin{proposition}\label{core:singular}
If invertible $A_j$ converge in norm to a nonzero noninvertible $A$, then
$\kappa(A_j)\to\infty$ and $\OLC(A_j)\to\infty$.
\end{proposition}
\begin{proof}
We have $\|A_j\|\to\|A\|>0$. If $\|A_j^{-1}\|$ had a bounded subsequence, then along it $\|A_j^{-1}(A-A_j)\|<1$ eventually. The Neumann series would make $A=A_j(I+A_j^{-1}(A-A_j))$ invertible, a contradiction.
\end{proof}
On an order interval $aI\le A,B\le bI$ with $a>0$, ordinary norm continuity is also quantitative: extremal eigenvalues are norm-Lipschitz, hence $|\OLC(A)-\OLC(B)|\le\|A-B\|/a$. A bound on condition numbers alone, without a lower scale bound, does not give this norm-Lipschitz conclusion.

\section{Products, equality and dimension}
\begin{theorem}\label{core:subadd}
For invertible $A,B$,
\[
\OLC(AB)\le\OLC(A)+\OLC(B),\qquad
\OMC(AB)\le\frac{\OMC(A)+\OMC(B)}{1+\OMC(A)\OMC(B)}.
\]
\end{theorem}
\begin{proof}
Submultiplicativity of the norm, applied to $AB$ and $B^{-1}A^{-1}$, gives
$\kappa(AB)\le\kappa(A)\kappa(B)$. Take half logarithms, then apply the increasing function $\tanh$ and its addition formula.
\end{proof}

\begin{lemma}\label{core:alignment}
For positive definite $X,Y\in M_n(\C)$, equality in Theorem~\ref{core:subadd} holds if and only if
\[
E_{\max}(X)\cap E_{\max}(Y)\ne\{0\},\qquad
E_{\min}(X)\cap E_{\min}(Y)\ne\{0\}.
\]
\end{lemma}
\begin{proof}
Equality of the product of the two norm bounds forces equality in both bounds. Choose a unit vector $v$ attaining $\|XY\|$. Equality throughout
$\|XYv\|\le\|X\|\|Yv\|\le\|X\|\|Y\|$ implies that $v$ is a top eigenvector of $Y$ and that $Yv$ is a top eigenvector of $X$. Positivity makes $Yv$ a positive scalar multiple of $v$. Apply the same argument to $Y^{-1}X^{-1}$ for the lower endpoint. Conversely, common top and bottom eigenvectors attain both bounds.
\end{proof}

\begin{theorem}\label{core:dimension}
If $X,Y\in M_n(\C)$ are positive definite and
$\OLC(XY)=\OLC(X)+\OLC(Y)$, then $XY=YX$ for $n\le3$.
Noncommuting equality pairs exist in every dimension $n\ge4$.
\end{theorem}
\begin{proof}
If either matrix is scalar there is nothing to prove. Otherwise common top and bottom eigenvectors supplied by Lemma~\ref{core:alignment} are orthogonal. Their spans reduce both self-adjoint matrices. The orthogonal complement has dimension at most one when $n\le3$, so both matrices are diagonal on a common orthonormal basis.
For $n=4$ take
\[
X=\begin{pmatrix}4&0&0&0\\0&3&0&0\\0&0&1&0\\0&0&0&1\end{pmatrix},\qquad
Y=\begin{pmatrix}4&0&0&0\\0&2&1&0\\0&1&2&0\\0&0&0&1\end{pmatrix}.
\]
Both have maximal eigenvalue $4$ and minimal eigenvalue $1$, sharing $e_1$ and $e_4$ at the two endpoints. Thus $\kappa(XY)=16$, while the middle $2\times2$ blocks do not commute. Direct sums with identity blocks preserve the example.
\end{proof}

Commutativity is necessary in low dimension, but is never sufficient by itself: the endpoint order must also align. For example, $X=\operatorname{diag}(a,a^{-1})$ and $Y=X^{-1}$ have scalar product and positive individual contrasts when $a>1$.

\subsection{Tensor products and local contrasts}
The spatial tensor norm and the inverse formula give
\[
\OLC(A\otimes B)=\OLC(A)+\OLC(B),\qquad
\OLC(A^{\otimes N})=N\OLC(A).
\]
\begin{proposition}\label{core:tensor-equality}
Let $X=\bigotimes_iX_i$, $Y=\bigotimes_iY_i$, where $X_i,Y_i$ are positive definite matrices. Equality for $X,Y$ holds if and only if it holds for every pair $X_i,Y_i$. In particular, if all factor dimensions are at most three, equality forces $XY=YX$, regardless of the total dimension.
\end{proposition}
\begin{proof}
The deficit $\OLC(X)+\OLC(Y)-\OLC(XY)$ is the sum of the nonnegative factor deficits. It vanishes precisely when every summand vanishes. Apply Theorem~\ref{core:dimension} in each factor.
\end{proof}
For instance, no noncommuting equality pair can have simultaneous pure tensor decompositions on $\C^3\otimes\C^2$. This is a restriction on the specified factorisation, not an assertion that every noncommuting equality pair is irreducible under arbitrary direct sums.

If both tensor factors are finite-dimensional and $Z$ is invertible, positivity and the order bounds for $|Z|$ imply
\[
\OLC(\Tr_2|Z|)\le\OLC(Z),\qquad
\OLC(\Tr_1|Z|)\le\OLC(Z).
\]
Indeed, $mI\le|Z|\le MI$ gives $d_2mI\le\Tr_2|Z|\le d_2MI$. There is no reverse upper estimate by the local contrasts: for $0<\varepsilon<1$, $Z=I\otimes I+\varepsilon\sigma_z\otimes\sigma_z$ has scalar partial traces and nonzero contrast.
For nearest-neighbour invertible factors on a finite chain,
\[
\OLC(A_{12}A_{23}\cdots A_{N-1,N})\le\sum_{i=1}^{N-1}\OLC(A_{i,i+1}).
\]
This is repeated subadditivity. Products occurring at intermediate steps need not be positive; the positive-matrix alignment lemma must not be applied to them without further hypotheses.

\section{Powers, means, sums and contractions}
For positive invertible $A$ and $z\in\C$, functional calculus gives
\[
|A^z|=A^{\Re z},\qquad \OLC(A^z)=|\Re z|\OLC(A).
\]
Hence, for $0\le t\le1$ and positive invertible $A,B$,
\[
\OLC(A^tB^{1-t})\le t\OLC(A)+(1-t)\OLC(B).
\]
For the weighted geometric mean
$A\#_tB=A^{1/2}(A^{-1/2}BA^{-1/2})^tA^{1/2}$ the corresponding estimate is
\[
\OLC(A\#_tB)\le (1-t)\OLC(A)+t\OLC(B).
\]
Indeed, monotonicity of the weighted geometric mean yields
$m_A^{1-t}m_B^t I\le A\#_tB\le M_A^{1-t}M_B^t I$, with the weights used in the definition of $\#_t$. More generally, this order argument applies to any positive mean with the stated endpoint bounds.

If a map on positive invertibles fixes the identity ray and is nonexpansive for $\delta_T$, Proposition~\ref{core:metric} gives $\OLC(F(A))\le\OLC(A)$. A separate sufficient condition is spectral range compression:
$m_{F(A)}\ge m_A$ and $M_{F(A)}\le M_A$. Unital positive linear maps have this property. Reversing both endpoint inequalities gives contrast expansion.

\begin{proposition}\label{core:sums}
For positive invertible $A_1,\ldots,A_N$ and positive scalars $w_i$,
\[
\OLC\Bigl(\sum_iw_iA_i\Bigr)\le\max_i\OLC(A_i).
\]
\end{proposition}
\begin{proof}
Put $K=\max_i\kappa(A_i)$. The sum lies between
$(\sum_iw_im_i)I$ and $(\sum_iw_iM_i)I$, and
$\sum_iw_iM_i\le K\sum_iw_im_i$.
\end{proof}
Thus Ces\`aro means have contrast bounded by the supremum of the individual contrasts. If such means converge in norm to an invertible limit, continuity passes the bound to the limit. Similarly,
\[
\OLC(A_1\cdots A_N)\le\sum_{i=1}^N\OLC(A_i);
\]
if the products converge in norm to an invertible $P$ and the series on the right converges, then $\OLC(P)\le\sum_{i\ge1}\OLC(A_i)$.

For $A>0$ one also has
\[
\OLC(A)\ge\frac12\log\bigl(1+\|I-A/\|A\|\|\bigr),
\]
with equality exactly for scalar $A$. To see this, set $r=m/M\in(0,1]$; the inequality is $r^{-1}\ge2-r$, equivalent to $(1-r)^2\ge0$.
For $X>0$ and invertible $Y$,
\[
\OLC(X^{1/2}YX^{1/2})\le\OLC(X)+\OLC(Y).
\]
If $Y>0$, the equality case is the same endpoint alignment as in Lemma~\ref{core:alignment}; see the Jordan-algebra proof in Chapter~\ref{ch:cones}.

\section{Quantum states: qubits and higher dimensions}
For a positive definite density matrix $\rho$ on $\C^d$,
$\OLC(\rho)=\frac12\log(\lambda_{\max}/\lambda_{\min})$ and
$\OLC(\rho)\le-\frac12\log\lambda_{\min}$. For $d\ge2$ equality in this last bound is unattainable at a full-rank state, but the bound is asymptotically attained along spectra approaching a pure state. The maximally mixed state has zero contrast; any nonzero rank-deficient state has infinite logarithmic contrast.

\begin{proposition}\label{core:qubit}
For a qubit $\rho=\frac12(I+\boldsymbol r\cdot\boldsymbol\sigma)$ with $r=|\boldsymbol r|<1$,
\[
\OMC(\rho)=r,\qquad \OLC(\rho)=\operatorname{artanh}r,
\]
and, with natural logarithms,
\[
S(\rho)=H(r):=\log2-\frac12\bigl((1+r)\log(1+r)+(1-r)\log(1-r)\bigr).
\]
This is a strictly decreasing bijection from $r\in[0,1)$ to $S\in(0,\log2]$.
Furthermore
\[
\frac12\|\rho-I/2\|_1=\frac r2,\qquad
\|\rho-I/2\|_2=\frac r{\sqrt2},\qquad
\Tr\rho^2=\frac{1+r^2}{2},
\]
where these two Schatten norms are unnormalised.
\end{proposition}
\begin{proof}
The Pauli identity gives $(\boldsymbol r\cdot\boldsymbol\sigma)^2=r^2I$, hence eigenvalues $(1\pm r)/2$ for $\rho$. All formulas follow by substitution. In particular $H'(r)=\frac12\log((1-r)/(1+r))<0$ for $r>0$.
\end{proof}
In terms of $\eta=\OLC(\rho)$, the equivalent formula
$S(\eta)=\log(2\cosh\eta)-\eta\tanh\eta$ gives
\[
S(\eta)=\log2-\frac12\eta^2+O(\eta^4)\quad(\eta\to0),\qquad
S(\eta)=(1+2\eta)e^{-2\eta}\bigl(1+O(e^{-2\eta})\bigr)\quad(\eta\to\infty).
\]
These follow from Taylor expansion at zero and the exact representation
$S(\eta)=\log(1+e^{-2\eta})+2\eta/(1+e^{2\eta})$. There is no exact symmetry under $\eta\mapsto1/\eta$.

\subsection{Entropy at fixed qutrit condition number}
Fix $\kappa>1$. Every ordered qutrit spectrum with this condition number is
\[
(\lambda_1,\lambda_2,\lambda_3)=\bigl(t,1-(1+1/\kappa)t,t/\kappa\bigr),
\quad \frac{\kappa}{2\kappa+1}\le t\le\frac{\kappa}{\kappa+2}.
\]
\begin{theorem}\label{core:qutrit-entropy}
The entropy range is a nondegenerate interval $[S_{\min}(\kappa),S_{\max}(\kappa)]$.
Its minimum and unique maximising spectrum are
\[
S_{\min}=\log(\kappa+2)-\frac{\kappa}{\kappa+2}\log\kappa,
\qquad
(\lambda_1^*,\lambda_2^*,\lambda_3^*)=
\frac{(1,q,\kappa^{-1})}{1+q+\kappa^{-1}},\quad q=\kappa^{-1/(\kappa+1)}.
\]
For every finite $\kappa>1$,
$\log2<S_{\max}(\kappa)<\log3$; as $\kappa\to\infty$, $S_{\min}\to0$ and $S_{\max}\to\log2$.
\end{theorem}
\begin{proof}
Entropy is strictly concave along the nonconstant affine family. Its derivative vanishes exactly when
$(\kappa+1)\log\lambda_2=\kappa\log\lambda_1+\log\lambda_3$,
giving $q$ above; $\kappa^{-1}<q<1$ places the maximum in the interior.
The two endpoint entropies are
$S_P=\log(2\kappa+1)-\frac{2\kappa}{2\kappa+1}\log\kappa$ and $S_Q=S_{\min}$.
Their difference vanishes at $1$ and has derivative
\[
(S_P-S_Q)'=\frac{6(\kappa^2-1)\log\kappa}{(2\kappa+1)^2(\kappa+2)^2}>0.
\]
Thus $S_Q<S_P$. Finally $S_P'(\kappa)=-2\log\kappa/(2\kappa+1)^2<0$ and $S_P\to\log2$, so $S_{\max}>S_P>\log2$; the strict upper bound $\log3$ follows because no admissible spectrum is uniform. The limits follow from the displayed spectra.
\end{proof}
For every $d>3$, entropy likewise cannot be determined by contrast. One direct construction fixes unnormalised extreme eigenvalues $\kappa$ and $1$, puts the other $d-2$ eigenvalues equal to $s\in[1,\kappa]$, and normalises by $Z=\kappa+1+(d-2)s$. Its entropy has derivative
\[
\frac{d-2}{Z^2}\bigl(\kappa\log\kappa-(\kappa+1)\log s\bigr),
\]
which is not identically zero, while its condition number is always $\kappa$.

\subsection{Purity at fixed qutrit condition number}
\begin{proposition}\label{core:qutrit-purity}
On the same family,
\[
P_{\min}(\kappa)=\frac{\kappa^2+1}{2(\kappa^2+\kappa+1)},\qquad
P_{\max}(\kappa)=\frac{\kappa^2+2}{(\kappa+2)^2}.
\]
The minimum is attained at $t_*={\kappa(\kappa+1)}/{[2(\kappa^2+\kappa+1)]}$, and the maximum at the endpoint with $\lambda_2=\lambda_3$.
\end{proposition}
\begin{proof}
Put $a=1+1/\kappa$ and $b=1+a^2+1/\kappa^2$. Then
$P(t)=bt^2-2at+1$ is strictly convex, with vertex $a/b$ strictly between the endpoints for $\kappa>1$. Substituting gives $P_{\min}=1-a^2/b$. The difference between the two endpoint purities, with the larger-eigenvalue multiplicity one first, is
\[
\frac{\kappa^2+2}{(\kappa+2)^2}-\frac{2\kappa^2+1}{(2\kappa+1)^2}
=\frac{2(\kappa-1)^3(\kappa+1)}{(\kappa+2)^2(2\kappa+1)^2}>0.
\]
\end{proof}
For $\kappa>1$, $1/3<P_{\min}<1/2$ and $1/3<P_{\max}<1$. Setting $\kappa=e^{2\eta}$ gives
\[
1-P_{\max}=4e^{-2\eta}(1+O(e^{-2\eta})),\quad
\tfrac12-P_{\min}=\tfrac12e^{-2\eta}(1+O(e^{-2\eta})).
\]
The second limit approaches a rank-two mixed state, not a pure state.

\subsection{Separability and ensemble mixing}
The matrices in Theorem~\ref{core:dimension}, normalised to trace one, give a noncommuting equality pair without requiring entanglement. Indeed, $X/9$ is diagonal in the computational two-qubit basis and the partial transpose of $Y/9$ is
\[
\frac19\begin{pmatrix}4&0&0&1\\0&2&0&0\\0&0&2&0\\1&0&0&1\end{pmatrix},
\]
with eigenvalues $(5\pm\sqrt{13})/18,2/9,2/9$, all positive. The two-qubit PPT criterion therefore implies separability \cite{Horodecki96}. The tensor restriction of Proposition~\ref{core:tensor-equality} concerns simultaneous pure tensor decompositions and is compatible with this separable example.

The classical link between Michelson contrast and Jensen--Shannon divergence was studied in \cite{BRV12}. The qubit formulas below use the same binary-entropy function after passing to eigenvalues.

For a binary ensemble of qubits with Bloch-vector lengths $r_1,r_2$ and angle $\theta$, put
\[
\bar r=\frac12\sqrt{r_1^2+r_2^2+2r_1r_2\cos\theta}.
\]
Then the quantum Jensen--Shannon divergence is exactly
\[
\operatorname{QJSD}(\rho_1,\rho_2)=H(\bar r)-\tfrac12H(r_1)-\tfrac12H(r_2).
\]
This follows because the Bloch vector of the average state is the average of the Bloch vectors. Parallel vectors reduce the formula to a classical Bernoulli divergence; antiparallel unit vectors give $\log2$. For the antipodal pair $\rho_\pm=(I\pm r\sigma_z)/2$, the formula becomes $\operatorname{QJSD}(\rho_+,\rho_-)=\log2-H(r)$, a strictly increasing function of their common Michelson contrast $r$ on $[0,1]$. This is the binary entropy--contrast correspondence. For a general pair, the angle remains an additional variable even for qubits.
In any dimension Proposition~\ref{core:sums} gives
$\OLC(\sum_ip_i\rho_i)\le\max_i\OLC(\rho_i)$.
The elementary inequality $\OLC(\rho)\ge\frac12(S(\rho)-\log d)$ is valid but vacuous as a lower bound, since $S\le\log d$; combining it with the mixing inequality gives no improvement over the standard binary information bound $\operatorname{QJSD}\le\log2$. Useful entropy-deficit bounds in terms of $C_p$ appear in Chapter~\ref{ch:norms}.

\section{Log-spectral measures and the scope of analogies}
For a positive $n\times n$ matrix define
$\nu_A=n^{-1}\sum_j\delta_{\log\lambda_j(A)}$. Then
\[
\OLC(A)=\tfrac12\operatorname{diam}\supp\nu_A,\qquad
\widehat\nu_A(t)=n^{-1}\Tr A^{it},\qquad
\nu_{A\otimes B}=\nu_A*\nu_B.
\]
The Fourier transform of this finite atomic measure is a finite exponential sum; it need not decay at infinity. For two eigenvalues, after subtracting their logarithmic midpoint, it is $\cos(t\OLC(A))$. Before that centring it has an additional scalar phase. These are exact spectral identities. They contain no information about relative eigenvectors and therefore cannot by themselves detect whether two operators commute. Chapter~\ref{ch:norms} develops the corresponding Laplace profile and its homometric ambiguity.

The scalar relation $\OMC=\tanh\OLC$ resembles the velocity--rapidity relation, but does not identify arbitrary operator geometry with spacetime geometry. A precise Lorentz-cone statement, with its product law and rank hypotheses, is given in Chapter~\ref{ch:cones}. Analogies with neural-network activations, colour-channel counts, representation dimensions, or spin-chain gap criteria are not used as proofs or as results of this preprint.

\chapter{Contrasts beyond the operator norm}\label{ch:norms}
\section{Symmetric gauges and normalised contrast}
Throughout this chapter $n\ge2$, unless a one-dimensional case is explicitly mentioned. Let $\Phi$ be a symmetric gauge function on $\R^n$, and write
$\tnorm{A}_\Phi=\Phi(s(A))$ for the associated unitarily invariant norm. Here
$s_1(A)\ge\cdots\ge s_n(A)>0$ are the singular values of an invertible matrix,
$\one=(1,\ldots,1)$, and $e_j$ is a coordinate vector. Put
\[
C_\Phi(A)=\frac12\log\frac{\Phi(s(A))\Phi(s(A)^{-1})}{\Phi(\one)^2},\qquad
d_\Phi=\log\frac{\Phi(\one)}{\Phi(e_1)}.
\]
The reciprocal vector is understood up to permutation. We use standard symmetric-gauge and majorisation facts \cite{Bhatia97}:
\begin{align*}
|x|\le|y|&\Longrightarrow\Phi(x)\le\Phi(y),\\
\Phi(e_1)\|x\|_\infty&\le\Phi(x)\le\Phi(\one)\|x\|_\infty,\\
\frac{\Phi(\one)}n\|x\|_1&\le\Phi(x)\le\Phi(e_1)\|x\|_1,\\
\Phi(x^t y^{1-t})&\le\Phi(x)^t\Phi(y)^{1-t}\quad(x,y\ge0, 0\le t\le1).
\end{align*}
The last inequality follows by homogeneity from the coordinatewise weighted arithmetic--geometric mean inequality. The lower $\ell_1$ bound follows by averaging all coordinate permutations. A symmetric convex function is monotone under majorisation because a majorised vector is a convex combination of permutations.

\begin{proposition}\label{norm:basic}
The contrast is unchanged by scaling $\Phi$, by $A\mapsto\lambda UAV$ with $U,V$ unitary and $\lambda\ne0$, and by adjoint, inversion or absolute value. Moreover
\[
0\le C_1(A)\le C_\Phi(A)\le\OLC(A),\qquad
C_\Phi(A)\ge\OLC(A)-d_\Phi.
\]
It vanishes exactly on scalar multiples of unitaries.
\end{proposition}
\begin{proof}
The invariances follow from the singular values. The gauge bounds give
$\Phi(s)\Phi(s^{-1})/\Phi(\one)^2\ge\mathrm{AM}(s)/\mathrm{HM}(s)\ge1$,
with equality in the latter inequality only for constant $s$. The upper and lower $\ell_\infty$ estimates give the two comparisons with $\OLC$.
\end{proof}

\section{The sharp defect and rigidity}\label{norm:sec:defect}
\begin{theorem}\label{norm:defect}
For $n\ge2$, invertible $A,B$, positive definite $P$, and $s\ge1$,
\begin{align*}
C_\Phi(AB)&\le C_\Phi(A)+C_\Phi(B)+d_\Phi,\\
sC_\Phi(P)&\le C_\Phi(P^s)\le sC_\Phi(P)+(s-1)d_\Phi.
\end{align*}
Both constants are optimal: for $P_t=\operatorname{diag}(t,1,\ldots,1,t^{-1})$,
\[
\lim_{t\to\infty}\bigl(C_\Phi(P_t^s)-sC_\Phi(P_t)\bigr)=(s-1)d_\Phi.
\]
\end{theorem}
\begin{proof}
Normalise $\Phi(\one)=1$ and set $a=\Phi(e_1)$. The ideal norm inequality gives
$\tnorm{AB}_\Phi\le\|A\|\tnorm{B}_\Phi\le a^{-1}\tnorm{A}_\Phi\tnorm{B}_\Phi$.
Apply it also to the inverses and take half logarithms.
For an eigenvalue vector $\lambda>0$, the gauge H\"older inequality with exponent $1/s$ gives $\Phi(\lambda^s)\ge\Phi(\lambda)^s$. Coordinatewise monotonicity gives
$\Phi(\lambda^s)\le\|\lambda\|_\infty^{s-1}\Phi(\lambda)\le a^{1-s}\Phi(\lambda)^s$.
Repeat for $\lambda^{-1}$. Finally $t^{-s}s(P_t^s)\to e_1$ and the corresponding inverse vector tends to $e_n$, so
$C_\Phi(P_t^s)=s\log t+\log a+o(1)$. Subtraction proves optimality.
\end{proof}
The norm-product estimate is the classical submultiplicativity of a gauge normalised by $\Phi(e_1)=1$. The power inequalities are related to Marshall--Olkin's condition-number inequalities \cite{MO65}; the normalisation and optimal additive constants are the point of this formulation. For $0<s\le1$, the equivalent bounds are
$sC_\Phi(P)-(1-s)d_\Phi\le C_\Phi(P^s)\le sC_\Phi(P)$.

\begin{corollary}\label{norm:rigidity}
The following are equivalent: product subadditivity of $C_\Phi$ for all invertibles; the inequality $C_\Phi(P^2)\le2C_\Phi(P)$ for all $P>0$; exact positive-power homogeneity; $d_\Phi=0$; and $\Phi$ being a positive multiple of the operator-norm gauge.
\end{corollary}
\begin{proof}
If $d_\Phi>0$, the family $P_t$ violates square subadditivity. If $d_\Phi=0$, then $\Phi(e_1)=\Phi(\one)$, and the two $\ell_\infty$ bounds force $\Phi(x)=\Phi(e_1)\|x\|_\infty$. All remaining implications follow from the operator-norm power law and Theorem~\ref{norm:defect}.
\end{proof}
Applying the power bound to $|A|^2=A^*A$ gives
$2C_\Phi(A)\le C_\Phi(A^*A)\le2C_\Phi(A)+d_\Phi$.
The full range of defects is $0\le d_\Phi\le\log n$, with the endpoints attained exactly by multiples of the operator and trace norms, respectively. For the latter assertion, $\Phi(\one)=n\Phi(e_1)$ makes the two $\ell_1$ bounds coincide.
At fixed individual contrasts, the product excess is also bounded by
$\min\{d_\Phi,\OLC(A)+\OLC(B)\}$. If
$\Phi(x)>\Phi(e_1)\|x\|_\infty$ whenever $x$ has at least two nonzero coordinates, the defect is not attained by invertibles. It can be attained for other gauges: for $n=4$,
$\Phi(x)=\max\{\|x\|_1/2,\|x\|_\infty\}$ and $P=\operatorname{diag}(3,1,1,1/3)$ give $C_\Phi(P^2)-2C_\Phi(P)=\log2=d_\Phi$.

\section{Schatten contrasts and logarithmic spectral profiles}\label{norm:profile}
Let $\tau=n^{-1}\Tr$, $x_i=\log s_i(A)$, and
\[
\nu_A=\frac1n\sum_i\delta_{x_i},\quad
L_A(t)=\log\int e^{tx}\,\dd\nu_A(x),\quad
C_p(A)=\frac{L_A(p)+L_A(-p)}{2p}\quad(p>0).
\]
For $p\ge1$ this is the contrast for the Schatten norm; the same moment formula defines it for $0<p<1$. Set $C_\infty=\OLC$.
The tilted probability measure is $\dd\nu_A^{(t)}=e^{tx-L_A(t)}\dd\nu_A$.
Standard cumulant identities \cite{DZ98} yield
\[
L_A'(t)=\int x\,\dd\nu_A^{(t)},\qquad
L_A''(t)=\Var(\nu_A^{(t)})=\tfrac12\iint(x-y)^2\,\dd\nu_A^{(t)}(x)\dd\nu_A^{(t)}(y).
\]
Thus $pC_p$ is the even part of $L_A$, unaffected by centring the log-spectrum, and $\OLC=\frac12(L_A'(+\infty)-L_A'(-\infty))$.

\begin{theorem}\label{norm:schatten}
For $p>0$ put $T_p=e^{2pC_p}$. Then
\begin{align*}
T_p(A)&=\iint\cosh(p(x-y))\,\dd\nu_A(x)\dd\nu_A(y),\\
\frac{\dd}{\dd p}C_p(A)&=\frac{D(\nu_A^{(p)}\|\nu_A)+D(\nu_A^{(-p)}\|\nu_A)}{2p^2}\ge0,\\
\lim_{p\downarrow0}\frac{C_p(A)}p&=\tfrac12\Var(\nu_A),\qquad
\lim_{p\to\infty}C_p(A)=\OLC(A),\\
C_p(|A|^s)&=|s|C_{|s|p}(A)\quad(s\ne0),\qquad
C_p(A\otimes B)=C_p(A)+C_p(B).
\end{align*}
The function $C_p$ is real-analytic, and its derivative is strictly positive unless $A$ is a scalar multiple of a unitary.
\end{theorem}
\begin{proof}
Symmetrise the product of the two moment integrals to obtain the cosh formula. Since
$D(\nu_A^{(t)}\|\nu_A)=tL_A'(t)-L_A(t)$, differentiation gives the second formula. Relative entropy vanishes here only for a point mass. Taylor expansion gives
$L_A(p)+L_A(-p)=p^2\Var(\nu_A)+O(p^4)$; at infinity the largest exponential in each sum gives the limit. The power map sends the measure by $x\mapsto sx$, while tensor products convolve the two log-spectral measures and add their cumulant profiles.
\end{proof}
In particular $C_p(P^s)>sC_p(P)$ for nonscalar $P>0$ and $s>1$. The product defect is $\log n/p$, and this remains valid for $0<p<1$: the singular-value estimate $s_j(AB)\le\|A\|s_j(B)$ and the bound $\|A\|\le n^{1/p}\|A\|_{p,\tau}$ give the same proof without using norm convexity.

\section{Endpoint asymptotics and sharp bounds}
Let $x_1>\cdots>x_r$ be the distinct log-singular values, of multiplicities $m_1,\ldots,m_r$, with $r\ge2$. Set $g_+=x_1-x_2$ and $g_-=x_{r-1}-x_r$.
\begin{theorem}\label{norm:rate}
For every $p>0$,
\[
C_p(A)=\OLC(A)-\frac1{2p}\log\frac{n^2}{m_1m_r}+\frac{R_A(p)}{2p},
\]
where
\begin{align*}
0<R_A(p)&\le\frac{n-m_1}{m_1}e^{-pg_+}+\frac{n-m_r}{m_r}e^{-pg_-},\\
R_A(p)&=\frac{m_2}{m_1}e^{-pg_+}(1+o(1))+\frac{m_{r-1}}{m_r}e^{-pg_-}(1+o(1)).
\end{align*}
Hence $p(\OLC-C_p)\to\log n-\frac12\log(m_1m_r)$.
\end{theorem}
\begin{proof}
Factor the top exponential out of each moment sum. The remaining logarithms are
\[
R_A(p)=\log\Bigl(1+\sum_{j\ge2}\frac{m_j}{m_1}e^{-p(x_1-x_j)}\Bigr)
+\log\Bigl(1+\sum_{j<r}\frac{m_j}{m_r}e^{-p(x_j-x_r)}\Bigr).
\]
Use $\log(1+u)\le u$ for the bound and $\log(1+u)\sim u$ for the asymptotics.
\end{proof}
The same elementary Laplace asymptotics give
\[
L_A''(t)\sim g_+^2\frac{m_2}{m_1}e^{-g_+t},\quad
x_1-L_A'(t)\sim g_+\frac{m_2}{m_1}e^{-g_+t}\quad(t\to\infty),
\]
and the analogous lower-edge formula. Thus
$g_\pm=-\lim_{t\to\pm\infty}|t|^{-1}\log L_A''(t)$.
Indeed, in the tilted variance the pair of levels $1,2$ supplies the leading term; every other pair has smaller exponential order.

\paragraph{Gibbs dictionary.}
For $A_\beta=e^{-\beta H}$ and $Z(s)=\Tr e^{-sH}$,
\[
L_{A_\beta}(t)=\log Z(\beta t)-\log n,\quad
L_{A_\beta}''(t)=\beta^2\Var_{\beta t}(H),\quad
\OLC(A_\beta)=\tfrac\beta2(E_{\max}-E_{\min}).
\]
Also $C_p(A_\beta)=[\log Z(p\beta)+\log Z(-p\beta)-2\log n]/(2p)$.
These are the classical partition-function and low-temperature gap formulas \cite{Pathria11}. Tensor powers multiply $L_A$ and the contrasts by $N$, but leave the edge gaps unchanged; the adjacent-to-extreme multiplicity ratios acquire a factor $N$.

\begin{theorem}\label{norm:bounds}
Let $z=\OLC(A)$, $p>0$, and
$\theta_n=4\lfloor n/2\rfloor\lceil n/2\rceil/n^2$.
Then
\begin{align*}
\frac1p\log\frac{2\cosh(pz)+n-2}{n}
&\le C_p(A)\le\frac1{2p}\log\bigl(1+\theta_n\sinh^2(pz)\bigr)\\
&\le\frac1p\log\cosh(pz).
\end{align*}
For $z>0$, equality in the first inequality holds exactly when all interior singular values equal $\sqrt{s_1s_n}$. Equality in the middle bound holds exactly when all singular values are endpoints with multiplicities $\lfloor n/2\rfloor,\lceil n/2\rceil$, in either order.
\end{theorem}
\begin{proof}
By scaling, take the log-endpoints to be $z,-z$. The product
$F(y)=(\sum e^{py_i})(\sum e^{-py_i})$ is strictly convex in each coordinate on $[-z,z]^n$, so its maximum is at a vertex. With $\ell$ top endpoints,
$F=n^2+4\ell(n-\ell)\sinh^2(pz)$, maximised precisely by the balanced multiplicities.
For the minimum, put $X=\sum_{i=2}^{n-1}e^{py_i}$ and $W=\sum_{i=2}^{n-1}e^{-py_i}$. Then $XW\ge(n-2)^2$ and $X+W\ge2(n-2)$, whence
$F=(2\cosh(pz)+X)(2\cosh(pz)+W)\ge(2\cosh(pz)+n-2)^2$.
Both inequalities are equalities precisely when all interior $y_i=0$. The $n=2$ case follows directly.
\end{proof}
The upper bound at $p=1$ is the sharp finite-dimensional Schweitzer--Kantorovich inequality, including the odd-dimensional correction due to Lupa\c{s} \cite{Schweitzer14,Kantorovich48,Lupas72}. General $p$ follows by applying it to $|A|^p$. Here $\theta_n=1$ for even $n$ and $1-n^{-2}$ for odd $n$. The lower bound is the companion fixed-endpoint minimisation.
Consequences include
\[
\sqrt{T_p}\le\cosh(p\OLC)\le\frac{n\sqrt{T_p}-(n-2)}2,
\quad C_p\le\OLC\le C_p+\frac{\log n}p,
\]
and $C_p\le\min\{\OLC,p\OLC^2/2\}$. In dimension two,
$C_p=p^{-1}\log\cosh(p\OLC)$ exactly.

For a positive definite density matrix $\rho$ and $p>1$, monotonicity of R\'enyi entropy \cite{Renyi61} and Jensen's inequality give
\[
\log d-S(\rho)\le\frac{2p}{p-1}C_p(\rho).
\]
Indeed $\tau(\rho^{-p})\ge d^p$, so $\Tr\rho^p\le d^{1-p}T_p(\rho)$ and
$S\ge S_p\ge\log d-2pC_p/(p-1)$. In particular
\[
\log d-S(\rho)\le4C_2(\rho)\le2\log\cosh(2\OLC(\rho))\le4\OLC(\rho)^2.
\]
Letting $p\to\infty$ gives the weaker linear bound $\log d-S\le2\OLC$.

\section{Ky Fan contrasts}\label{norm:kyfan}
For the sum-of-$k$-largest-coordinates gauge $\Phi_{(k)}$, write $C_{(k)}=C_{\Phi_{(k)}}$. Thus
\[
e^{2C_{(k)}(A)}=\Bigl(\frac1k\sum_{i=1}^ks_i\Bigr)
\Bigl(\frac1k\sum_{i=1}^ks_{n+1-i}^{-1}\Bigr).
\]
This is the ratio of the arithmetic mean of the top $k$ values to the harmonic mean of the bottom $k$ values.
\begin{theorem}\label{norm:kyfan-bounds}
The contrasts decrease as $k$ increases:
$\OLC=C_{(1)}\ge\cdots\ge C_{(n)}=C_1\ge0$.
Their product defect is $\log k$ and their power defect is $(s-1)\log k$ for $s\ge1$, both sharp.
Writing $z=\OLC(A)$, $a=\min\{k,\lceil n/2\rceil\}$ and $b=\min\{k,\lfloor n/2\rfloor\}$,
\[
e^{2C_{(k)}}\le\frac{(ae^z+(k-a)e^{-z})(be^z+(k-b)e^{-z})}{k^2}.
\]
For $z>0$ and $2k\le n$, equality holds exactly when the top $k$ singular values coincide and the bottom $k$ coincide; the right side is $e^{2z}$. For $2k>n$, equality requires precisely the two endpoints with balanced multiplicities.
The sharp lower bounds are
\[
C_{(k)}\ge\log\frac{e^z+k-1}{k}\ (k<n),\qquad
C_{(n)}\ge\log\frac{2\cosh z+n-2}{n}.
\]
For $n\ge3$, $2\le k\le n$ and $z>0$, lower equality holds exactly when $s_2=\cdots=s_{n-1}=\sqrt{s_1s_n}$.
\end{theorem}
\begin{proof}
The average of the top $k$ coordinates is decreasing in $k$; apply this to $s$ and $s^{-1}$. The defects follow from Theorem~\ref{norm:defect}.
Scale $s_1=e^z$, $s_n=e^{-z}$. Let $T$ and $B$ be the sums of the largest $k$ entries of $s$ and $s^{-1}$. Holding all but one coordinate $x$ fixed, their product is the maximum of
\[
\alpha\beta,\quad \alpha(\beta'+x^{-1}),\quad
(\alpha'+x)\beta,\quad(\alpha'+x)(\beta'+x^{-1}),
\]
where $\alpha,\alpha'$ are the top-$k$ and top-$(k-1)$ sums of the other coordinates, and $\beta,\beta'$ the reciprocal counterparts. At $k=n$ only the last expression occurs. Each is convex in $x>0$, so a maximum occurs at cube vertices. A vertex with $\ell$ top endpoints has
\[
T=\min(k,\ell)e^z+(k-\min(k,\ell))e^{-z},\quad
B=\min(k,n-\ell)e^z+(k-\min(k,n-\ell))e^{-z}.
\]
For $2k\le n$, the maximal product is $k^2e^{2z}$ and its equality conditions follow directly from $T,B\le ke^z$. For $2k>n$, the product is strictly concave in $\ell$ on $[n-k,k]$ and decreases beyond that interval, giving the balanced vertices.
There are no other maximisers in this case. An interior coordinate at a maximum would make the coordinatewise convex function constant throughout its interval. Since the nonconstant branches above are strictly monotone or strictly convex, the constant branch $\alpha\beta$ would dominate throughout. Evaluating at both endpoints would require at least $k$ other entries equal to $e^z$ and $k$ equal to $e^{-z}$, impossible when $2k>n$.

For the lower bound with $2\le k<n$, let $u_1\ge\cdots\ge u_{n-2}$ be the interior values. Write $T=e^z+X$, $B=e^z+W$, where $X$ sums the top $k-1$ interior values and $W$ sums reciprocals of the bottom $k-1$. If $Y$ sums those bottom values, then $Y\le X$ and $YW\ge(k-1)^2$. Therefore
$TB\ge(e^z+X)(e^z+(k-1)^2/X)\ge(e^z+k-1)^2$.
Equality forces the bottom values equal, $X=Y$, and $X=k-1$; the ordered top and bottom blocks then force every interior value to be one. The case $k=n$ is Theorem~\ref{norm:bounds} at $p=1$.
\end{proof}
For $k<n$ this also yields
$C_{(k)}\le\OLC\le\log(ke^{C_{(k)}}-k+1)$.
The Ky Fan family and the Schatten family retain different information. For
$A=\operatorname{diag}(11,11,7,1)$ and $B=\operatorname{diag}(12,9,2,1)$, all four Ky Fan contrasts of $A$ are smaller than those of $B$, while $C_3(A)\approx0.8734>0.8595\approx C_3(B)$.

\section{Homometric spectra and identifiability}\label{norm:homometry}
For a finitely supported probability measure $\nu$, let $\check\nu(E)=\nu(-E)$. Measures are homometric if $\nu*\check\nu$ agrees. Translation and reflection give the trivial examples. Matrices are isocontrast if their $C_p$ agree for every $p>0$.
\begin{theorem}\label{norm:homometric}
For $A\in\GL_n(\C)$ and $B\in\GL_m(\C)$, the following are equivalent: isocontrast; equality of $C_p$ on a subset of $(0,\infty)$ having a finite accumulation point; and homometry of $\nu_A,\nu_B$. In that case their operator logarithmic contrasts also agree.
\end{theorem}
\begin{proof}
$T_p(A)=\int e^{pu}\,\dd(\nu_A*\check\nu_A)(u)$ is an entire finite exponential sum. Agreement on an accumulating set implies agreement everywhere, and distinct real exponentials are linearly independent. The largest point in the support of the autocorrelation is $2\OLC$.
\end{proof}
For matrices of the same size, this is equality of the multisets of pairwise nonnegative logarithmic differences. This is the classical turnpike problem \cite{Patterson44,Piccard39,Bloom77,BG07,LSS03}; the spectral interpretation should be distinguished from its combinatorial ingredients.

\begin{theorem}\label{norm:threshold}
If $n\le5$ and $A,B\in\GL_n(\C)$ are isocontrast, their singular-value multisets agree up to scaling or scaling and inversion. For every $n\ge6$ there are positive definite counterexamples. They can have simple spectrum for $6\le n\le11$ and $n=13$.
\end{theorem}
\begin{proof}
For the first assertion see the complete small-cardinality argument in Section~\ref{norm:small-proof} below.
For the second, exponentiate a nontrivially homometric pair $X,X'$ of integer multisets: $A=\operatorname{diag}(q^x)_{x\in X}$, $B=\operatorname{diag}(q^x)_{x\in X'}$, $q>1$.
The following simple sets have equal difference multisets and are neither translates nor reflected translates:
\begin{center}\small
\begin{tabular}{c l l}
\toprule $n$&$X$&$X'$\\\midrule
6&$0,1,4,10,12,17$&$0,1,8,11,13,17$\\
7&$0,1,2,5,7,9,12$&$0,1,5,7,8,10,12$\\
8&$0,1,2,5,6,8,10,13$&$0,1,5,6,8,9,11,13$\\
9&$0,1,4,5,6,7,8,9,11$&$0,1,2,3,4,6,7,8,11$\\
10&$0,1,3,6,7,8,9,10,12,14$&$0,1,2,3,5,7,8,9,11,14$\\
11&$0,1,2,3,5,6,7,8,11,14,15$&$0,1,2,3,4,6,7,9,10,14,15$\\
13&$0,1,2,3,5,6,7,8,12,13,15,16,17$&$0,1,2,3,4,5,7,10,11,12,15,16,17$\\
\bottomrule
\end{tabular}
\end{center}
The first pair is Bloom's example \cite{Bloom77}, with difference multiset $\{1,2,\ldots,13,16,17\}$.
For the remaining sizes use Laurent polynomials. If $f=PQ$ and $f'=P(z)Q(z^{-1})$ have nonnegative integer coefficients, then $ff^*=f'f'^*$, where $f^*(z)=f(z^{-1})$. If neither $P$ nor $Q$ is reciprocal up to a signed monomial, cancellation in $\mathbb Z[z^{\pm1}]$ rules out trivial homometry.
For $n=rs$, $r,s\ge3$, use $P=1+(r-1)z$, $Q=1+(s-1)z$.
For $n=3a+2b$, $1\le a<b\le2a$, use
$P=1-z+z^3$, $Q=a+az+bz^2+bz^3+az^4$; the coefficient vectors of $f,z^4f'$ are
\[
(a,0,b-a,a,2a-b,b-a,b,a),\quad
(a,b-a,0,2a-b,b,b-a,a,a).
\]
For $n=23$, use $P=1-z+z^3$, $Q=3+4z+7z^2+5z^3+4z^4$, producing vectors
$(3,1,3,1,3,3,5,4)$ and $(4,1,2,1,4,4,4,3)$.
These cover all sizes. For $n\ge35$, the half-open interval $[n/7,n/5)$ has length at least two, so it contains an integer $a\equiv n\pmod2$; then $b=(n-3a)/2$ satisfies $a<b\le2a$. For $6\le n\le34$, direct enumeration leaves only $6,8,10,11,13,23$ outside the product and two-parameter families, and these were supplied explicitly.
\end{proof}

\begin{proposition}\label{norm:generic}
Suppose $A$ has simple singular values whose successive logarithmic gaps are linearly independent over $\mathbb Q$. If $B\in\GL_m(\C)$ is isocontrast with $A\in\GL_n(\C)$, then $n$ divides $m$, and the singular values of $B$ are those of $cA$ or $(cA)^{-1}$, each repeated $m/n$ times. The gap hypothesis holds outside a countable union of hyperplanes.
\end{proposition}
\begin{proof}
The gap group $\Gamma$ is free abelian of rank $n-1$ and its real group algebra is the Laurent polynomial UFD $\R[t_1^{\pm1},\ldots,t_{n-1}^{\pm1}]$. After translation the measure of $A$ corresponds to
\[
f=\frac1n(1+t_{n-1}+t_{n-1}t_{n-2}+\cdots+t_{n-1}\cdots t_1).
\]
It is irreducible: as a polynomial in $t_1$ it has degree one with unit leading coefficient; in any factorisation the factor independent of $t_1$ must divide that unit. The same holds after Laurent localisation, and for $f^*$.
Every difference in $\supp\nu_B$ belongs to $\Gamma$, so after translation $\nu_B$ is represented by a nonnegative Laurent polynomial $h$ of mass one. Homometry says $hh^*=ff^*$. Unique factorisation forces $h$ to be a unit times $f$ or $f^*$: taking both irreducible factors, or neither, would make the other member of $hh^*$ a unit, an impossibility. Mass one fixes the scalar of the unit at one. Thus the measures agree up to translation/reflection, and each of the $n$ atoms has multiplicity $m/n$. Rational gap dependence is exactly a countable union of integer-coefficient hyperplanes.
\end{proof}
Factor reversal is the classical algebraic mechanism for integer homometry \cite{RS82}. For example, $P=Q=1+2z$ gives multiplicities $(1,4,4)$ and $(2,5,2)$ on three successive levels, a nine-dimensional isocontrast pair. Even the inverse ambiguity can change von Neumann entropy: normalised spectra proportional to $(1,1,4)$ and $(4,4,1)$ are isocontrast but have different entropy. For Bloom's pair with $q=2^{1/4}$, the normalised entropies are approximately $1.3272$ and $1.4076$, while $C_{(2)}$ is approximately $1.2604$ and $1.2876$. Thus Ky Fan contrasts can distinguish a pair invisible to all $C_p$.
The full profile $L_A(t)$, unlike its even part, determines $\nu_A$ uniquely: the top level is $\lim_{t\to\infty}L_A(t)/t$, its mass is $\lim e^{L_A(t)-x_1t}$, and subtraction of that exponential permits iteration. This is a classical, potentially ill-conditioned exponential reconstruction problem.

\section{Small-cardinality reconstruction}\label{norm:small-proof}
\begin{proof}[Proof of Theorem~\ref{norm:threshold} for $n\le5$]
Write the nonnegative successive gap vector in increasing-point order, and let $S$ be its sum, the largest pairwise distance. Orient it so that its last gap is at most its first. For $n=1$ the claim is trivial, and for $n=2$ the gap is $S$. For $n=3$, with gaps $(a,b)$ and $b\le a$, removal of one copy of $S$ leaves $\{a,b\}$, so both are determined.
For $n=4$, with gaps $(a,b,c)$, $c\le a$, the largest distance after removing one $S$ is $a+b$; this determines $c$. The sum of all six distances is $3S+b$, determining $b$ and then $a$.
For $n=5$, write $(a,b,c,d)$, $d\le a$. The largest distance after one $S$ is removed is $a+b+c$, determining $d$. The total of all ten distances is $4S+2(b+c)$; hence $K=b+c$ and $a=S-d-K$ are determined. Removing the known six entries $S,a+K,K+d,K,a,d$ leaves
\[
R(b)=\{b,K-b,a+b,K+d-b\}.
\]
If another normalised vector $(a,b',K-b',d)$ has the same distances, equality of the sums of squares of $R(b),R(b')$ gives either $b'=b$ or
$b+b'=K+(d-a)/2$. In the latter case $b'\in R(b)$.
If $b'=K-b$, this equation forces $a=d$ and the vectors are reversals. If $b'=a+b\ne b$, then $a>0$ and $K=2b+(3a-d)/2$. Remove the common entry $a+b$ and subtract $b$ from the remaining entries to get
\[
\{(a-d)/2,2a,(a+d)/2\}=\{0,(3a-d)/2,(3a+d)/2\}.
\]
Since $0\le d\le a$ and $a>0$, the zero can match only $(a-d)/2$, so again $a=d$ and the vectors reverse. Finally, $b'=K+d-b$ forces $d=-a$, hence $a=d=0$, reducing to the first case. All reasoning is multiset reasoning and allows repeated points and zero gaps.
\end{proof}

\section{Geodesic contrasts}\label{norm:geodesic}
An alternative to a product of norms is
\[
G_\Phi(A)=\min_{c\in\R}\frac{\Phi(\log s(A)-c\one)}{\Phi(\one)}.
\]
The minimum exists by coercivity. Geometrically it is the distance from $|A|$ to the scalar ray in the invariant Finsler metric associated with $\Phi$, divided by $\Phi(\one)$ \cite{Bhatia03,Conde09,CL10}.
\begin{theorem}\label{norm:G}
The functional $G_\Phi$ has the scaling, unitary, adjoint and inverse invariances of $C_\Phi$, and vanishes exactly on scalar multiples of unitaries. For real $s$,
\[
G_\Phi(|A|^s)=|s|G_\Phi(A),\qquad
G_\Phi(AB)\le G_\Phi(A)+G_\Phi(B).
\]
Moreover
\[
\frac{\Phi(e_1+e_2)}{\Phi(\one)}\OLC(A)\le G_\Phi(A)\le\OLC(A),
\]
and the lower constant is sharp. Thus $G_\Phi=\OLC$ when $n=2$, or for the operator norm.
\end{theorem}
\begin{proof}
The invariances, zero set and homogeneity follow immediately in log-singular coordinates. Horn's log-majorisation inequality \cite{Horn50} states
$\log s(AB)\prec\log s(A)+\log s(B)$, with the vectors ordered decreasingly; equality of the full sums follows from determinants. Subtracting $(c_1+c_2)\one$, applying gauge monotonicity under majorisation, and then the triangle inequality gives
\[
\Phi(\log s(AB)-(c_1+c_2)\one)
\le\Phi(\log s(A)-c_1\one)+\Phi(\log s(B)-c_2\one).
\]
Minimise over $c_1,c_2$.
For the comparison, put $x=\log s(A)$ and $z=(x_1-x_n)/2$. Centring at the midpoint gives the upper bound. For any $c$, the endpoint deviations $\alpha=|x_1-c|$, $\beta=|x_n-c|$ satisfy $\alpha+\beta\ge2z$. Monotonicity and symmetry imply
\[
\Phi(x-c\one)\ge\Phi(\alpha e_1+\beta e_n)
\ge\frac{\alpha+\beta}2\Phi(e_1+e_n)\ge z\Phi(e_1+e_2).
\]
The vector $(z,0,\ldots,0,-z)$, centred at zero, attains the bound.
\end{proof}
For Schatten gauges with $p\ge1$, write $G_p=\min_c(\int|x-c|^p\dd\nu_A)^{1/p}$. Then $G_p$ increases to $\OLC$: minimise the monotone $L^p$ norms over the compact interval of spectral endpoints and use uniform convergence there (or Dini's theorem). The minimising centre is a median for $p=1$, and the mean for $p=2$. Consequently
\[
G_2=\sqrt{\Var(\nu_A)},\qquad G_2(A\otimes B)^2=G_2(A)^2+G_2(B)^2.
\]
The cosh representation also gives, for every $p>0$,
\[
C_p(A)\ge\frac1{2p}\log(1+p^2G_2(A)^2),\qquad
\lim_{p\downarrow0}\frac{C_p(A)}p=\tfrac12G_2(A)^2.
\]
Indeed $\cosh u\ge1+u^2/2$ and the mean squared pair difference is twice the variance.

\section{Connections and remaining questions}
The two constructions have different tensor laws: $C_p$ is additive, while $G_2^2$ is additive. Product subadditivity and power homogeneity hold for every $G_\Phi$, but for $C_\Phi$ only the operator norm has both properties without a defect. The unresolved questions include the maximal defect at fixed individual contrasts, the exceptional homometric configurations beyond the generic reconstruction theorem, and whether the complete Ky Fan family determines singular values up to scaling and inversion. The finite-trace continuation of these questions is developed in Chapter~\ref{ch:tracial}.

\chapter{Schatten contrasts in finite von Neumann algebras}\label{ch:tracial}
\section{Introduction}

For an invertible bounded operator $A$ the \emph{operator Michelson contrast}
and the \emph{operator logarithmic contrast} are
\begin{gather*}
\OMC(A)=\frac{\kappa(A)-1}{\kappa(A)+1},\qquad
\OLC(A)=\operatorname{artanh}\OMC(A)=\tfrac12\ln\kappa(A),\\
\kappa(A)=\|A\|\,\|A^{-1}\|,
\end{gather*}
see \cite{ANN24} and Nikolaeva and Novikov \cite{NN26}. Replacing the operator norm on $M_n(\C)$ by the
normalised Schatten norms
$\|X\|_p=\tau(|X|^p)^{1/p}$, $\tau=\frac1n\Tr$, one obtains the contrast
\[
C_p(A)=\tfrac12\ln\bigl(\|A\|_p\,\|A^{-1}\|_p\bigr),
\]
which is nonnegative, vanishes exactly on scalar multiples of unitaries, increases in
$p$ to $\OLC(A)$ (Proposition~\ref{tracial:prop:vn}), and is subadditive under products up to the sharp defect
$\frac{\ln n}p$ (Theorem~\ref{tracial:thm:vn-defect} below), while
$\OLC(A)-C_p(A)=O(\frac1p)$ as $p\to\infty$ (Theorem~\ref{tracial:thm:edge}). The defect $\frac{\ln n}p$ is unbounded in the dimension, and the definition
uses the trace in an essential way. The present paper asks what remains of
the Schatten contrast in an infinite-dimensional algebra, and what it detects
there. Three questions organise the chapter.

\begin{enumerate}[label=\textup{(Q\arabic*)}]
\item \emph{Which functionals give a contrast?} If the trace is replaced by a
positive functional $\varphi$, when is the quantity
$\frac1{2p}\ln\bigl(\varphi(|A|^p)\varphi(|A^{-1}|^p)\bigr)$ nonnegative?
\item \emph{What do the constants of the matrix case become?}
What replaces the defect $\frac{\ln n}p$ and the rate $\OLC-C_p\sim\frac cp$ in
a von Neumann algebra?
\item \emph{What happens beyond finite traces}: for infinite tensor products
of matrix algebras, which may be of type $\mathrm{III}$ and then carry no
trace?
\end{enumerate}
\subsection*{Background}

\emph{Characterisations of the trace.} Gardner \cite{Gardner79} proved that the
inequality $|\varphi(x)|\le\varphi(|x|)$ characterises traces among positive
functionals on a $C^*$-algebra. Since then many classical trace inequalities,
among them those of Young, H\"older, Cauchy--Schwarz, Golden--Thompson,
Peierls--Bogoliubov and Araki--Lieb--Thirring, have been shown to hold for a
positive functional only if it is tracial; see \cite{Bik11,Bik12,MB25}. In a
parallel line, central elements of operator algebras, rather than tracial
functionals, are characterised by inequalities \cite{NT15}. The
Michelson contrast of operators arose in exactly this context, as a
quantitative measure of centrality \cite{ANN24}; its logarithmic form $\OLC$
and the associated geometry were studied in Nikolaeva and Novikov \cite{NN26}. Passing to Schatten
norms puts a trace into the definition of the contrast itself, and
question (Q1) asks whether the trace property is forced. Theorem~\ref{tracial:thm:trace}
shows that it is.

\emph{Noncommutative integration.} The quantities $\|X\|_p=\tau(|X|^p)^{1/p}$
are the norms of the noncommutative $L^p$-spaces $L^p(\cM,\tau)$ of a
semifinite von Neumann algebra \cite{Nelson74,FK86,PX03}, and
$e^{2C_p(A)}=\|A\|_p\|A^{-1}\|_p$ is the condition number of $A$ in
$L^p(\cM,\tau)$. Banach spaces of $L_1$-type associated with a single positive
operator affiliated with a von Neumann algebra were studied in
\cite{Nov17,NT16}, and invertibility of a positive element is characterised by
the equivalence of the associated norms with the original ones in \cite{Nov19};
the contrast is a quantitative form of such an equivalence. For algebras
without a trace the corresponding spaces are
those of Haagerup \cite{Haagerup79}; we do not treat them here and return to
this point in the concluding remarks.

\emph{Dimension and types.} Murray and von Neumann \cite{MvN36} classified
finite factors by the range of the trace on projections,
$\{0,\frac1n,\dots,1\}$ in type $\mathrm I_n$ and $[0,1]$ in type
$\mathrm{II}_1$. Theorems~\ref{tracial:thm:vn-defect} and \ref{tracial:thm:edge} show that the
two basic constants of the matrix case, the subadditivity defect
and the first-order correction to $\OLC$, are functions of this range. For
factors of type $\mathrm{III}$, the work of Powers \cite{Powers67}, Araki and
Woods \cite{AW68} and Connes \cite{Connes73} showed that the type of an
infinite tensor product of matrix algebras is governed by the ratios of the
eigenvalues of the component states, that is, for $2\times2$ components, by
their condition numbers.
Section~\ref{tracial:sec:itpfi} makes this dictionary explicit.

\subsection*{Main results}

\noindent\textbf{(A) Positivity of the contrast characterises traces}
(Theorem~\ref{tracial:thm:trace}). For a nonzero positive functional $\varphi$ on a
unital $C^*$-algebra and $p\in(0,\infty)$, the inequality
$\varphi(|A|^p)\varphi(|A^{-1}|^p)\ge\varphi(1)^2$ for all invertible $A$ holds
if and only if $\varphi$ is tracial; the same is true for the symmetry
$A\mapsto A^*$ of the contrast.

\medskip\noindent\textbf{(B) Monotonicity} (Proposition~\ref{tracial:prop:mono}). If
$\Psi$ is a positive unital trace-preserving map between tracial von Neumann
algebras and $A$ is positive and invertible, then $C_p(\Psi(A))\le C_p(A)$ for
$1\le p<\infty$, and this fails for $p<1$. In particular the contrast decreases
under conditional expectations and under partial traces, and for $p\ge1$
every contrast in the hyperfinite factor of type $\mathrm{II}_1$ is a monotone
limit of matrix contrasts (Example~\ref{tracial:ex:R}).

\medskip\noindent\textbf{(C) The trace must be finite; the defect is a
dimension} (Proposition~\ref{tracial:prop:BH}, Theorem~\ref{tracial:thm:vn-defect},
Corollary~\ref{tracial:cor:types}). If $\tau(1)=\infty$, then
$\|A\|_p\|A^{-1}\|_p=\infty$ for every invertible $\tau$-measurable $A$. For a
faithful normal tracial state,
\[
C_p(AB)\le C_p(A)+C_p(B)+\frac1p\ln\frac1{\varepsilon_\cM},
\qquad\varepsilon_\cM=\inf\{\tau(e):\ e\ne0\ \text{a projection}\},
\]
and for factors the constant is optimal: $\frac{\ln n}p$ in type $\mathrm I_n$
and $+\infty$ in type $\mathrm{II}_1$. Thus the defect is $\frac1p$ times the
logarithm of the reciprocal of the smallest Murray--von Neumann dimension of a
nonzero projection, and the types $\mathrm I_n$, $\mathrm{II}_1$,
$\mathrm I_\infty$, $\mathrm{II}_\infty$ of semifinite factors are
distinguished by the behaviour of $C_p$.

\medskip\noindent\textbf{(D) The rate of convergence to $\OLC$}
(Theorem~\ref{tracial:thm:edge}, Proposition~\ref{tracial:prop:logcorr}). For invertible $A$
with $|A|$ not scalar,
\[
\lim_{p\to\infty}2p\bigl(\OLC(A)-C_p(A)\bigr)=-\ln\tau(e_+)-\ln\tau(e_-)\in[2\ln2,\infty],
\]
where $e_\pm$ are the spectral projections of $|A|$ at the top and the bottom
of its spectrum. Hence $\OLC-C_p=O(1/p)$ for all invertible elements if and
only if $\cM$ is finite-dimensional; in $M_n(\C)$ the limits form a finite set,
in a factor of type $\mathrm{II}_1$ they fill $[2\ln2,\infty]$, and when
$e_\pm=0$ and the power-law edge hypotheses of Proposition~\ref{tracial:prop:logcorr} hold, the correction is of order $\frac{\ln p}p$ with a coefficient given by
the edge exponents of the spectral measure: $\frac32$ for the exponential of a
semicircular element, $\frac12$ for Haar unitaries
(Examples~\ref{tracial:ex:semicircular}, \ref{tracial:ex:arcsine}). The second derivative of
the profile, $L_A''$, decays exponentially exactly when the top of the spectrum
is isolated, at the rate given by the gap; this last fact is classical and is
recalled in Proposition~\ref{tracial:prop:edge-exp}.

\medskip\noindent\textbf{(E) Araki--Woods factors}
(Section~\ref{tracial:sec:itpfi}). For an infinite tensor product of $2\times2$ density
matrices $\rho_n$ we reformulate the classical criteria of Araki--Woods and
Connes: types $\mathrm{II}_1$ and $\mathrm I$ correspond to summable saturation
of the two lower bounds $C_1\ge0$ and $C_1\ge\OLC-\ln2$, type $\mathrm{III}$ to
the divergence of $\sum_n\min\{a_n,b_n\}$, where $a_n,b_n$ are the two
saturation gaps, type $\mathrm{II}_\infty$ to the remaining case, and Connes'
invariant is
$T=\{t:\sum_n(1-\OMC(\rho_n)^2)\sin^2(t\,\OLC(\rho_n))<\infty\}$
(Proposition~\ref{tracial:prop:itpfi}; this is a reformulation, not a new theorem). For
components of arbitrary size with uniformly bounded contrast we prove that the
factor is of type $\mathrm{II}_1$ if and only if the Schatten contrasts of the
partial products $\rho_1\otimes\dots\otimes\rho_N$ stay bounded
(Theorem~\ref{tracial:thm:itpfi-bounded}). Without the bound, boundedness of $C_p$ for
some $p\ge1$ still implies type $\mathrm{II}_1$ (Proposition~\ref{tracial:prop:itpfi-p1}),
but explicit examples show that this fails for $p<1$ and that a factor of type
$\mathrm{II}_1$ may have unbounded $C_p$ for all $p\ge2$, or for all $p>0$
(Example~\ref{tracial:ex:itpfi-pc}).

\subsection*{Relation to classical results}
The trace characterisation, the factor defect and the finite-dimensionality criterion are proved below. Monotonicity follows from established Jensen trace inequalities, and the spectral-edge estimates use classical Laplace asymptotics. The two-dimensional Araki--Woods criteria are explicitly presented as reformulations of the classical classification; the bounded-component contrast criterion is obtained by uniform comparison with its square-root criterion. No classification of factors beyond those cited results is claimed.

\subsection*{Scope}
Unbounded affiliated operators, for which $\OLC=\infty$, and finite
sections of unbounded operators on a Hilbert space are not treated here; the
definitions of Section~\ref{tracial:subsec:def} are, however, stated for them, to fix the notation.

\section{Schatten contrasts in a tracial von Neumann algebra}\label{tracial:sec:vn}

Throughout, $\cM$ is a von Neumann algebra with a faithful normal semifinite
trace $\tau$, and for a $\tau$-measurable operator $X$ we put
$\|X\|_p=\tau(|X|^p)^{1/p}\in[0,\infty]$, $p\in(0,\infty)$ \cite{Nelson74,FK86}.
We first show that the trace has to be finite, and then set up the contrast for
a tracial state.

\subsection{The Schatten contrast requires a finite trace}

\begin{proposition}\label{tracial:prop:BH}
Let $\cM$ be a von Neumann algebra with a faithful normal semifinite trace
$\tau$ such that $\tau(1)=\infty$, and let $A$ be a $\tau$-measurable operator
with $\tau$-measurable inverse. Then
$\|A\|_p\,\|A^{-1}\|_p=\infty$ for every $p\in(0,\infty)$, where
$\|X\|_p=\tau(|X|^p)^{1/p}\in[0,\infty]$. In particular this holds in
$B(\cH)$, $\dim\cH=\infty$, with the usual trace and the Schatten classes.
\end{proposition}

\begin{proof}
In the polar decomposition $A=U|A|$ the partial isometry $U$ is unitary because
$A$ is invertible, and $|A^{-1}|=U|A|^{-1}U^*$. Since the generalised singular
values are unitarily invariant, $\mu(UXU^*)=\mu(X)$ \cite{FK86}, we have
$\tau(|A^{-1}|^p)=\tau(|A|^{-p})$. As $t^p+t^{-p}\ge2$ for $t>0$, the
functional calculus gives $|A|^p+|A|^{-p}\ge2\cdot1$, hence
$\tau(|A|^p)+\tau(|A^{-1}|^p)=\tau(|A|^p+|A|^{-p})\ge2\tau(1)=\infty$, so at
least one of the two factors is infinite, while the other is nonzero because
$\tau$ is faithful.
\end{proof}

So the Schatten contrasts require a \emph{finite} trace. This is consistent
with the matrix case, in which the defect $\frac{\ln n}p$ is
unbounded in $n$ (Theorem~\ref{tracial:thm:vn-defect}).

\subsection{Contrast in a finite von Neumann algebra}\label{tracial:subsec:def}

Let $\cM$ be a von Neumann algebra with a faithful normal tracial state $\tau$,
and $\widetilde\cM$ the $*$-algebra of closed densely defined operators
affiliated with $\cM$ (all of them are $\tau$-measurable since $\tau$ is finite)
\cite{Nelson74,FK86}. If $A\in\widetilde\cM$ is injective, then its right
support is $1$; the left support is equivalent to it, hence equal to $1$ since
$\cM$ is finite, so $A$ has dense range and $A^{-1}\in\widetilde\cM$; moreover the polar decomposition $A=U|A|$ has $U$ unitary, and
$|A^{-1}|=U|A|^{-1}U^*$. The \emph{logarithmic spectral measure} of $A$ is the
Borel probability measure
\[
\nu_A(E)=\tau\bigl(\one_E(\ln|A|)\bigr),\qquad E\subset\R ,
\]
so that $\|A\|_p^p=\tau(|A|^p)=\int e^{px}\dd\nu_A$ and
$\|A^{-1}\|_p^p=\tau(|A|^{-p})=\int e^{-px}\dd\nu_A$, with values in
$(0,\infty]$. Every Borel probability measure on $\R$ arises in this way, already
for $\cM=L^\infty[0,1]$ with the Lebesgue integral as trace and $A$ the
multiplication operator by $\exp(Q_\nu)$, $Q_\nu$ being the quantile function of
$\nu$.

\begin{definition}
For a Borel probability measure $\nu$ on $\R$ and $p\in(0,\infty)$ put
\begin{gather*}
C_p(\nu)=\frac1{2p}\Bigl(\ln\int e^{px}\dd\nu(x)+\ln\int e^{-px}\dd\nu(x)\Bigr)\in[0,\infty],\\
C_\infty(\nu)=\tfrac12\operatorname{diam}\supp\nu\in[0,\infty],
\end{gather*}
and $C_p(A)=C_p(\nu_A)=\frac12\ln\bigl(\|A\|_p\|A^{-1}\|_p\bigr)$,
$\OLC(A)=C_\infty(\nu_A)$ for injective $A\in\widetilde\cM$. The \emph{critical
exponent} is $p_c(\nu)=\sup\{p>0:\ C_p(\nu)<\infty\}\in[0,\infty]$ (with
$\sup\emptyset=0$).
\end{definition}

For $\cM=M_n(\C)$, $\tau=\frac1n\Tr$, this is the Schatten contrast of the
introduction. We write $L_\nu(t)=\ln\int e^{tx}\dd\nu(x)\in(-\infty,\infty]$ for
the cumulant generating function and $L_A=L_{\nu_A}$, so that
$2p\,C_p(A)=L_A(p)+L_A(-p)$, and we recall that $\OMC=\tanh\OLC$. When $\nu_A$
has a first moment, $2p\,C_p(A)=K_A(p)+K_A(-p)$, where
$K_A(t)=L_A(t)-t\int x\,\dd\nu_A$ is the centred cumulant generating function.
We call $L_A$ the \emph{logarithmic spectral profile} of $A$;
the quantity $pC_p(A)=\tfrac12(K_A(p)+K_A(-p))$ is the even part of the centred profile, which is why the
edge exponents and the masses of the spectral measure at its edges appear in
Section~\ref{tracial:sec:rate}. We recall the classical facts that $L_A$ is convex and
that, in the interior of the set where $L_A$ is finite, $L_A'(t)=\int x\,\dd\nu_A^{(t)}$ and
$L_A''(t)=\Var(\nu_A^{(t)})$ for the tilted measure
$\dd\nu_A^{(t)}=e^{tx-L_A(t)}\dd\nu_A$.

\begin{proposition}\label{tracial:prop:vn}
Let $\nu$ be a Borel probability measure on $\R$ and $0<p<\infty$.
\begin{enumerate}
\item $C_p(\nu)\ge0$, with equality if and only if $\nu$ is a point mass.
\item $C_p(\nu)<\infty$ if and only if $\int e^{p|x|}\dd\nu<\infty$.
\item $p\mapsto C_p(\nu)$ is nondecreasing and
$\lim_{p\to\infty}C_p(\nu)=C_\infty(\nu)$.
\item $C_p(\nu)\le\frac1p\ln\cosh\bigl(p\,C_\infty(\nu)\bigr)$ whenever
$C_\infty(\nu)<\infty$.
\item If $\nu_s$ is the image of $\nu$ under $x\mapsto sx$, $s\ne0$, then
$C_p(\nu_s)=|s|\,C_{|s|p}(\nu)$; and $C_p(\nu*\nu')=C_p(\nu)+C_p(\nu')$.
\end{enumerate}
\end{proposition}

\begin{proof}
(a) By the Cauchy--Schwarz inequality,
\[1=\bigl(\int e^{px/2}e^{-px/2}\dd\nu\bigr)^2\le\int e^{px}\dd\nu\int e^{-px}\dd\nu,\]
with equality iff $e^{px}$ is $\nu$-a.e.\ constant.
(b) Both integrals are positive, so their product is finite iff both are, and
$e^{p|x|}\le e^{px}+e^{-px}\le2e^{p|x|}$.
(c) On a probability space the $L^p$-norms of the functions $e^{x}$ and $e^{-x}$
are nondecreasing in $p$ and converge to their essential suprema.
(d) The case $C_\infty(\nu)=0$ follows from (a). Let $z=C_\infty(\nu)>0$. Since $C_p$ is invariant under translations of
$\nu$, we may assume $\supp\nu\subset[-z,z]$. By convexity,
$e^{\pm py}\le\cosh(pz)\pm\frac yz\sinh(pz)$ for $y\in[-z,z]$; integrating, with
$\mu=\int y\,\dd\nu$,
\[
\int e^{py}\dd\nu\int e^{-py}\dd\nu\le\cosh^2(pz)-\frac{\mu^2}{z^2}\sinh^2(pz)\le\cosh^2(pz).
\]
(e) $\int e^{\pm px}\dd\nu_s=\int e^{\pm psx}\dd\nu$, and
$\int e^{tx}\dd(\nu*\nu')=\int e^{tx}\dd\nu\int e^{tx}\dd\nu'$.
\end{proof}

\begin{remark}\label{tracial:rem:nosub}
For injective elements of $\widetilde\cM$ there is no subadditivity estimate at a fixed exponent: if
$p<p_c(\nu_A)<2p$ for a positive $A$, then $C_p(A)<\infty$ but
$C_p(A^2)=2C_{2p}(A)=\infty$. For bounded invertible elements the failure of subadditivity is quantified
in Theorem~\ref{tracial:thm:vn-defect}.
\end{remark}

\subsection{Monotonicity under trace-preserving unital positive maps}\label{tracial:subsec:mono}

The contrast of positive elements decreases under coarse-graining. This is a
consequence of a form of Jensen's trace inequality (cf.\ \cite{Petz87,BK90});
since the function $t^{-p}$ is not defined at $0$ and $\Psi$ need not be a
self-map, we include a short proof.

\begin{proposition}\label{tracial:prop:mono}
Let $(\cM,\tau)$ and $(\mathcal N,\tau')$ be von Neumann algebras with faithful
normal tracial states, and let $\Psi\colon\cM\to\mathcal N$ be a positive,
unital linear map with $\tau'\circ\Psi=\tau$. Let $A\in\cM$ be positive and
invertible. Then $\Psi(A)$ is positive and invertible,
\[
\OLC(\Psi(A))\le\OLC(A),\qquad C_p(\Psi(A))\le C_p(A)\quad\text{for }1\le p<\infty .
\]
In particular this holds for $\tau$-preserving conditional expectations onto
von Neumann subalgebras.
\end{proposition}

\begin{proof}
Let $m\le A\le M$ with $0<m\le M$ the extreme points of the spectrum of $A$.
Since $\Psi$ is positive and unital, $m\le\Psi(A)\le M$, which gives the first
two statements. We claim that for every continuous convex function $f$ on
$[m,M]$,
\begin{equation}\label{tracial:eq:jensen}
\tau'\bigl(f(\Psi(A))\bigr)\le\tau\bigl(f(A)\bigr).
\end{equation}
Applied to $f(t)=t^p$, $p\ge1$, and to $f(t)=t^{-p}$ this gives
$\tau'(\Psi(A)^{\pm p})\le\tau(A^{\pm p})$, hence the claim for $C_p$.

To prove \eqref{tracial:eq:jensen} put $B=\Psi(A)$, let $\mathcal B\subset\mathcal N$ be
the abelian von Neumann algebra generated by $B$, and let $E$ be the
$\tau'$-preserving conditional expectation of $\mathcal N$ onto $\mathcal B$.
The map $\Theta(g)=E\bigl(\Psi(g(A))\bigr)$, $g\in C[m,M]$, is positive and
unital with values in $\mathcal B$, and $\Theta(\mathrm{id})=E(B)=B$. If
$\ell(t)=\alpha t+\beta$ is an affine function with $\ell\le f$ on $[m,M]$, then
$\Theta(f)\ge\Theta(\ell)=\alpha B+\beta=\ell(B)$. As $f$ is the supremum of
countably many such $\ell$, and $\mathcal B$ is an algebra of functions, we get
$\Theta(f)\ge f(B)$. Therefore
$\tau'(f(B))\le\tau'(\Theta(f))=\tau'(\Psi(f(A)))=\tau(f(A))$.
\end{proof}

\begin{remark}\label{tracial:rem:mono}
(i) The restriction $p\ge1$ is necessary. In $M_3(\C)$ let
$\Psi(X)=\frac12(X+UXU^*)$, where $U$ is the permutation matrix exchanging the
first two basis vectors, and $A=\operatorname{diag}(16,1,\frac1{16})$. Then
$\Psi(A)=\operatorname{diag}(\frac{17}2,\frac{17}2,\frac1{16})$ and, for
$p=\frac12$, $e^{C_{1/2}(A)}=\tau(A^{1/2})\tau(A^{-1/2})=\frac{49}{16}=3.0625$,
whereas $e^{C_{1/2}(\Psi(A))}=3.1661\ldots$ More generally, for
$A=\operatorname{diag}(s,1,s^{-1})$ one has
$\Psi(A)=\operatorname{diag}(\frac{s+1}2,\frac{s+1}2,s^{-1})$ and
\[
\frac{\tau\bigl(\Psi(A)^p\bigr)\,\tau\bigl(\Psi(A)^{-p}\bigr)}{\tau(A^p)\,\tau(A^{-p})}
\longrightarrow2^{1-p}\qquad(s\to\infty),
\]
so that $C_p(\Psi(A))>C_p(A)$ for all large $s$, for every $p\in(0,1)$.

(ii) Unitality is necessary: for a non-scalar positive $S$ with $\tau(S)=1$,
the map $X\mapsto\tau(X)S$ is positive and trace preserving and sends $1$, of
contrast $0$, to $S$.

(iii) For matrices, Proposition~\ref{tracial:prop:mono} applies to pinchings, to
mixtures of unitary conjugations and to the normalised partial trace
$d_2^{-1}\Tr_2\colon M_{d_1}\otimes M_{d_2}\to M_{d_1}$. Since $C_p$ is invariant
under scaling, for a positive definite density matrix $\rho$ on
$\C^{d_1}\otimes\C^{d_2}$ this gives
\[
C_p(\Tr_2\rho)\le C_p(\rho),\qquad\OLC(\Tr_2\rho)\le\OLC(\rho)\qquad(1\le p<\infty)
\]
(the second inequality also follows from Chapter~\ref{ch:core}):
the contrast of a marginal does not exceed the contrast of the joint state.
We do not pursue quantum-information interpretations here; in particular we
make no claim that the defect of Section~\ref{tracial:sec:defect-vn} is a measure of
entanglement.
\end{remark}

\section{Positivity of the contrast characterises traces}\label{tracial:sec:trace}

In \cite{ANN24} the Michelson contrast was related to inequalities that
characterise tracial functionals, in the tradition of Gardner's theorem
\cite{Gardner79}. The Schatten contrast gives a characterisation of the same
kind. Let $\mathcal A$ be a unital $C^*$-algebra and $\varphi$ a positive linear
functional on $\mathcal A$. For invertible $A\in\mathcal A$ and $p\in(0,\infty)$
put
\[
\kappa_p^\varphi(A)=\bigl(\varphi(|A|^p)\,\varphi(|A^{-1}|^p)\bigr)^{1/p}.
\]
If $\varphi=\tau$ is a tracial state this is $e^{2C_p(A)}$.

\begin{theorem}\label{tracial:thm:trace}
Let $\varphi\ne0$ be a positive linear functional on a unital $C^*$-algebra
$\mathcal A$ and $p\in(0,\infty)$. The following are equivalent:
\begin{enumerate}
\item $\varphi$ is tracial, i.e.\ $\varphi(XY)=\varphi(YX)$ for all
$X,Y\in\mathcal A$;
\item $\kappa_p^\varphi(A)\ge\kappa_p^\varphi(1)$, i.e.\
$\varphi(|A|^p)\,\varphi(|A^{-1}|^p)\ge\varphi(1)^2$, for all invertible
$A\in\mathcal A$;
\item $\kappa_p^\varphi(A^*)=\kappa_p^\varphi(A)$ for all invertible
$A\in\mathcal A$.
\end{enumerate}
\end{theorem}

\begin{proof}
Let $A$ be invertible with polar decomposition $A=UP$, $P=|A|$, $U=AP^{-1}$
unitary. Then
\begin{equation}\label{tracial:eq:moduli}
|A|=P,\qquad|A^{-1}|=UP^{-1}U^*,\qquad|A^*|=UPU^*,\qquad|(A^*)^{-1}|=P^{-1}.
\end{equation}

(a)$\Rightarrow$(b),(c). If $\varphi$ is tracial, then by \eqref{tracial:eq:moduli}
$\kappa_p^\varphi(A)^p=\varphi(P^p)\varphi(P^{-p})=\kappa_p^\varphi(A^*)^p$,
and by the Cauchy--Schwarz inequality for $\varphi$,
$\varphi(1)^2=\varphi(P^{p/2}P^{-p/2})^2\le\varphi(P^p)\,\varphi(P^{-p})$.

(b)$\Rightarrow$(a). Fix $X=X^*\in\mathcal A$ and a unitary $U\in\mathcal A$.
For real $\varepsilon$ with $|\varepsilon|\,\|X\|<1$ let
$P_\varepsilon=(1+\varepsilon X)^{1/p}$ and $A_\varepsilon=UP_\varepsilon$. By
\eqref{tracial:eq:moduli},
\begin{align*}
f(\varepsilon):=\kappa_p^\varphi(A_\varepsilon)^p
&=\varphi(1+\varepsilon X)\,\varphi\bigl(U(1+\varepsilon X)^{-1}U^*\bigr)\\
&=\varphi(1)^2+\varepsilon\,\varphi(1)\bigl(\varphi(X)-\varphi(UXU^*)\bigr)+O(\varepsilon^2),
\end{align*}
because $(1+\varepsilon X)^{-1}=1-\varepsilon X+O(\varepsilon^2)$ in norm.
By (b), $f$ has a minimum at $\varepsilon=0$, so the linear term vanishes;
since $\varphi(1)>0$, $\varphi(UXU^*)=\varphi(X)$. By linearity this holds for
all $X\in\mathcal A$; replacing $X$ by $XU$ gives $\varphi(UX)=\varphi(XU)$, and
since the unitaries span $\mathcal A$, $\varphi$ is tracial.

(c)$\Rightarrow$(a). With the same notation, by \eqref{tracial:eq:moduli},
\[
0=\kappa_p^\varphi(A_\varepsilon)^p-\kappa_p^\varphi(A_\varepsilon^*)^p
=2\varepsilon\,\varphi(1)\bigl(\varphi(X)-\varphi(UXU^*)\bigr)+O(\varepsilon^2),
\]
and we conclude as before.
\end{proof}

Thus, for a non-tracial state, the normalised contrast
$\frac1{2p}\ln\bigl(\varphi(|A|^p)\varphi(|A^{-1}|^p)\bigr)$ takes negative
values and distinguishes $A$ from $A^*$: nonnegativity of the contrast and the symmetry
$C_p(A)=C_p(A^*)$ are each equivalent to the trace property.

\begin{remark}\label{tracial:rem:trace-lit}
Many trace inequalities (Cauchy--Schwarz, H\"older, Young, Golden--Thompson,
Peierls--Bogoliubov, inequalities for the modulus and others) are known to
characterise tracial functionals among positive ones, often already when
restricted to projections; see
\cite{Gardner79,PZ88,STS02,SY06,Bik11,Bik12,MB25,BM26} and the references in
\cite{ANN24}. Condition (b) involves the two operators $|A|^p$ and
$|A^{-1}|^p=U|A|^{-p}U^*$, which do not commute in general; only for a tracial
$\varphi$ does it reduce to the Cauchy--Schwarz inequality used in the proof.
We have not found (b) or (c) in this product form in the literature. The proof
shows that it suffices to require (b) or (c) for invertible elements in an
arbitrarily small norm neighbourhood of the unitary group.
\end{remark}

\section{The defect of a finite von Neumann algebra and the type of a factor}\label{tracial:sec:defect-vn}

Let $(\cM,\tau)$ be a von Neumann algebra with a faithful normal tracial state,
$\GL(\cM)$ its group of invertible elements, and
\begin{align*}
\varepsilon_\cM&=\inf\{\tau(e):\ e\in\cM\ \text{a nonzero projection}\}\in[0,1],\\
D_p(\cM)&=\sup\bigl\{C_p(P^2)-2C_p(P):\ P\in\GL(\cM),\ P>0\bigr\}\in[0,\infty].
\end{align*}
For $\cM=M_n(\C)$, $\varepsilon_\cM=\frac1n$. More generally, for a
finite-dimensional algebra $\cM=\bigoplus_kM_{n_k}(\C)$ with the tracial state
$\tau=\sum_kw_k\,\frac1{n_k}\Tr$, where $w_k>0$ and $\sum_kw_k=1$, one has
$\varepsilon_\cM=\min_kw_k/n_k$.

\begin{theorem}\label{tracial:thm:vn-defect}
Let $(\cM,\tau)$ be as above and $p\in(0,\infty)$.
\begin{enumerate}
\item\label{tracial:vd:sub} For all $A,B\in\GL(\cM)$,
$\displaystyle C_p(AB)\le C_p(A)+C_p(B)+\frac1p\ln\frac1{\varepsilon_\cM}$
(with no content if $\varepsilon_\cM=0$).
\item\label{tracial:vd:sharp} If $\cM\ne\C$ is a factor, then
$D_p(\cM)=\frac1p\ln\frac1{\varepsilon_\cM}$. Thus $D_p(\cM)=\frac{\ln n}p$ if
$\cM$ is of type $\mathrm{I}_n$, and $D_p(\cM)=+\infty$ if $\cM$ is of type
$\mathrm{II}_1$.
\end{enumerate}
\end{theorem}

\begin{proof}
\ref{tracial:vd:sub} We may assume $\varepsilon_\cM>0$. Let $X\in\cM$, $X\ne0$, and
$0<\eta<\|X\|$. The spectral projection
$e_\eta=\one_{(\|X\|-\eta,\|X\|]}(|X|)$ is nonzero and
$|X|^p\ge(\|X\|-\eta)^pe_\eta$, so
$\tau(|X|^p)\ge(\|X\|-\eta)^p\varepsilon_\cM$. Hence
$\|X\|\le\varepsilon_\cM^{-1/p}\|X\|_p$. Together with
$\|XY\|_p\le\|X\|\,\|Y\|_p$ \cite{FK86} this gives
$\|AB\|_p\le\varepsilon_\cM^{-1/p}\|A\|_p\|B\|_p$, and the same for
$(AB)^{-1}=B^{-1}A^{-1}$. Multiply and take half the logarithm.

\ref{tracial:vd:sharp} The upper bound is \ref{tracial:vd:sub} with $A=B=P$. For the lower
bound let $\delta\in(0,\frac12]$ be such that there are orthogonal projections
$e,f\in\cM$ with $\tau(e)=\tau(f)=\delta$, and put
$P_t=te+(1-e-f)+t^{-1}f$, $t>1$. For every $q>0$,
$\tau(P_t^{q})=\tau(P_t^{-q})=\delta t^q+(1-2\delta)+\delta t^{-q}$, so
$C_q(P_t)=\frac1q\ln\bigl(\delta t^q+1-2\delta+\delta t^{-q}\bigr)=\ln t+\frac1q\ln\delta+o(1)$
as $t\to\infty$. By the power law of Proposition~\ref{tracial:prop:vn},
\[
C_p(P_t^2)-2C_p(P_t)=2C_{2p}(P_t)-2C_p(P_t)\longrightarrow\frac1p\ln\delta-\frac2p\ln\delta=\frac1p\ln\frac1\delta .
\]
In the factor of type $\mathrm{I}_n$, $n\ge2$, one can take
$\delta=\frac1n=\varepsilon_\cM$; in a factor of type $\mathrm{II}_1$ every
$\delta\in(0,\frac12]$ is admissible and $\varepsilon_\cM=0$.
\end{proof}

\begin{corollary}[The type of a semifinite factor]\label{tracial:cor:types}
Let $\cM$ be a factor with a faithful normal semifinite trace $\tau$,
normalised by $\tau(1)=1$ if $\tau$ is finite, and let $p\in(0,\infty)$.
For $\cM=\C$ we have $D_p(\cM)=0$, and (a) holds with $n=1$.
\begin{enumerate}
\item $\cM$ is of type $\mathrm{I}_n$, $n<\infty$, if and only if
$\|A\|_p\|A^{-1}\|_p<\infty$ on $\GL(\cM)$ and $D_p(\cM)<\infty$; in this case
$n=\exp\bigl(p\,D_p(\cM)\bigr)$.
\item $\cM$ is of type $\mathrm{II}_1$ if and only if
$\|A\|_p\|A^{-1}\|_p<\infty$ on $\GL(\cM)$ and $D_p(\cM)=\infty$.
\item $\cM$ is of type $\mathrm{I}_\infty$ or $\mathrm{II}_\infty$ if and only if
$\|A\|_p\|A^{-1}\|_p=\infty$ for all $A\in\GL(\cM)$. The two types are
distinguished by their corners: for a nonzero projection $e$ with
$\tau(e)<\infty$, the finite factor $e\cM e$ with the tracial state
$\tau/\tau(e)$ has $D_p(e\cM e)=\frac1p\ln\operatorname{rank}(e)<\infty$ in type
$\mathrm{I}_\infty$ and $D_p(e\cM e)=\infty$ in type $\mathrm{II}_\infty$.
\end{enumerate}
Factors of type $\mathrm{III}$ carry no normal semifinite trace, and the tracial Schatten contrast used here is unavailable; the operator-norm contrast $\OLC$ remains defined on bounded invertible elements.
\end{corollary}

\begin{proof}
If $\tau(1)=1$ and $A\in\GL(\cM)$, then
$\|A\|_p\|A^{-1}\|_p\le\|A\|\,\|A^{-1}\|<\infty$; if $\tau(1)=\infty$,
Proposition~\ref{tracial:prop:BH} applies. The rest follows from
Theorem~\ref{tracial:thm:vn-defect} and the classification of semifinite factors
\cite{Takesaki79}; in type $\mathrm{I}_\infty$ the corner $e\cM e$ is isomorphic
to $M_r(\C)$, $r=\operatorname{rank}(e)$, and in type $\mathrm{II}_\infty$ it is
of type $\mathrm{II}_1$.
\end{proof}

In other words, for a finite factor the sharp defect is $\frac1p$ times the
logarithm of the reciprocal of the smallest Murray--von Neumann dimension
$\tau(e)$ of a nonzero projection $e$; for matrices it is the constant
$\frac{\ln n}p$.

\begin{example}[The hyperfinite factor]\label{tracial:ex:R}
Let $R=\bigotimes_{k\ge1}(M_2(\C),\frac12\Tr)$ be the hyperfinite factor of type
$\mathrm{II}_1$, with its increasing sequence of subfactors
$R_k\cong M_{2^k}(\C)$, whose union is dense, and the trace-preserving
conditional expectations $E_k\colon R\to R_k$. By Theorem~\ref{tracial:thm:vn-defect},
\[
D_p(R_k)=\frac{k\ln2}{p}\ \longrightarrow\ \infty=D_p(R)\qquad(k\to\infty):
\]
the defect grows linearly in $k$, that is, logarithmically in the matrix size,
and becomes infinite in the limit. The contrast itself, on the other hand, is
continuous along this approximation: for positive invertible $A\in R$ and
$1\le p<\infty$,
\[
C_p\bigl(E_k(A)\bigr)\ \uparrow\ C_p(A)\qquad(k\to\infty).
\]
Indeed, $E_k=E_k\circ E_{k+1}$, so $C_p(E_k(A))\le C_p(E_{k+1}(A))\le C_p(A)$ by
Proposition~\ref{tracial:prop:mono}. Moreover $E_k(A)\to A$ in $\|\cdot\|_2$, all
$E_k(A)$ have spectrum in the interval $[m,M]$ spanned by the spectrum of $A$,
and for uniformly bounded self-adjoint $X_k\to X$ in $\|\cdot\|_2$ one has
$\tau(g(X_k))\to\tau(g(X))$ for polynomials $g$, hence for continuous $g$ on
$[m,M]$; with $g(t)=t^{\pm p}$ this gives $C_p(E_k(A))\to C_p(A)$. Since
$C_p(A)=C_p(|A|)$ for invertible $A$, every Schatten contrast $C_p$, $p\ge1$, in
$R$ is a monotone limit of matrix contrasts, although no subadditivity estimate
survives the limit.
\end{example}

\begin{remark}[Symmetric operator spaces]\label{tracial:rem:symmetric}
The same argument works for a symmetric Banach function space $E$ on $(0,1)$
and the associated noncommutative space $E(\cM,\tau)$ over a factor of type
$\mathrm{II}_1$ \cite{DDP89}. Recall that $L^\infty(0,1)\subset E$ with
$\|x\|_E\le\|\one\|_E\|x\|_\infty$. Let $\phi_E(t)=\|e\|_E$ for a projection
$e$ with $\tau(e)=t$ be the fundamental function, which is nondecreasing, and
$C_E(A)=\frac12\ln\bigl(\|A\|_E\|A^{-1}\|_E/\phi_E(1)^2\bigr)$. If $X\ne0$ and
$0<\eta<\|X\|$, the spectral projection $e_\eta$ of $|X|$ for
$(\|X\|-\eta,\|X\|]$ is nonzero and $|X|\ge(\|X\|-\eta)e_\eta$, so
$\|X\|_E\ge(\|X\|-\eta)\phi_E(\tau(e_\eta))\ge(\|X\|-\eta)\phi_E(0+)$. Hence
$\phi_E(0+)\,\|X\|\le\|X\|_E\le\phi_E(1)\|X\|$, and together with
$\|XY\|_E\le\|X\|\,\|Y\|_E$ this gives
$C_E(AB)\le C_E(A)+C_E(B)+\ln\frac{\phi_E(1)}{\phi_E(0+)}$. The constant is
optimal: for the operators $P_t$ above, $\|P_t^{\pm k}\|_E=t^k\phi_E(\delta)(1+o(1))$
as $t\to\infty$, $k=1,2$, so
$C_E(P_t^2)-2C_E(P_t)\to\ln\frac{\phi_E(1)}{\phi_E(\delta)}$, and
$\delta\downarrow0$. Hence in type $\mathrm{II}_1$ the defect is finite if and
only if $\phi_E(0+)>0$, i.e.\ if and only if $E=L^\infty$ with an equivalent
norm, and it vanishes if and only if $\phi_E$ is constant, i.e.\
$\|\cdot\|_E$ is a multiple of the operator norm. Thus in a factor of type
$\mathrm{II}_1$ the operator norm is, up to equivalence, the \emph{only}
symmetric norm with finite defect. This is strictly more rigid than the matrix
case: in $M_n(\C)$, $n\ge2$, the same argument, whose inputs are now a rank-one
projection $e\le\mathbf 1_{\{\|X\|\}}(|X|)$, so that $|X|\ge\|X\|e$, the monotonicity of a unitarily
invariant norm $\tnorm{\cdot}$ on positive operators, and the inequality
$\tnorm{XY}\le\|X\|\,\tnorm{Y}$ \cite[Ch.~IV]{Bhatia97}, shows that $\tnorm{\cdot}$, with the
contrast $\frac12\ln\bigl(\tnorm{A}\,\tnorm{A^{-1}}/\tnorm{1}^2\bigr)$, has defect
exactly $\ln\frac{\tnorm{1}}{\tnorm{e}}$, $e$ a rank-one projection. This is at most $\ln n$, since
$\tnorm{1}\le n\tnorm{e}$, and it vanishes if and only if $\tnorm{1}=\tnorm{e}$, i.e.\ (as
$\tnorm{e}\,\|X\|\le\tnorm{X}\le\tnorm{1}\,\|X\|$) if and only if
$\tnorm{\cdot}$ is a multiple of the operator norm. So in the matrix case every unitarily
invariant norm has finite defect, and only the vanishing of the defect singles out the
operator norm.
\end{remark}

\section{The rate of convergence to \texorpdfstring{$\OLC$}{OLC} and the dimension range}\label{tracial:sec:rate}

\begin{theorem}\label{tracial:thm:edge}
Let $(\cM,\tau)$ be a von Neumann algebra with a faithful normal tracial state.
\begin{enumerate}
\item\label{tracial:e:lim} Let $A\in\GL(\cM)$ with $|A|$ not scalar, let $a<b$ be the
minimum and maximum of the spectrum of $\ln|A|$, and
$e_+=\one_{\{b\}}(\ln|A|)$, $e_-=\one_{\{a\}}(\ln|A|)$. Then
\[
\lim_{p\to\infty}2p\bigl(\OLC(A)-C_p(A)\bigr)=-\ln\tau(e_+)-\ln\tau(e_-)\in[2\ln2,\infty],
\]
with the convention $-\ln0=\infty$.
\item\label{tracial:e:fd} $\OLC(A)-C_p(A)=O(1/p)$ as $p\to\infty$ for every
$A\in\GL(\cM)$ if and only if $\cM$ is finite-dimensional.
\item\label{tracial:e:range} If $\cM\cong M_n(\C)$, the limits in \ref{tracial:e:lim} form the
finite set $\{\ln\frac{n^2}{jk}:\ j,k\ge1,\ j+k\le n\}$. If $\cM$ is a factor of
type $\mathrm{II}_1$, they fill the whole interval $[2\ln2,\infty]$.
\end{enumerate}
\end{theorem}

\begin{proof}
\ref{tracial:e:lim} We have
$L_A(p)-pb=\ln\int e^{-p(b-x)}\dd\nu_A(x)$, and the integrand decreases to
$\one_{\{b\}}$ as $p\to\infty$, so $L_A(p)-pb\to\ln\nu_A(\{b\})=\ln\tau(e_+)$ by
monotone convergence. Similarly $L_A(-p)+pa\to\ln\tau(e_-)$. Since
$2p\,\OLC(A)=p(b-a)$, the formula follows. As $e_+\perp e_-$,
$\tau(e_+)+\tau(e_-)\le1$ and $\tau(e_+)\tau(e_-)\le\frac14$.

\ref{tracial:e:fd} If $\cM$ is finite-dimensional, $\ln|A|$ has finite spectrum,
$e_\pm\ne0$, and the limit in \ref{tracial:e:lim} is finite. If $\cM$ is
infinite-dimensional, it contains a sequence $(q_k)_{k\ge1}$ of mutually
orthogonal nonzero projections. Let $H=\sum_{k\ge1}(1-\frac1k)q_k$ and
$A=e^{H}\in\GL(\cM)$. The spectrum of $H$ has maximum $1$ and minimum $0$, but
$\one_{\{1\}}(H)=0$, so the limit in \ref{tracial:e:lim} is $+\infty$.

\ref{tracial:e:range} In $M_n(\C)$, $\tau(e_+)=\frac jn$ and $\tau(e_-)=\frac kn$ with
$j,k\ge1$, $j+k\le n$, and all such pairs occur. In a factor of type
$\mathrm{II}_1$, for any $s,t>0$ with $s+t\le1$ there are orthogonal projections
$e,f$ with $\tau(e)=s$, $\tau(f)=t$, and $A=\exp(e-f)$ has limit $-\ln(st)$;
these values fill $[2\ln2,\infty)$, and $+\infty$ occurs by \ref{tracial:e:fd}.
\end{proof}

Part
\ref{tracial:e:range} shows that the first-order correction to $\OLC$ reads off the
range of the trace on projections, which is the classical invariant
distinguishing type $\mathrm{I}_n$ from type $\mathrm{II}_1$. When an extreme
spectral projection vanishes, the correction is of larger order, determined by
the behaviour of $\nu_A$ at the edges.

\begin{proposition}\label{tracial:prop:logcorr}
Let $A\in\GL(\cM)$ and let $a<b$ be the extreme points of $\supp\nu_A$. Suppose
that, as $u\downarrow0$,
\[
\nu_A\bigl((b-u,b]\bigr)=c_+u^{\beta_+}(1+o(1)),\qquad
\nu_A\bigl([a,a+u)\bigr)=c_-u^{\beta_-}(1+o(1)),
\]
with $c_\pm>0$ and $\beta_\pm>0$. Then, as $p\to\infty$,
\[
C_p(A)=\OLC(A)-\frac{\beta_++\beta_-}{2}\,\frac{\ln p}{p}
+\frac1{2p}\ln\bigl(c_+c_-\Gamma(1+\beta_+)\Gamma(1+\beta_-)\bigr)+o\Bigl(\frac1p\Bigr).
\]
\end{proposition}

\begin{proof}
The hypothesis forces $\nu_A(\{b\})=\nu_A(\{a\})=0$. Let
$F(u)=\nu_A((b-u,b])$, so that $F(0)=0$. Integrating by parts and substituting
$v=pu$,
\[
e^{L_A(p)-pb}=\int_{[0,\infty)}e^{-pu}\,\dd F(u)=\int_0^\infty e^{-v}F(v/p)\,\dd v .
\]
Choose $u_0>0$ with $F(u)\le2c_+u^{\beta_+}$ for $u\le u_0$; since $F\le1$,
$F(u)\le Cu^{\beta_+}$ for all $u>0$ with $C=\max\{2c_+,u_0^{-\beta_+}\}$. Hence
$p^{\beta_+}F(v/p)\le Cv^{\beta_+}$, and by dominated convergence
$p^{\beta_+}e^{L_A(p)-pb}\to c_+\int_0^\infty e^{-v}v^{\beta_+}\dd v=c_+\Gamma(1+\beta_+)$.
The lower edge is treated in the same way, and
$2pC_p(A)=L_A(p)+L_A(-p)$, $2p\,\OLC(A)=p(b-a)$.
\end{proof}

\begin{example}[Exponential of a semicircular element]\label{tracial:ex:semicircular}
Let $s$ be a semicircular element of variance $1$ in a factor of type
$\mathrm{II}_1$ \cite{NS06}, with density $\frac1{2\pi}\sqrt{4-x^2}$ on
$[-2,2]$, and $A=e^{s}$. Then $\tau(A^{p})=\tau(A^{-p})=I_1(2p)/p$, where $I_1$
is the modified Bessel function, and
\[
C_p(A)=\frac1p\ln\frac{I_1(2p)}{p}=2-\frac32\,\frac{\ln p}{p}-\frac{\ln(2\sqrt\pi)}{p}+o\Bigl(\frac1p\Bigr),
\qquad\OLC(A)=2,
\]
in agreement with Proposition~\ref{tracial:prop:logcorr}, where
$\beta_\pm=\frac32$ and $c_\pm=\frac2{3\pi}$. No matrix has this behaviour:
by Theorem~\ref{tracial:thm:edge}\ref{tracial:e:lim},\ref{tracial:e:fd}, in $M_n(\C)$ the quantity
$2p(\OLC-C_p)$ has a finite limit, so the correction is of order $\frac1p$.
\end{example}

\begin{example}[Haar unitaries]\label{tracial:ex:arcsine}
Let $u$ be a Haar unitary in $(\cM,\tau)$, i.e.\ $\tau(u^k)=0$ for $k\ne0$; for
instance a canonical generator of a free group factor $L(\mathbb F_r)$, or of
$L(\mathbb Z)$. For $0<c<1$ let $A=1+\frac c2(u+u^*)$. Then $\frac12(u+u^*)$
has the arcsine distribution on $[-1,1]$, $A$ is positive and invertible, and
\begin{gather*}
\tau(A)=1,\qquad\tau(A^{-1})=\frac1{\sqrt{1-c^2}},\\
C_1(A)=-\tfrac14\ln(1-c^2),\qquad\OLC(A)=\tfrac12\ln\frac{1+c}{1-c}.
\end{gather*}
The arcsine density has inverse square-root singularities at both edges, so
Proposition~\ref{tracial:prop:logcorr} applies with $\beta_\pm=\frac12$ and
$c_\pm=\frac1\pi\sqrt{2(1\pm c)/c}$:
\[
C_p(A)=\OLC(A)-\frac12\,\frac{\ln p}{p}+\frac1{2p}\ln\frac{\sqrt{1-c^2}}{2\pi c}+o\Bigl(\frac1p\Bigr).
\]
Together with Example~\ref{tracial:ex:semicircular} ($\beta_\pm=\frac32$) this shows
that the coefficient of $\frac{\ln p}p$ records the edge exponents of the
spectral measure, here those of the two basic laws of free probability.
\end{example}

\subsection*{The second derivative of the profile: the gap at the edge}

Theorem~\ref{tracial:thm:edge} concerns the first-order behaviour of the profile $L_A$
at infinity, which sees the \emph{mass} of the spectral measure at the edge.
The second derivative $L_A''(t)=\Var(\nu_A^{(t)})$ sees whether the edge is
\emph{isolated}. The following statement consists of standard Laplace-type
asymptotics \cite{BGT87} and is not new. We include it because together with Theorem~\ref{tracial:thm:edge} it gives
the complete picture of the upper edge.

\begin{proposition}[Edge exponent; classical]\label{tracial:prop:edge-exp}
Let $A\in\GL(\cM)$ with $|A|$ not scalar, let $[a,b]$ be the smallest interval
containing the spectrum of $\ln|A|$, and $e_+=\one_{\{b\}}(\ln|A|)$.
\begin{enumerate}
\item If $b$ is an isolated point of the spectrum of $\ln|A|$, at distance
$\delta>0$ from the rest of the spectrum, then
\[
L_A''(t)=\frac{1+o(1)}{\tau(e_+)}\int_{\{x<b\}}(b-x)^2e^{-t(b-x)}\,\dd\nu_A(x)
\qquad\text{and}\qquad-\lim_{t\to\infty}\frac1t\ln L_A''(t)=\delta .
\]
\item If $b$ is not an isolated point of the spectrum, then
$-\lim_{t\to\infty}\frac1t\ln L_A''(t)=0$.
\item If $\nu_A((b-u,b])=c\,u^{\beta}(1+o(1))$ as $u\downarrow0$ with $c,\beta>0$,
then $L_A''(t)=\beta\,t^{-2}(1+o(1))$.
\end{enumerate}
In particular, $L_A''$ decays exponentially if and only if the top of the
spectrum is isolated, and then the rate of decay is the gap.
\end{proposition}

\begin{proof}
Write $\nu=\nu_A$ and $M_k(t)=\int_{\{x<b\}}(b-x)^ke^{-t(b-x)}\dd\nu(x)$ for
$k=0,1,2$. Since the variance is invariant under $x\mapsto b-x$,
\begin{equation}\label{tracial:eq:var-edge}
L_A''(t)=\frac{M_2(t)}{\tau(e_+)+M_0(t)}-\Bigl(\frac{M_1(t)}{\tau(e_+)+M_0(t)}\Bigr)^2 .
\end{equation}

(a) Here $\tau(e_+)>0$ and $\delta\le b-x\le b-a$ on the support of the
restriction of $\nu$ to $\{x<b\}$. Hence $M_0(t)\le e^{-\delta t}$, and by the
Cauchy--Schwarz inequality $M_1(t)^2\le M_2(t)M_0(t)\le M_2(t)e^{-\delta t}$, so
the second term in \eqref{tracial:eq:var-edge} is $o(M_2(t))$ and
$L_A''(t)=M_2(t)(1+o(1))/\tau(e_+)$. Moreover
$M_2(t)\le(b-a)^2e^{-\delta t}$, and for every $\varepsilon>0$ the set
$\{b-\delta-\varepsilon<x\le b-\delta\}$ has positive measure $\eta$, so
$M_2(t)\ge\delta^2\eta\,e^{-(\delta+\varepsilon)t}$. Thus
$\frac1t\ln M_2(t)\to-\delta$.

(b) Since $\nu^{(t)}$ is supported in $[a,b]$,
$L_A''(t)\le\frac14(b-a)^2$, so $\liminf(-\frac1t\ln L_A''(t))\ge0$. For the
converse let $\varepsilon>0$. As $b$ is not isolated, there are two distinct
points of the spectrum in $(b-\varepsilon,b]$, hence disjoint intervals
$U,V\subset(b-\varepsilon,b]$ at positive distance $d$ with $\nu(U),\nu(V)>0$.
Since $\int e^{tx}\dd\nu\le e^{tb}$, $\nu^{(t)}(U)\ge\nu(U)e^{-\varepsilon t}$,
and likewise for $V$. Therefore
$L_A''(t)=\frac12\iint(x-y)^2\dd\nu^{(t)}\dd\nu^{(t)}\ge d^2\nu(U)\nu(V)e^{-2\varepsilon t}$
and $\limsup(-\frac1t\ln L_A''(t))\le2\varepsilon$.

(c) Here $\tau(e_+)=0$. Let $F(u)=\nu((b-u,b])$ and
$G_k(u)=\int_{[0,u]}v^k\dd F(v)=u^kF(u)-k\int_0^uv^{k-1}F(v)\dd v$, so that
$G_k(u)=\frac{c\beta}{\beta+k}u^{\beta+k}(1+o(1))$. The Abelian argument in the
proof of Proposition~\ref{tracial:prop:logcorr}, applied to $G_k$, gives
$M_k(t)=c\,\beta\,\Gamma(\beta+k)\,t^{-(\beta+k)}(1+o(1))$. Hence
$M_1/M_0\sim\beta/t$ and $M_2/M_0\sim\beta(\beta+1)/t^2$, and
\eqref{tracial:eq:var-edge} gives $L_A''(t)\sim\beta/t^2$.
\end{proof}

Combining Theorem~\ref{tracial:thm:edge} and Proposition~\ref{tracial:prop:edge-exp}, there are
three regimes at the top of the spectrum of $|A|$, read off from the top-edge
term $L_A(p)-pb$ of $2p(\OLC-C_p)=(pb-L_A(p))+(-pa-L_A(-p))$. If the top is an
isolated eigenvalue, then $L_A(p)-pb$ has the finite limit $\ln\tau(e_+)$ and
$L_A''$ decays exponentially at the rate given by the gap. If it is an
eigenvalue which is not isolated, the limit is still finite but $L_A''$ decays
subexponentially. If it is not an eigenvalue, then $L_A(p)-pb\to-\infty$, and
for a regular edge with exponent $\beta$ one has
$L_A(p)-pb=-\beta\ln p+O(1)$, a contribution $\frac\beta2\frac{\ln p}p$ to
$\OLC-C_p$, and $L_A''(t)\sim\beta/t^2$; for instance
$L_A''(t)\sim\frac3{2t^2}$ for the exponential of a semicircular element and
$\frac1{2t^2}$ for $1+\frac c2(u+u^*)$ with a Haar unitary $u$. The same
applies to the bottom edge. Only the first regime occurs in finite dimension;
the second requires an eigenvalue that is a limit of other spectral points,
and the third a spectral edge without an eigenvalue.

\section{Araki--Woods factors and sequences of contrasts}\label{tracial:sec:itpfi}

The invariants of Sections~\ref{tracial:sec:defect-vn} and \ref{tracial:sec:rate} depend only
on the range of the trace on projections and say nothing about factors of type
$\mathrm{III}$. For hyperfinite factors the natural object is not the contrast
of one operator but the sequence of contrasts of the component states of an
infinite tensor product.

\subsection{Infinite tensor products of matrix algebras}

Let $n_k\ge2$ be integers and $\rho_k\in M_{n_k}(\C)$ positive definite density
matrices, with states $\varphi_k=\Tr(\rho_k\,\cdot\,)$. On the algebraic
infinite tensor product $\mathfrak A=\bigotimes_{k\ge1}M_{n_k}(\C)$, the union
of the finite tensor products
$\mathfrak A_N=M_{n_1}(\C)\otimes\dots\otimes M_{n_N}(\C)$, consider the product
state $\varphi=\bigotimes_k\varphi_k$ and its GNS representation $\pi_\varphi$.
The von Neumann algebra
\[
\cM=\bigotimes_{k\ge1}\bigl(M_{n_k}(\C),\rho_k\bigr):=\pi_\varphi(\mathfrak A)''
\]
is a hyperfinite factor, called an \emph{Araki--Woods} or \emph{ITPFI} factor
\cite{AW68}; it is called $\mathrm{ITPFI}_2$ if all $n_k=2$. The cyclic vector
defines a faithful normal state on $\cM$ whose modular group is
$\bigotimes_k\operatorname{Ad}\rho_k^{it}$ \cite{Connes73}. For tracial
components $\rho_k=\frac1{n_k}$ one obtains the hyperfinite factor $R$ of type
$\mathrm{II}_1$, and for $\rho_k=\frac1{1+\lambda}\operatorname{diag}(1,\lambda)$,
$0<\lambda<1$, the Powers factor $R_\lambda$ \cite{Powers67}. The restriction of
$\varphi$ to $\mathfrak A_N$ has the density matrix
\[
\rho^{(N)}=\rho_1\otimes\dots\otimes\rho_N\in M_{d_N}(\C),\qquad d_N=n_1\cdots n_N .
\]
For a positive definite matrix $\rho$ we write $\nu_\rho$ for its logarithmic
spectral measure with respect to the normalised trace, $L_\rho$ for its cumulant
generating function and $V(\rho)=\Var(\nu_\rho)$ for the variance of its
logarithmic spectrum.

\begin{lemma}\label{tracial:lem:additive}
For every $N$ and every $p\in(0,\infty)$,
\begin{gather*}
\OLC\bigl(\rho^{(N)}\bigr)=\sum_{k\le N}\OLC(\rho_k),\qquad
C_p\bigl(\rho^{(N)}\bigr)=\sum_{k\le N}C_p(\rho_k),\\
V\bigl(\rho^{(N)}\bigr)=\sum_{k\le N}V(\rho_k).
\end{gather*}
\end{lemma}

\begin{proof}
The eigenvalues of a tensor product are the products of the eigenvalues, so
$\nu_{\rho^{(N)}}=\nu_{\rho_1}*\dots*\nu_{\rho_N}$. Supports and variances add
under convolution, and $C_p$ is additive by Proposition~\ref{tracial:prop:vn}. (For $\OLC$
this is also (Nikolaeva and Novikov, arXiv:2609.06997v2, Prop.~4.7).)
\end{proof}

Thus all contrasts of the partial products are nondecreasing in $N$, and the
question is which of them stay bounded and how fast the others grow. For the
Powers factor with $\OLC(\rho_k)=z$ one finds
$\OLC(\rho^{(N)})=Nz$ and $C_p(\rho^{(N)})=\frac Np\ln\cosh(pz)$: all contrasts
grow linearly in $N$, that is, logarithmically in the dimension $d_N=2^N$.

\subsection{Components of size two: the Araki--Woods criteria in contrast form}

Let all $n_k=2$ and put
\[
z_n=\OLC(\rho_n),\qquad
a_n=C_1(\rho_n),\qquad
b_n=C_1(\rho_n)-\OLC(\rho_n)+\ln2 .
\]
For a positive definite $2\times2$ density matrix $\rho$ one has
$C_1(\rho)=\ln\cosh\OLC(\rho)$, and therefore
\begin{equation}\label{tracial:eq:two-bounds}
\max\{0,\ \OLC(\rho)-\ln2\}\le C_1(\rho)\le\OLC(\rho)
\end{equation}
(in dimension $n$ the same holds with $\ln n$ in place of $\ln2$, by
Proposition~\ref{tracial:prop:vn}(a),(c) and
$\tau(\rho)\tau(\rho^{-1})\ge\frac{\|\rho\|}n\cdot\frac{\|\rho^{-1}\|}n$),
so $a_n\ge0$ and $b_n\ge0$ measure how far $\rho_n$ is from saturating the two
lower bounds for the contrast: $a_n=0$ for the tracial state and $b_n\to0$ as
$\rho_n$ approaches a pure state. Explicitly, $a_n=\ln\cosh z_n$ and
$b_n=\ln(1+e^{-2z_n})$.

The following statement is not new: it is a reformulation of the criteria of
Araki--Woods and of Connes' formula for the invariant $T$. We include it
because in the language of contrasts these classical results take a
strikingly simple form, in which the types $\mathrm{II}_1$ and $\mathrm I$
correspond exactly to the two lower bounds in \eqref{tracial:eq:two-bounds}, and the
Michelson and the logarithmic contrast play the roles of amplitude and
frequency.

\begin{proposition}[Araki--Woods and Connes, reformulated]\label{tracial:prop:itpfi}
Let $\cM=\bigotimes_n(M_2(\C),\rho_n)$. With the above notation:
\begin{enumerate}
\item $\cM$ is of type $\mathrm{II}_1$ if and only if $\sum_na_n<\infty$,
equivalently $\sum_n\OMC(\rho_n)^2<\infty$.
\item $\cM$ is of type $\mathrm{I}$ if and only if $\sum_nb_n<\infty$,
equivalently $\sum_n\bigl(1-\OMC(\rho_n)\bigr)<\infty$, equivalently
$\sum_ne^{-2\OLC(\rho_n)}<\infty$.
\item $\cM$ is of type $\mathrm{III}$ if and only if
$\sum_n\min\{a_n,b_n\}=\infty$. Consequently $\cM$ is of type
$\mathrm{II}_\infty$ if and only if $\sum_na_n=\sum_nb_n=\infty$ and
$\sum_n\min\{a_n,b_n\}<\infty$.
\item Connes' invariant is
$\displaystyle T(\cM)=\Bigl\{t\in\R:\ \sum_{n\ge1}\bigl(1-\OMC(\rho_n)^2\bigr)\,\sin^2\bigl(t\,\OLC(\rho_n)\bigr)<\infty\Bigr\}$.
\end{enumerate}
Moreover, by Lemma~\ref{tracial:lem:additive},
\[
\sum_{n\le N}a_n=C_1\bigl(\rho^{(N)}\bigr),\qquad
\sum_{n\le N}b_n=C_1\bigl(\rho^{(N)}\bigr)-\Bigl(\OLC\bigl(\rho^{(N)}\bigr)-\ln d_N\Bigr),
\]
so that $\cM$ is of type $\mathrm{II}_1$ (resp.\ of type $\mathrm{I}$) if and
only if $C_1(\rho^{(N)})$ stays within bounded distance of the lower bound $0$
(resp.\ of the lower bound $\OLC(\rho^{(N)})-\ln d_N$).
\end{proposition}

\begin{proof}
The eigenvalues of $\rho_n$ are $\lambda_n=e^{z_n}/(2\cosh z_n)$ and
$1-\lambda_n$, with ratio $\kappa_n=e^{2z_n}$, and $\OMC(\rho_n)=\tanh z_n$.
We use \cite[Lemma~2.14]{AW68}.

(a) $\cM$ is of type $\mathrm{II}_1$ iff
$\sum_n\bigl(|\lambda_n^{1/2}-2^{-1/2}|^2+|(1-\lambda_n)^{1/2}-2^{-1/2}|^2\bigr)<\infty$.
The $n$-th term equals $2\bigl(1-\cosh(z_n/2)/\sqrt{\cosh z_n}\bigr)$. This
function of $z_n$, the function $\ln\cosh z_n$ and the function $\tanh^2z_n$
are increasing in $z_n\ge0$, vanish only at $z_n=0$, and are asymptotic to
constant multiples of $z_n^2$ as $z_n\to0$. Convergence of any of the three
series forces its terms, hence $z_n$, to tend to $0$, and then the terms are
comparable; hence the three series converge simultaneously.

(b) $\cM$ is of type $\mathrm{I}$ iff $\sum_n(1-\lambda_n)<\infty$, and
$1-\lambda_n=u_n/(1+u_n)$ with $u_n=e^{-2z_n}\in(0,1]$, while
$b_n=\ln(1+u_n)$ and $1-\tanh z_n=2u_n/(1+u_n)$; all three are between
$\frac12u_n$ and $2u_n$.

(c) Since $\lambda_n\ge\frac12$ is the larger eigenvalue, the criterion of
\cite[Lemma~2.14]{AW68} says that $\cM$ is of type $\mathrm{III}$ iff
$\sum_n(1-\lambda_n)\min\{(1-\kappa_n^{-1})^2,C\}=\infty$ for some (equivalently,
every) constant $C>0$; as $(1-\kappa_n^{-1})^2\le1$, we may take $C=1$, and
the minimum is then $(1-\kappa_n^{-1})^2$. The $n$-th term is a
continuous positive function of $z_n>0$ which is asymptotic to $2z_n^2$ as
$z_n\to0$ and to $e^{-2z_n}$ as $z_n\to\infty$. The same holds, with the
constants $\frac12$ and $1$, for $\min\{\ln\cosh z_n,\ln(1+e^{-2z_n})\}$, so
the two are bounded by constant multiples of each other (terms with $z_n=0$
vanish in both series). The statement on type $\mathrm{II}_\infty$ follows from
(a)--(c), since an $\mathrm{ITPFI}$ factor is of one of the types $\mathrm I$,
$\mathrm{II}_1$, $\mathrm{II}_\infty$, $\mathrm{III}$.

(d) By \cite[Cor.~1.3.9]{Connes73},
$T(\cM)=\{t:\sum_n(1-|\lambda_n^{1+it}+(1-\lambda_n)^{1+it}|)<\infty\}$, and
\[
\bigl|\lambda_n^{1+it}+(1-\lambda_n)^{1+it}\bigr|^2
=1-4\lambda_n(1-\lambda_n)\sin^2(t\,z_n)=1-\bigl(1-\tanh^2z_n\bigr)\sin^2(t\,z_n),
\]
while $\frac u2\le1-\sqrt{1-u}\le u$ for $u\in[0,1]$.
\end{proof}

\begin{corollary}\label{tracial:cor:dichotomy}
If $\sup_n\OLC(\rho_n)<\infty$, then $\bigotimes_n(M_2(\C),\rho_n)$ is of type
$\mathrm{II}_1$ when $\sum_n\OLC(\rho_n)^2<\infty$ and of type $\mathrm{III}$
otherwise.
\end{corollary}

\begin{proof}
If $z_n\le Z$ for all $n$, then $b_n\ge\ln(1+e^{-2Z})>0$, so
$\min\{a_n,b_n\}$ is comparable with $\min\{a_n,1\}$, and
$\sum_n\min\{a_n,1\}<\infty$ iff $\sum_na_n<\infty$ iff $\sum_nz_n^2<\infty$.
\end{proof}

\subsection{Components of arbitrary size with bounded contrast}

For components of arbitrary size the criterion for type $\mathrm{II}_1$ can
still be expressed through the contrasts of the partial products, provided the
contrasts of the components are uniformly bounded.

\begin{theorem}\label{tracial:thm:itpfi-bounded}
Let $\cM=\bigotimes_k(M_{n_k}(\C),\rho_k)$ with $n_k<\infty$ arbitrary and
$\sup_k\OLC(\rho_k)<\infty$. The following are equivalent:
\begin{enumerate}
\item $\cM$ is of type $\mathrm{II}_1$;
\item $\sup_NC_p(\rho^{(N)})=\sum_kC_p(\rho_k)<\infty$ for some, equivalently
for every, $p\in(0,\infty)$;
\item $\sup_NV(\rho^{(N)})=\sum_kV(\rho_k)<\infty$.
\end{enumerate}
Moreover, under the hypothesis $\sup_k\OLC(\rho_k)<\infty$ the factor $\cM$ is
never of type $\mathrm{I}$.
\end{theorem}

\begin{proof}
Let $Z=\sup_k\OLC(\rho_k)$, let $\mu_{k,j}=n_k\lambda_{k,j}$ be the eigenvalues
of $h_k=n_k\rho_k$, and let $X_k$ be a random variable with law $\nu_{h_k}$,
i.e.\ taking the values $\ln\mu_{k,j}$ with equal probabilities. Then
$\mathbb E\,e^{X_k}=\tau(h_k)=1$, so the smallest and the largest of the
$\mu_{k,j}$ lie on different sides of $1$, and since their ratio is at most
$e^{2Z}$, $|X_k|\le2Z$. By \cite[Lemma~2.14]{AW68}, $\cM$ is of type
$\mathrm{II}_1$ if and only if $\sum_kS_k<\infty$, where
\[
S_k=\sum_j\bigl|\lambda_{k,j}^{1/2}-n_k^{-1/2}\bigr|^2=\frac1{n_k}\sum_j\bigl(\mu_{k,j}^{1/2}-1\bigr)^2
=\mathbb E\bigl(e^{X_k/2}-1\bigr)^2 .
\]
We show that $S_k$, $V_k=V(\rho_k)=\Var X_k$ and $C_p(\rho_k)$ are bounded by
constant multiples of each other, with constants depending only on $Z$ and $p$;
together with Lemma~\ref{tracial:lem:additive} this proves the equivalences.

Since $(e^{x/2}-1)^2$ is bounded above and below by constant multiples of $x^2$
on $[-2Z,2Z]$, $S_k$ is comparable with
$\mathbb EX_k^2=V_k+(\mathbb EX_k)^2$. Put $m=\mathbb EX_k$ and $Y=X_k-m$, so
$|Y|\le4Z$. By Jensen's inequality $e^m\le\mathbb Ee^{X_k}=1$, i.e.\ $m\le0$,
and from $e^y\le1+y+\frac12e^{4Z}y^2$ for $|y|\le4Z$ we get
$1=e^m\,\mathbb Ee^{Y}\le e^m(1+\frac12e^{4Z}V_k)$, hence
$0\le-m\le\frac12e^{4Z}V_k$. As $V_k\le4Z^2$, this gives
$m^2\le Z^2e^{8Z}V_k$, and $\mathbb EX_k^2$ is comparable with $V_k$.

For the contrast, let $X_k'$ be an independent copy of $X_k$. Then
$e^{2pC_p(\rho_k)}=\mathbb E\cosh\bigl(p(X_k-X_k')\bigr)$,
$\mathbb E(X_k-X_k')^2=2V_k$ and $|X_k-X_k'|\le4Z$, so
\[
1+p^2V_k\le e^{2pC_p(\rho_k)}\le1+p^2\cosh(4pZ)\,V_k ,
\]
and $\ln(1+u)$ is comparable with $u$ for $0\le u\le4p^2Z^2\cosh(4pZ)$.

Finally, $\lambda_{k,1}\le e^{2Z}\lambda_{k,n_k}$ and
$1-\lambda_{k,1}\ge\lambda_{k,n_k}$ give $\lambda_{k,1}\le(1+e^{-2Z})^{-1}$, so
$\sum_k(1-\lambda_{k,1})=\infty$ and $\cM$ is not of type $\mathrm{I}$
\cite[Lemma~2.14]{AW68}.
\end{proof}

When the conditions of Theorem~\ref{tracial:thm:itpfi-bounded} fail, $\cM$ is of type
$\mathrm{II}_\infty$ or $\mathrm{III}$, and the theorem does not decide
between the two; for components of size two Corollary~\ref{tracial:cor:dichotomy}
shows that the type is $\mathrm{III}$.

In words: an Araki--Woods factor with uniformly bounded component contrasts is
of type $\mathrm{II}_1$ exactly when the Schatten contrasts of the partial
products $\rho^{(N)}$ stay bounded, or equivalently when the logarithmic
spectra of the $\rho^{(N)}$ have bounded variance, while their operator-norm
contrasts $\OLC(\rho^{(N)})=\sum_{k\le N}\OLC(\rho_k)$ may well diverge.

\medskip\noindent\textbf{Arbitrary component contrasts.} Without the bound on
the component contrasts, the contrasts of the partial products still give a
sufficient condition for type $\mathrm{II}_1$, but only for $p\ge1$.

\begin{proposition}\label{tracial:prop:itpfi-p1}
Let $\cM=\bigotimes_k(M_{n_k}(\C),\rho_k)$ with arbitrary $n_k$. If
$\sum_kC_p(\rho_k)<\infty$ for some $p\ge1$ (it is enough to assume this for $p=1$), then $\cM$ is of type
$\mathrm{II}_1$.
\end{proposition}

\begin{proof}
Let $h_k=n_k\rho_k$, so that $\tau(h_k)=1$. As in the proof of
Theorem~\ref{tracial:thm:itpfi-bounded}, $\cM$ is of type $\mathrm{II}_1$ iff
$\sum_kS_k<\infty$, where
$S_k=\sum_j\bigl|\lambda_{k,j}^{1/2}-n_k^{-1/2}\bigr|^2=2-2\tau(h_k^{1/2})$. Write
$0=\theta\cdot\frac12+(1-\theta)(-p)$ with $\theta=\frac{2p}{1+2p}$. By H\"older's
inequality (log-convexity of $s\mapsto\tau(h^s)$),
$1=\tau(h_k^0)\le\tau(h_k^{1/2})^\theta\,\tau(h_k^{-p})^{1-\theta}$, i.e.\
$\ln\tau(h_k^{1/2})\ge-\frac1{2p}\ln\tau(h_k^{-p})$. By Jensen's inequality
$\tau(h_k^p)\ge\tau(h_k)^p=1$ for $p\ge1$, so
$\frac1{2p}\ln\tau(h_k^{-p})\le C_p(h_k)=C_p(\rho_k)$. Hence
$\tau(h_k^{1/2})\ge e^{-C_p(\rho_k)}$ and
$S_k\le2\bigl(1-e^{-C_p(\rho_k)}\bigr)\le2C_p(\rho_k)$.
\end{proof}

Put $p_c(\rho_\bullet)=\sup\{p>0:\ \sum_kC_p(\rho_k)<\infty\}$ (with
$\sup\emptyset=0$); since $C_p$ is nondecreasing in $p$, the sum is finite on
$(0,p_c)$ and infinite on $(p_c,\infty)$. By Theorem~\ref{tracial:thm:itpfi-bounded},
if $\sup_k\OLC(\rho_k)<\infty$ then $p_c\in\{0,\infty\}$, and $p_c=\infty$ exactly
for type $\mathrm{II}_1$. In general the threshold is not determined by the
type, and the condition $p\ge1$ in Proposition~\ref{tracial:prop:itpfi-p1} cannot be
dropped. In the following example $h_k=n_k\rho_k$ and $S_k$ are as in the
proof above, and we use that $C_p(\rho_k)=C_p(h_k)$.

\begin{example}\label{tracial:ex:itpfi-pc}
Put $t_k=2^{-k}$.
\begin{enumerate}
\item \emph{Type $\mathrm{II}_1$ with $p_c=2$.} Let $n_k=4^k$ and
$h_k=c_k\operatorname{diag}(2^k,1,\dots,1)$, $c_k=4^k/(2^k+4^k-1)$. Then
$c_k\ge1-t_k$, so
$\tau(h_k^{1/2})\ge c_k^{1/2}(1-4^{-k})\ge1-2t_k$ and $S_k\le4t_k$: $\cM$ is of
type $\mathrm{II}_1$. Further
\[
2p\,C_p(\rho_k)=\ln\bigl(1+(2^{pk}-1)4^{-k}\bigr)+\ln\bigl(1-(1-2^{-pk})4^{-k}\bigr).
\]
For $p<2$ this is at most $2^{(p-2)k}$, so $\sum_kC_p(\rho_k)<\infty$. For $p=2$
and $k\ge1$ the product of the two arguments is at least
$\frac74\cdot\frac34>1$, so $C_2(\rho_k)$ is bounded below by a positive
constant and $\sum_kC_2(\rho_k)=\infty$.
\item \emph{Not of type $\mathrm{II}_1$, with $p_c=1$.} Let $n_k=2^k$ and let
$h_k$ have the eigenvalue $\frac1{2t_k}$ once and the eigenvalue
$\frac1{2(1-t_k)}$ with multiplicity $2^k-1$. Then
$\tau(h_k^{1/2})=\sqrt{t_k/2}+\sqrt{(1-t_k)/2}\to2^{-1/2}$, so $S_k\to2-\sqrt2$ and
$\sum_kS_k=\infty$. For $0<p<1$,
\begin{align*}
\tau(h_k^p)\,\tau(h_k^{-p})
&=\bigl(t_k^{1-p}+(1-t_k)^{1-p}\bigr)\bigl(t_k^{1+p}+(1-t_k)^{1+p}\bigr)\\
&\le(1+t_k^{1-p})(1+t_k^{1+p})\le1+3t_k^{1-p},
\end{align*}
so $C_p(\rho_k)\le\frac3{2p}2^{-(1-p)k}$ and the sum converges; for $p=1$ and
$k\ge2$, $\tau(h_k)\tau(h_k^{-1})=2-4t_k+4t_k^2\ge\frac54$, so
$\sum_kC_1(\rho_k)=\infty$. In view of Proposition~\ref{tracial:prop:itpfi-p1} this
example is sharp.
\item \emph{Type $\mathrm{II}_1$ with $p_c=0$.} Let $n_k=2^k$ and let $h_k$
have the eigenvalue $\varepsilon_k=\exp(-e^k)$ once and the eigenvalue
$\eta_k=(1-t_k\varepsilon_k)/(1-t_k)\ge1$ with multiplicity $2^k-1$. Then
$\tau(h_k^{1/2})\ge(1-t_k)\eta_k^{1/2}\ge1-t_k$, so $S_k\le2t_k$ and $\cM$ is of
type $\mathrm{II}_1$; but for every $p>0$,
$2p\,C_p(\rho_k)\ge\ln\bigl((1-t_k)\eta_k^p\bigr)+\ln\bigl(t_k\varepsilon_k^{-p}\bigr)\ge-(k+1)\ln2+p\,e^k\to\infty$,
so $\sum_kC_p(\rho_k)=\infty$.
\end{enumerate}
\end{example}

Thus, without the bound on the component contrasts, finiteness of
$\sum_kC_p(\rho_k)$ for some $p\ge1$ implies type $\mathrm{II}_1$, finiteness
for $p<1$ does not, and a factor of type $\mathrm{II}_1$ may have $p_c=2$ or
$p_c=0$.

\begin{remark}\label{tracial:rem:modular}
The contrasts of the partial products are contrasts of modular operators. For
a positive definite density matrix $\rho\in M_d(\C)$, let $M_d(\C)$ act by left
multiplication on itself with the inner product $\langle x,y\rangle=\Tr(x^*y)$
and cyclic vector $\rho^{1/2}$; the modular operator of $\Tr(\rho\,\cdot\,)$ is
then $\Delta_\rho(x)=\rho x\rho^{-1}$, with eigenvalues $\lambda_i/\lambda_j$
and the matrix units as eigenvectors. With the normalised trace of
$B(M_d(\C))\cong M_{d^2}(\C)$ its logarithmic spectral measure is
$\nu_\rho*\check\nu_\rho$, so $C_p(\Delta_\rho)=2C_p(\rho)$ by
Proposition~\ref{tracial:prop:vn}(e), and $\OLC(\Delta_\rho)=2\OLC(\rho)$ because
diameters of supports add under convolution. Hence
$\sum_{k\le N}C_p(\rho_k)=\frac12C_p(\Delta_{\rho^{(N)}})$, where
$\Delta_{\rho^{(N)}}$ is the modular operator of the restriction of $\varphi$ to
$\mathfrak A_N$. These quantities are defined through the chosen tensor
product decomposition; we do not know whether they are determined by
$(\cM,\varphi)$.
\end{remark}

\subsection{Connes' invariant and the spectral profile}

For an arbitrary ITPFI factor Connes' formula \cite[Cor.~1.3.9]{Connes73} reads
$T(\cM)=\{t:\sum_k(1-|\Tr\rho_k^{1+it}|)<\infty\}$. Since
$\Tr\rho^{1+it}=\tau(\rho^{1+it})/\tau(\rho)$ is the characteristic function of
the tilted logarithmic spectral measure $\nu_\rho^{(1)}$,
\[
T(\cM)=\Bigl\{t\in\R:\ \sum_k\bigl(1-|\chi_k(t)|\bigr)<\infty\Bigr\},
\qquad
\chi_k(t)=\int e^{itx}\,\dd\nu_{\rho_k}^{(1)}(x)=\frac{\tau(\rho_k^{1+it})}{\tau(\rho_k)} .
\]
In other words, the invariant $T$ is determined by the analytic continuation of
the functions $e^{L_{\rho_k}(s)}=\tau(\rho_k^{s})$, i.e.\ of the logarithmic
spectral profiles of the components, to the line
$\operatorname{Re}s=1$, just as the Schatten contrasts are determined by them
on the real axis. For components
of size two this is Proposition~\ref{tracial:prop:itpfi}(d). Since
$1-|\chi_k(t)|\le1-|\chi_k(t)|^2\le t^2\Var(\nu_{\rho_k}^{(1)})$, summability of
these variances (of the tilted measures $\nu^{(1)}_{\rho_k}$, which are comparable with the
variances $V(\rho_k)$ when the contrasts are bounded) implies $T(\cM)=\R$, in accordance with
Theorem~\ref{tracial:thm:itpfi-bounded}. We also note that Connes' invariant $S(\cM)$
of an ITPFI factor coincides with the asymptotic ratio set $r_\infty(\cM)$ of
Araki and Woods \cite{AW68,Connes73}. We do not know whether $S(\cM)$ can be
read off from the logarithmic spectral profiles of the tail products
$\bigotimes_{k\ge m}\rho_k$, $m\to\infty$.

Recall that $T(\cM)=\R$ for semifinite factors,
$T(\cM)=\frac{2\pi}{|\ln\lambda|}\mathbb Z$ in type $\mathrm{III}_\lambda$,
$0<\lambda<1$, and $T(\cM)=\{0\}$ in type $\mathrm{III}_1$ \cite{Connes73}.

\begin{example}[Growth classes of the contrast]\label{tracial:ex:growth-classes}
Let all $n_k=2$ and $z_n=\OLC(\rho_n)$.

(i) \emph{Bounded contrast.} If $z_n=z>0$ for all $n$, then $\cM$ is the Powers
factor of type $\mathrm{III}_\lambda$ with $\lambda=e^{-2z}=\kappa(\rho_n)^{-1}$,
and $T(\cM)=\frac\pi z\mathbb Z$.

(ii) \emph{Polynomial condition numbers.} If $\kappa(\rho_n)=n^\alpha$, i.e.\
$z_n=\frac\alpha2\ln n$, then $b_n$ is comparable with $n^{-\alpha}$ and
$a_n\to\infty$. Hence $\cM$ is of type $\mathrm{I}$ for $\alpha>1$ and of type
$\mathrm{III}$ for $0<\alpha\le1$. In the latter case the series in
Proposition~\ref{tracial:prop:itpfi}(d) is comparable with
$\sum_nn^{-\alpha}\sin^2(\frac{t\alpha}2\ln n)$, which diverges for every
$t\ne0$; hence $T(\cM)=\{0\}$, and $\cM$ is not of type
$\mathrm{III}_\lambda$ for any $\lambda\in(0,1)$. The invariant $T$ alone does
not decide between types $\mathrm{III}_0$ and $\mathrm{III}_1$ here.

(iii) \emph{Factorial contrasts.} Let the contrast $z_j=\frac{\ln2}2\,j!$ occur
with multiplicity $2^{j!}$, $j\ge1$; factors of this kind appear in
\cite{ST10}. Since $b$ is comparable with
$e^{-2z_j}=2^{-j!}$ on the $j$-th block, $\sum\min\{a_n,b_n\}=\infty$ and $\cM$
is of type $\mathrm{III}$. Moreover $T(\cM)$ is the set of $t$ with
$\sum_j\sin^2(\frac{t\ln2}2\,j!)<\infty$, which contains
$\frac{2\pi}{\ln2}\mathbb Q$ because $r\,j!\in\mathbb Z$ for every
$r\in\mathbb Q$ and all large $j$. Thus $T(\cM)\ne\{0\}$, and $T(\cM)$ is not of
the form $\frac{2\pi}{|\ln\lambda|}\mathbb Z$ since it contains a dense
subgroup; by the values of $T$ recalled above, $\cM$ is of type
$\mathrm{III}_0$.
\end{example}

Thus bounded contrasts give type $\mathrm{II}_1$ or $\mathrm{III}$
(Corollary~\ref{tracial:cor:dichotomy}), with type $\mathrm{III}_\lambda$ for constant
contrast; polynomially growing condition numbers give type $\mathrm{I}$ or a
type $\mathrm{III}$ factor with trivial $T$-invariant; and the factorially
growing contrasts of Example~\ref{tracial:ex:growth-classes}(iii), taken with matching
multiplicities, give type $\mathrm{III}_0$. The type is, however, not a
function of the growth rate alone: the arithmetic of the sequence of contrasts
enters through the factor $\sin^2(t\,\OLC(\rho_n))$ in
Proposition~\ref{tracial:prop:itpfi}(d). In Example~\ref{tracial:ex:growth-classes}(iii) it is
the fact that $r\,j!$ is eventually an integer for every rational $r$, and not
the size of $j!$, that makes $T(\cM)$ dense.

\section{Concluding remarks}

In a tracial von Neumann algebra the Schatten contrast turns out to measure
three things: whether the functional is a trace, how small projections can be,
and how the spectral measure of $\ln|A|$ behaves at its edges; for infinite
tensor products, the sequence of contrasts of the components reproduces the
Araki--Woods type criteria for components of size two, and the type $\mathrm{II}_1$
criterion for components of bounded contrast. We mention
open questions.

\begin{enumerate}[label=\textup{(\arabic*)}]
\item \emph{Finer invariants.} The invariants of Sections~\ref{tracial:sec:defect-vn}
and \ref{tracial:sec:rate} depend only on the range of the trace on projections, and
therefore cannot separate non-isomorphic factors of type $\mathrm{II}_1$. Do
joint contrasts of several noncommuting operators, as in the equality problem
of Nikolaeva and Novikov \cite{NN26}, carry finer information?
\item \emph{Type $\mathrm{III}$.} Is there an analogue of
Theorem~\ref{tracial:thm:trace} for weights on factors of type $\mathrm{III}$, with
Haagerup $L^p$-spaces in place of $L^p(\cM,\tau)$? Which growth and arithmetic
conditions on the contrasts and on the profiles $L_{\rho_k}$ of the component
states characterise the types $\mathrm{III}_0$, $\mathrm{III}_\lambda$,
$\mathrm{III}_1$ of Araki--Woods factors with components of arbitrary size?
\end{enumerate}

\chapter{Infinite contrasts and finite sections}\label{ch:infinite}
\section{Introduction}

Let $A$ be a bounded invertible operator with condition number
$\kappa(A)=\|A\|\,\|A^{-1}\|$. Its \emph{operator Michelson contrast} and its
\emph{operator logarithmic contrast} \cite{ANN24,NN26} are
\[
\OMC(A)=\frac{\kappa(A)-1}{\kappa(A)+1}\qquad\text{and}\qquad
\OLC(A)=\tfrac12\ln\kappa(A)=\operatorname{artanh}\OMC(A);
\]
for positive $A$,
$\OMC(A)=(\sup\sigma(A)-\inf\sigma(A))/(\sup\sigma(A)+\inf\sigma(A))$.
Replacing the operator norm by the Schatten norms of a von Neumann algebra
with a faithful normal tracial state $\tau$ gives the Schatten contrasts
\[
C_p(A)=\tfrac12\ln\bigl(\|A\|_p\,\|A^{-1}\|_p\bigr),\qquad \|X\|_p=\tau(|X|^p)^{1/p}.
\]
All these quantities are finite for bounded invertible operators. If
$A$ or $A^{-1}$ is unbounded, then $\OLC(A)=\infty$ and $\OMC(A)=1$. The question
of this chapter is \emph{how} such a contrast is infinite: which spectral
information survives, and in which invariant form. We give three answers, one
for operators affiliated with a finite von Neumann algebra (through the
function $p\mapsto C_p(A)$), one for unbounded operators on a Hilbert space
(through the growth of the contrasts of finite sections, and through upper and
lower limits of that growth), and one for bounded operators whose spectrum is
read modulo compact operators.

The finite-section phase transition and the bounds for oscillating Weyl constants are the main quantitative results. The parts on moment thresholds and essential spectra organise classical Laplace-transform and quotient-algebra facts in the same notation. Precise assumptions and proofs are given in the corresponding sections.

\section{Preliminaries}\label{infinite:sec:prelim}

Throughout, $\tau$ is a faithful normal tracial state on a von Neumann algebra
$\cM$, and $\widetilde\cM$ denotes the $*$-algebra of closed, densely defined
operators affiliated with $\cM$; as $\tau$ is finite, each of them is
$\tau$-measurable \cite{Nelson74,FK86}. Let $A\in\widetilde\cM$ be injective.
Its right support equals $1$, and its left support, being equivalent to the right
one in the finite algebra $\cM$, equals $1$ as well. Hence $A$ has dense range,
$A^{-1}\in\widetilde\cM$, the partial isometry $U$ in the polar decomposition
$A=U|A|$ is unitary, and $|A^{-1}|=U|A|^{-1}U^*$. We call the Borel probability
measure $\nu_A(E)=\tau\bigl(\one_E(\ln|A|)\bigr)$ on $\R$ the
\emph{logarithmic spectral measure} of $A$; thus
$\tau(|A|^p)=\int e^{px}\dd\nu_A$ and $\tau(|A^{-1}|^p)=\int e^{-px}\dd\nu_A$,
both in $(0,\infty]$. Each Borel probability measure on $\R$ is of the form
$\nu_A$, already for $\cM=L^\infty[0,1]$. For a Borel probability measure $\nu$
on $\R$ let $L_\nu(t)=\ln\int e^{tx}\dd\nu(x)\in(-\infty,\infty]$, and for
$0<p<\infty$ set
\[
C_p(\nu)=\frac{L_\nu(p)+L_\nu(-p)}{2p}\in[0,\infty],\qquad
C_\infty(\nu)=\tfrac12\operatorname{diam}\supp\nu\in[0,\infty].
\]
For injective $A\in\widetilde\cM$ we put $C_p(A)=C_p(\nu_A)$ and
$\OLC(A)=C_\infty(\nu_A)$; with $L_A=L_{\nu_A}$ this reads
$2p\,C_p(A)=L_A(p)+L_A(-p)$. The \emph{critical exponent} of $\nu$ is
$p_c(\nu)=\sup\{p>0:\ C_p(\nu)<\infty\}$, where $\sup\emptyset=0$.

Proposition~\ref{tracial:prop:vn} applies without a bounded-support assumption: it gives positivity, the exponential-moment criterion, monotonicity in $p$, the bounded-support estimate, and the dilation and convolution laws. We use those properties below for extended-valued contrasts.

\begin{remark}\label{infinite:rem:nosub}
At a fixed exponent, $C_p$ admits no subadditivity estimate on injective
elements of $\widetilde\cM$: for a positive $A$ with $p<p_c(\nu_A)<2p$ we have $C_p(A)<\infty$ but
$C_p(A^2)=2C_{2p}(A)=\infty$ by Proposition~\ref{tracial:prop:vn}(e),(c) and the definition of $p_c$. If
the exponents are combined harmonically, subadditivity holds with no additive correction term
(Proposition~\ref{infinite:prop:holder}).
\end{remark}

\section{Unbounded operators: critical exponents and growth}\label{infinite:sec:unbounded}

In this section $A$ is an injective element of $\widetilde\cM$, possibly
unbounded and with unbounded inverse, and $C_p(A)=C_p(\nu_A)$ as in
Section~\ref{infinite:sec:prelim}. We first treat the case $\OLC(A)=\infty$, $p_c(\nu_A)=\infty$, where the rate at which $C_p(A)$ diverges is a quantitative substitute for the contrast, governed by the tails of $\nu_A$, and give first examples with finite $p_c$. In
Subsections~\ref{infinite:subsec:cfun}--\ref{infinite:subsec:transition} we then consider the
whole contrast function, its threshold $p_c$, and the transition at $p_c$.

\subsection{Growth of the contrast, and first examples}

\begin{theorem}[Growth of the contrast]\label{infinite:thm:growth}
Let $\nu$ be a Borel probability measure on $\R$ with
$\int e^{p|x|}\dd\nu<\infty$ for all $p>0$. Put
$M(p)=\ln\int e^{p|x|}\dd\nu$ and $\Theta(u)=\nu(\{x:|x|\ge u\})$.
\begin{enumerate}
\item\label{infinite:gr:sand} For all $p>0$,
\[
\frac{M(p)-\ln2}{2p}-\frac12\int|x|\,\dd\nu\;\le\;C_p(\nu)\;\le\;\frac{M(p)}p .
\]
\item\label{infinite:gr:tail} Suppose there are $q>1$, $c_1\ge c_2>0$ and $u_0\ge0$ such
that $e^{-c_1u^q}\le \Theta(u)\le e^{-c_2u^q}$ for all $u\ge u_0$. Then, with
$\eta(c)=(q-1)\,q^{-q/(q-1)}\,c^{-1/(q-1)}$,
\[
\tfrac12\eta(c_1)\;\le\;\liminf_{p\to\infty}\frac{C_p(\nu)}{p^{1/(q-1)}}
\;\le\;\limsup_{p\to\infty}\frac{C_p(\nu)}{p^{1/(q-1)}}\;\le\;\eta(c_2).
\]
\end{enumerate}
\end{theorem}

\begin{proof}
\ref{infinite:gr:sand} Since $\int e^{\pm px}\dd\nu\le\int e^{p|x|}\dd\nu$, the upper
bound is clear. For the lower bound,
$e^{M(p)}\le\int e^{px}\dd\nu+\int e^{-px}\dd\nu\le2\max\bigl(\int e^{px}\dd\nu,\int e^{-px}\dd\nu\bigr)$,
while by Jensen's inequality each of the two integrals is at least
$e^{-p\int|x|\dd\nu}$. Hence
$\int e^{px}\dd\nu\int e^{-px}\dd\nu\ge\frac12e^{M(p)}e^{-p\int|x|\dd\nu}$.

\ref{infinite:gr:tail} Let $q'=q/(q-1)$ and
$S_c(p)=\sup_{u\ge0}(pu-cu^q)=\eta(c)\,p^{q'}$, attained at
$u_p=(p/(cq))^{1/(q-1)}$. For the lower bound, by Chebyshev's inequality
$e^{M(p)}\ge e^{pu}\Theta(u)$ for every $u\ge0$; taking $u=u_p$ with $c=c_1$, which
exceeds $u_0$ for large $p$, gives $M(p)\ge\eta(c_1)p^{q'}$. For the upper
bound, by Fubini's theorem
\[
e^{M(p)}=1+p\int_0^\infty e^{pu}\,\nu(\{|x|>u\})\,\dd u
\le e^{pu_0}+p\int_0^\infty e^{pu-c_2u^q}\,\dd u .
\]
Let $U=(2p/c_2)^{1/(q-1)}$. For $u\ge U$ we have $c_2u^q\ge2pu$, so the
integrand is at most $e^{-pu}$ and this part of the integral is at most $1/p$.
For $u\le U$ the integrand is at most $e^{S_{c_2}(p)}$. Hence
$e^{M(p)}\le e^{pu_0}+1+pU\,e^{\eta(c_2)p^{q'}}$, and since $q'>1$,
$M(p)\le\eta(c_2)p^{q'}+O(\ln p)$. Inserting both estimates into
\ref{infinite:gr:sand} and using $q'-1=1/(q-1)$ gives the claim.
\end{proof}

Thus a logarithmic spectrum with tails of order $e^{-cu^q}$ gives
$C_p\asymp p^{1/(q-1)}$: the lighter the tails, the slower the growth, and
bounded spectra ($q=\infty$) give bounded contrast. Sharper asymptotics for
regularly varying tail exponents follow from Kasahara's Tauberian theorem
\cite[Section~4.12]{BGT87}.

\begin{example}[Log-normal operators]\label{infinite:ex:lognormal}
If $\nu_A$ is the nondegenerate Gaussian law $N(m,\sigma^2)$, $\sigma>0$, then
$L_A(t)=mt+\frac12\sigma^2t^2$ and
\[
C_p(A)=\frac{\sigma^2p}{2}\qquad(0<p<\infty),\qquad\OLC(A)=\infty .
\]
Here $A$ and $A^{-1}$ are unbounded, all contrasts are finite, and the growth is
linear. This is consistent with Theorem~\ref{infinite:thm:growth}: for $q=2$ the tail
hypothesis holds with any $c_1>\frac1{2\sigma^2}>c_2$, and
$\eta(\frac1{2\sigma^2})=\frac{\sigma^2}2$, so the theorem gives
$\frac{\sigma^2}4\le\liminf_pC_p/p\le\limsup_pC_p/p\le\frac{\sigma^2}2$; the
upper bound is attained.
\end{example}

\begin{example}[Finite critical exponent]\label{infinite:ex:laplace}
If $\nu_A$ has density $\frac1{2b}e^{-|x|/b}$, then
$\int e^{\pm px}\dd\nu_A=(1-b^2p^2)^{-1}$ for $p<1/b$, so
$C_p(A)=-\frac1p\ln(1-b^2p^2)$ for $p<1/b$ and $p_c(\nu_A)=1/b$.
\end{example}

\begin{example}[Free Poisson elements]\label{infinite:ex:MP}
Let $A\in\cM$ be positive with the Marchenko--Pastur distribution (a scaled free Poisson law)
of parameter $y\in(0,1]$ \cite{MP67,NS06}, i.e.\ with density
$\frac{\sqrt{(b-\lambda)(\lambda-a)}}{2\pi y\lambda}$ on $[a,b]$,
$a=(1-\sqrt y)^2$, $b=(1+\sqrt y)^2$. This is the limiting spectral distribution
of sample covariance matrices $\frac1NXX^*$ with $X$ of size $n\times N$,
$n/N\to y$. For $y<1$, using $\tau(A)=1$, $\tau(A^2)=1+y$,
$\tau(A^{-1})=(1-y)^{-1}$ and $\tau(A^{-2})=(1-y)^{-3}$,
\[
C_1(A)=-\tfrac12\ln(1-y),\qquad
C_2(A)=\tfrac14\ln\frac{1+y}{(1-y)^3},\qquad
\OLC(A)=\ln\frac{1+\sqrt y}{1-\sqrt y}.
\]
As $y\uparrow1$ all three diverge, at the different rates
$\frac12|\ln(1-y)|$, $\frac34|\ln(1-y)|$ and $|\ln(1-y)|$. At $y=1$ the
operator $A$ is injective and bounded with unbounded inverse, the density is
$\frac1{2\pi}\sqrt{(4-\lambda)/\lambda}$ on $(0,4]$, and by the Beta integral
and the duplication formula
\begin{align*}
\tau(A^{s})&=\frac{\Gamma(2s+1)}{\Gamma(s+1)\Gamma(s+2)}&&(s>-\tfrac12),\\
C_p(A)&=\frac1{2p}\ln\frac{\Gamma(1+2p)\,\Gamma(1-2p)}{\Gamma(1+p)\Gamma(2+p)\Gamma(1-p)\Gamma(2-p)}
&&(0<p<\tfrac12).
\end{align*}
Hence $\OLC(A)=\infty$ and $p_c(\nu_A)=\frac12$: the critical exponent records
the square-root singularity of the spectral density at the hard edge. More
generally, for a bounded positive $A$ whose spectral density behaves like
$\lambda^{\alpha}$ at $0$, $\alpha>-1$, one has $p_c=\alpha+1$. For $y<1$ the
largest and the smallest eigenvalue of $\frac1NXX^*$ converge almost surely to
$b$ and to $a>0$ under the moment assumptions of \cite{YBK88,BY93}, so the
values above are the almost sure limits of the matrix contrasts.
\end{example}

\subsection{The contrast function, and orders of growth and decay}\label{infinite:subsec:cfun}

If $A$ or $A^{-1}$ is unbounded, then $\OLC(A)=\infty$, but the whole function
$p\mapsto C_p(A)$ still carries spectral information. We collect this
information in a few numbers defined through traces of powers only. Since
$C_p(A)=C_p(\nu_A)$, they are characteristics of the probability measure
$\nu_A$, and most of the facts below are classical properties of two-sided
Laplace transforms \cite{Widder41}; we state them because, read as statements
about the contrast, they describe how an infinite contrast is infinite.

\begin{definition}\label{infinite:def:cfun}
Let $A\in\widetilde\cM$ be injective. The \emph{contrast function} of $A$ is
$p\mapsto C_p(A)\in[0,\infty]$, $p\in(0,\infty]$. Put
\[
p_+(A)=\sup\Bigl\{p\ge0:\ \tau(|A|^p)<\infty\Bigr\},\qquad
p_-(A)=\sup\Bigl\{p\ge0:\ \tau(|A^{-1}|^p)<\infty\Bigr\},
\]
and call $\omega_+(A)=1/p_+(A)$ and $\omega_-(A)=1/p_-(A)$ (with $\frac10=\infty$,
$\frac1\infty=0$) the \emph{order of growth} and the \emph{order of decay} of
$A$.
\end{definition}

Thus $\tau(|A|^s)=\int e^{sx}\dd\nu_A$ is finite for $-p_-(A)<s<p_+(A)$ and
infinite for $s>p_+(A)$ and for $s<-p_-(A)$: the two orders are the reciprocal
abscissae of convergence of the two-sided Laplace transform of $\nu_A$, i.e.\ of
the Mellin transform of the spectral distribution of $|A|$. By
Proposition~\ref{tracial:prop:vn}(b),
\begin{equation}\label{infinite:eq:pc-orders}
p_c(\nu_A)=\min\{p_+(A),p_-(A)\},\qquad\frac1{p_c(\nu_A)}=\max\{\omega_+(A),\omega_-(A)\}.
\end{equation}
Bounded operators have $\omega_+=0$, operators with bounded inverse have
$\omega_-=0$, and $\omega_+(A^{-1})=\omega_-(A)$, since $\nu_{A^{-1}}$ is the
reflection of $\nu_A$.

\begin{proposition}\label{infinite:prop:orders}
Let $A\in\widetilde\cM$ be injective.
\begin{enumerate}
\item One has
\[
\begin{aligned}
p_+(A)&=\liminf_{x\to\infty}\frac{-\ln\nu_A((x,\infty))}{x},\\
p_-(A)&=\liminf_{x\to\infty}\frac{-\ln\nu_A((-\infty,-x))}{x},
\end{aligned}
\]
where $-\ln0=+\infty$.
\item For $\lambda\in\C\setminus\{0\}$, unitaries $U,V\in\cM$ and $s>0$:
$p_\pm(\lambda UAV)=p_\pm(A)$, $p_\pm(A^{-1})=p_\mp(A)$ and
$p_\pm(|A|^s)=p_\pm(A)/s$.
\end{enumerate}
\end{proposition}

\begin{proof}
(a) This is the classical expression of the abscissa of convergence of a
Laplace--Stieltjes transform through the decay of the tail \cite[Ch.~II]{Widder41};
we include the short argument. Write $\nu=\nu_A$, $M(s)=\int e^{sx}\dd\nu$ and
let $\sigma$ be the lower limit in the first formula. If $0<s<\sigma$, choose
$s<s'<\sigma$ and $x_0>0$ with $\nu((x,\infty))\le e^{-s'x}$ for $x\ge x_0$. Then
\[
M(s)\le e^{sx_0}+\int_{(x_0,\infty)}e^{sx}\dd\nu(x)
=e^{sx_0}+e^{sx_0}\nu((x_0,\infty))+\int_{x_0}^\infty se^{sy}\nu((y,\infty))\,\dd y<\infty ,
\]
so $p_+(A)\ge\sigma$. Conversely, if $M(s)<\infty$, then
$e^{sx}\nu((x,\infty))\le M(s)$ for $x>0$ by Chebyshev's inequality, hence
$-\frac1x\ln\nu((x,\infty))\ge s-\frac1x\ln M(s)$ and $\sigma\ge s$; thus
$\sigma\ge p_+(A)$. The second formula follows by applying the first to
$A^{-1}$.
(b) $|\lambda UAV|=|\lambda|V^*|A|V$ has the spectral distribution of
$|\lambda|\,|A|$, which translates $\nu_A$; $|A^{-1}|$ has the spectral
distribution of $|A|^{-1}$; and $\nu_{|A|^s}$ is the image of $\nu_A$ under
$x\mapsto sx$.
\end{proof}

With the distribution function $d_A(\lambda)=\tau(\one_{(\lambda,\infty)}(|A|))$,
part (a) says that
\[
p_+(A)=\liminf_{\lambda\to\infty}\frac{\ln(1/d_A(\lambda))}{\ln\lambda}.
\]
Inverting the distribution function one finds equivalently
\[
\omega_+(A)=\limsup_{t\to0}\frac{\ln\mu_t(A)}{\ln(1/t)},\qquad
\omega_-(A)=\limsup_{t\to0}\frac{\ln\mu_t(A^{-1})}{\ln(1/t)},
\]
in terms of the generalised singular values $\mu_t$ of \cite{FK86}; this
explains the names.

\subsection{Subadditivity with harmonic exponents}

By Remark~\ref{infinite:rem:nosub}, $C_p$ is not subadditive at a fixed exponent. It is
subadditive, with no additive correction term, if the exponents are combined as in
H\"older's inequality.

\begin{proposition}\label{infinite:prop:holder}
Let $A,B\in\widetilde\cM$ be injective and $p,q\in(0,\infty]$, and define $r$ by
$\frac1r=\frac1p+\frac1q$. Then
\[
C_r(AB)\le C_p(A)+C_q(B)\qquad\text{in }[0,\infty],
\]
with $C_\infty=\OLC$. Equality holds if $A\ge0$, $B=A^\theta$ for some $\theta>0$
and $q=p/\theta$; in particular if $A=B\ge0$ and $p=q$. Consequently
\[
\omega_\pm(AB)\le\omega_\pm(A)+\omega_\pm(B),\qquad
\frac1{p_c(\nu_{AB})}\le\frac1{p_c(\nu_A)}+\frac1{p_c(\nu_B)},
\]
and, for $p=q=\infty$, $\OLC(AB)\le\OLC(A)+\OLC(B)$.
\end{proposition}

\begin{proof}
For $\tau$-measurable $X,Y$ and $0<r<\infty$, the majorisation theorem of Fack
and Kosaki \cite[Theorem~4.2]{FK86}, applied to $f(t)=t^r$, gives
$\tau(|XY|^r)=\int_0^1\mu_t(XY)^r\dd t\le\int_0^1\mu_t(X)^r\mu_t(Y)^r\dd t$, and
the classical H\"older inequality on $(0,1)$ yields
$\|XY\|_r\le\|X\|_p\|Y\|_q$ for all $0<p,q\le\infty$ (see also \cite{PX03});
for $p=q=\infty$ this is $\|XY\|\le\|X\|\|Y\|$.
Here $AB$ is the strong product, which is injective with
$(AB)^{-1}=B^{-1}A^{-1}$ in $\widetilde\cM$. Applying the inequality to $AB$
and to $B^{-1}A^{-1}$,
\begin{align*}
2C_r(AB)&=\ln\|AB\|_r+\ln\|B^{-1}A^{-1}\|_r\\
&\le\ln\|A\|_p+\ln\|B\|_q+\ln\|B^{-1}\|_q+\ln\|A^{-1}\|_p=2C_p(A)+2C_q(B).
\end{align*}
If $A\ge0$, $B=A^\theta$ and $q=p/\theta$, then $AB=A^{1+\theta}$ and
$r=p/(1+\theta)$, and by Proposition~\ref{tracial:prop:vn}(e) both sides equal
$(1+\theta)C_p(A)$. For the orders, if $p<p_+(A)$ and $q<p_+(B)$ then
$\tau(|AB|^r)<\infty$, so $p_+(AB)\ge r$, that is,
$\omega_+(AB)\le\frac1p+\frac1q$; letting $p\uparrow p_+(A)$, $q\uparrow p_+(B)$
gives the inequality for $\omega_+$. The one for $\omega_-$ follows by applying it
to $B^{-1}A^{-1}$, and the last one from \eqref{infinite:eq:pc-orders}.
\end{proof}

The bound on the orders is the contrast form of the classical inclusion
$L^p\cdot L^q\subset L^r$. Thus $\omega_\pm$ and $\max\{\omega_+,\omega_-\}=1/p_c$
have the formal properties of $\OLC$ recorded in \cite{NN26}: they are
subadditive under products and homogeneous under powers of positive operators
(Proposition~\ref{infinite:prop:orders}(b)), $\max\{\omega_+,\omega_-\}$ is invariant
under $A\mapsto A^{-1}$, and $\omega_+$ and $\omega_-$ both vanish, i.e.\
$\max\{\omega_+,\omega_-\}=0$, exactly when $p_c=\infty$, in particular on
invertible elements of $\cM$. For all injective $A,B$, the case
$p=q=2r$ gives $C_r(AB)\le C_{2r}(A)+C_{2r}(B)$.

\subsection{The transition at the threshold}\label{infinite:subsec:transition}

\begin{proposition}\label{infinite:prop:cfun}
Let $A\in\widetilde\cM$ be injective and $p_c=p_c(\nu_A)$.
\begin{enumerate}
\item The function $p\mapsto C_p(A)$ is finite and real-analytic on $(0,p_c)$,
equal to $+\infty$ on $(p_c,\infty)$, nondecreasing on $(0,\infty]$, and
left-continuous at every $p\in(0,\infty]$ as a $[0,\infty]$-valued function. If
$0<p_c<\infty$, the restriction of $p\mapsto C_p(A)$ to $(0,p_c)$ has no
real-analytic extension to any interval $(p_c-\varepsilon,p_c+\varepsilon)$.
\item If $0<p_c<\infty$, exactly one of the following holds:
\begin{enumerate}[label=\textup{(\roman*)}]
\item \emph{blow-up}: $C_{p_c}(A)=\infty$, and $C_p(A)\uparrow\infty$ as $p\uparrow p_c$;
\item \emph{jump}: $C_{p_c}(A)<\infty$, while $C_p(A)=\infty$ for every $p>p_c$.
\end{enumerate}
The jump occurs if and only if $|A|$ and $|A^{-1}|$ both belong to
$L^{p_c}(\cM,\tau)$.
\item In the case of a jump, let $k\ge1$. The function $p\mapsto C_p(A)$ has
finite left derivatives of all orders $\le k$ at $p_c$ if and only if
$\int|x|^ke^{p_c|x|}\dd\nu_A(x)<\infty$. If $k_0$ is the smallest $k$ for which
this integral is infinite, then $\frac{\dd^{k_0}}{\dd p^{k_0}}C_p(A)\to+\infty$ as
$p\uparrow p_c$, and the left derivative of order $k_0$ at $p_c$ is $+\infty$.
\end{enumerate}
\end{proposition}

\begin{proof}
Write $\nu=\nu_A$, $G_\pm(p)=\int e^{\pm px}\dd\nu$ and $F=\ln G_++\ln G_-$, so that
$C_p(A)=F(p)/(2p)$.

(a) Finiteness on $(0,p_c)$, infiniteness on $(p_c,\infty)$ and monotonicity
are Proposition~\ref{tracial:prop:vn}(b),(c) and the definition of $p_c$. The two-sided
Laplace transform $s\mapsto\int e^{sx}\dd\nu$ is holomorphic in the strip
$-p_-(A)<\operatorname{Re}s<p_+(A)$ and positive on its real points
\cite[Ch.~VI]{Widder41}, so $F$ is real-analytic on $(0,p_c)$. Moreover
$e^{F(p)}=\int e^{px}\dd\nu\int e^{-py}\dd\nu=\int e^{pz}\dd(\nu*\check\nu)(z)$,
where $\check\nu$ is the reflection of $\nu$, is the Laplace transform of the
positive measure $\nu*\check\nu$, and its abscissa of convergence is $p_c$. By
Landau's theorem \cite[Ch.~II, Thm.~5b]{Widder41}, the real point of the
abscissa of convergence of the Laplace transform of a positive measure is a
singular point (apply it to the restriction of $\nu*\check\nu$ to $[0,\infty)$;
the transform of the restriction to $(-\infty,0)$ is holomorphic near $p_c$);
hence $e^F$, and with it $C_p=F/(2p)$, has no analytic
continuation to a neighbourhood of $p_c$. For
left-continuity at $p_0\in(0,\infty)$ let $p\uparrow p_0$: on $\{x>0\}$ the
integrand $e^{px}$ increases to $e^{p_0x}$, and on $\{x\le0\}$ it is bounded by
$1$ and converges; by monotone and dominated convergence
$G_+(p)\to G_+(p_0)\in(0,\infty]$, and likewise for $G_-$. At $p_0=\infty$ this
is Proposition~\ref{tracial:prop:vn}(c).

(b) This is (a) at $p=p_c$ together with Proposition~\ref{tracial:prop:vn}(b):
$C_{p_c}(A)<\infty$ iff $\int e^{p_c|x|}\dd\nu<\infty$ iff
$\tau(|A|^{p_c})<\infty$ and $\tau(|A^{-1}|^{p_c})<\infty$.

(c) Assume $C_{p_c}(A)<\infty$ and put $I=[p_c/2,p_c]$. For $p\in I$, $j\ge0$ and
all $x$, $|x|^je^{\pm px}\le|x|^j(1+e^{p_c|x|})$, and $\int|x|^j\dd\nu<\infty$
because $e^{p_c|x|/2}$ is integrable. Hence, if
$\int|x|^ke^{p_c|x|}\dd\nu<\infty$, dominated convergence shows that
$G_\pm^{(j)}(p)=\int(\pm x)^je^{\pm px}\dd\nu$ extend continuously to $I$ for
$j\le k$; since $G_\pm>0$ on $I$, so do the derivatives of $F$ and of
$F(p)/(2p)$ up to order $k$, and by the mean value theorem these limits are the
left derivatives at $p_c$. Conversely, let $k_0$ be the smallest order with
$\int|x|^{k_0}e^{p_c|x|}\dd\nu=\infty$. By the first part applied with $k_0-1$,
the derivatives of order $<k_0$ extend continuously to $I$. Split
\[
G_+^{(k_0)}(p)=\int_{x>0}x^{k_0}e^{px}\dd\nu+\int_{x<0}x^{k_0}e^{px}\dd\nu,
\]
\[
G_-^{(k_0)}(p)=\int_{x<0}|x|^{k_0}e^{p|x|}\dd\nu+\int_{x>0}(-x)^{k_0}e^{-px}\dd\nu .
\]
The second integrals are bounded on $I$ by $\int|x|^{k_0}\dd\nu$. The first ones
are nonnegative and nondecreasing in $p$, and by monotone convergence their sum
tends, as $p\uparrow p_c$, to $\int|x|^{k_0}e^{p_c|x|}\dd\nu=\infty$. Now
$F^{(k_0)}=G_+^{(k_0)}/G_++G_-^{(k_0)}/G_-+R$, where $R$ is a polynomial in
$G_\pm^{(j)}/G_\pm$, $j<k_0$, and therefore bounded on $I$; since $G_\pm$ stay
between positive constants on $I$, $F^{(k_0)}\to+\infty$. Further
$\frac{\dd^{k_0}}{\dd p^{k_0}}\frac{F(p)}{2p}=\frac{F^{(k_0)}(p)}{2p}+\sum_{j<k_0}c_{j}\,p^{-(k_0-j)-1}F^{(j)}(p)$
with bounded second term, so it tends to $+\infty$. Finally, the derivative of
order $k_0-1$ extends continuously to $p_c$, the extension is its left
derivative there, and by the mean value theorem its left difference quotients
at $p_c$ tend to $+\infty$.
\end{proof}

So, when $0<p_c<\infty$, the passage from finite to infinite contrast happens at
a single point, which is always a singularity, and comes in a graded scale: a blow-up, a jump with
infinite slope, a jump with finite slope but infinite curvature, and so on. This
is the behaviour of a power series with nonnegative coefficients, such as
$\sum_nn^{-\delta}z^n$, at the positive point of its circle of convergence. In
terms of noncommutative $L^p$-spaces: a jump means that $|A|$ and $|A^{-1}|$ lie
in $L^{p_c}$, while for every $p>p_c$ at least one of them is not in $L^p$; a
blow-up means that both lie in every $L^p$ with $p<p_c$, while at least one of
them is not in $L^{p_c}$. The relation between tail decay and singular
behaviour at the abscissa is studied in detail for Laplace--Stieltjes
transforms in \cite{Nakagawa}.

\begin{example}[The scale of transitions]\label{infinite:ex:transitions}
Let $\delta\ge0$ and let $A\ge1$ have the logarithmic spectral measure
$\nu=\nu_A$ on $[0,\infty)$ with $\nu((x,\infty))=e^{-x}(1+x)^{-\delta}$, $x\ge0$
(such $A$ exists in every diffuse finite von Neumann algebra). Then $A^{-1}$ is
bounded, $p_-(A)=\infty$, and $p_c=p_+(A)=1$ by
Proposition~\ref{infinite:prop:orders}(a), so $\omega_+(A)=1$ and $\omega_-(A)=0$.
Integrating by parts,
\[
\int e^{px}\dd\nu=1+p\int_0^\infty e^{-(1-p)x}(1+x)^{-\delta}\dd x ,
\]
\[
1\ge\int e^{-px}\dd\nu\ge\int e^{-x}\dd\nu>0\quad(0<p\le1).
\]
\begin{enumerate}
\item If $0\le\delta<1$, the substitution $x=y/(1-p)$ and monotone convergence
give $\int_0^\infty e^{-(1-p)x}(1+x)^{-\delta}\dd x\sim\Gamma(1-\delta)(1-p)^{\delta-1}$,
so there is a blow-up with
$C_p(A)=\frac{1-\delta}2\ln\frac1{1-p}+O(1)$ as $p\uparrow1$.
\item If $\delta=1$, the integral is $\ln\frac1{1-p}+O(1)$, and
$C_p(A)=\frac12\ln\ln\frac1{1-p}+O(1)$: still a blow-up, but a slower one.
\item If $\delta>1$, there is a jump, with
$C_1(A)=\frac12\ln\frac\delta{\delta-1}+\frac12\ln\int e^{-x}\dd\nu<\infty$.
Since, integrating by parts,
$\int x^ke^x\dd\nu=\int_0^\infty(kx^{k-1}+x^k)(1+x)^{-\delta}\dd x$ for $k\ge1$,
Proposition~\ref{infinite:prop:cfun}(c) shows that the left derivatives of order $k$
at $p_c=1$ are finite exactly for $k<\delta-1$.
\end{enumerate}
In terms of generalised singular values,
$\mu_t(A)\sim t^{-1}(\ln\frac1t)^{-\delta}$ as $t\to0$: the power $t^{-1}$ fixes
the threshold, and the logarithmic correction decides between blow-up and jump
and grades the jumps. In Example~\ref{infinite:ex:MP} with $y=1$ the factor
$\Gamma(1-2p)$ produces a blow-up at $p_c=\frac12$ with
$C_p(A)=\ln\frac1{\frac12-p}+O(1)$.
\end{example}

\begin{remark}\label{infinite:rem:invariants}
Proposition~\ref{infinite:prop:cfun} and Example~\ref{infinite:ex:transitions} single out the
following numerical characteristics of an operator with infinite contrast.
They depend only on $\nu_A$, hence are unitary invariants of $A$, and are
defined through traces of powers, without any truncation: the threshold $p_c$
and the orders $\omega_\pm$; the type of the transition at $p_c$; in the case of
a jump the value $C_{p_c}(A)$ and the number of finite left derivatives at
$p_c$; in the case of a blow-up its order
$\limsup_{p\uparrow p_c}2p_c\,C_p(A)/\ln\frac1{p_c-p}\in[0,\infty]$, which equals
$1-\delta$ ($0\le\delta\le1$) in Example~\ref{infinite:ex:transitions} and $1$ in Example~\ref{infinite:ex:MP} with
$y=1$; and, when $p_c=\infty$ but $\OLC(A)=\infty$, the growth order
$\limsup_{p\to\infty}\ln C_p(A)/\ln p$, which equals $\frac1{q-1}$ under the
hypothesis of Theorem~\ref{infinite:thm:growth}(b).
\end{remark}

\section{Finite sections of unbounded operators: growth of the contrast in the dimension}\label{infinite:sec:sections}

In Section~\ref{infinite:sec:unbounded} the operator was unbounded but affiliated with a
finite algebra. We now consider an unbounded operator on a Hilbert space, where
the Schatten contrast with respect to the trace of $B(\cH)$ is not available
(already for bounded invertible operators on an infinite-dimensional $\cH$, since $|A^{-1}|$ is unitarily equivalent to $|A|^{-1}$ and $\Tr(|A|^p)+\Tr(|A^{-1}|^p)=\Tr(|A|^p+|A|^{-p})\ge2\Tr(1)=\infty$, so that one of the two positive traces is infinite and $C_p(A)=\infty$), and study the
contrasts of its finite sections as functions of their size $n$. Throughout,
$L$ is a positive self-adjoint operator on a separable Hilbert space with
compact resolvent, with eigenvalues
\[
0<\lambda_1\le\lambda_2\le\dots,\qquad\lambda_n\to\infty,
\]
counted with multiplicity, and an orthonormal basis of eigenvectors $(e_j)$. The
\emph{spectral section} $L_n$ is the restriction of $L$ to
$\operatorname{span}\{e_1,\dots,e_n\}$, a positive invertible element of
$M_n(\C)$ with the normalised trace. Thus
\begin{equation}\label{infinite:eq:OLC-section}
\OLC(L_n)=\tfrac12\ln\frac{\lambda_n}{\lambda_1},\qquad\kappa(L_n)=\frac{\lambda_n}{\lambda_1},
\end{equation}
and the growth of the operator-norm contrast is the growth of the eigenvalues.
The Schatten contrasts $C_p(L_n)$ behave differently, and their behaviour
separates spectra of power type, $\kappa(L_n)=O(n^\alpha)$, from lacunary
spectra, $\lambda_{j+1}/\lambda_j\ge\vartheta>1$, such as $\lambda_j=j!$ with
$\kappa(L_n)=n!$.

\subsection{Spectral sections are optimal}

We first show that \eqref{infinite:eq:OLC-section} is, up to an additive constant, the
slowest possible growth among all nested finite sections, so that the growth
class of $n\mapsto\OLC(L_n)$ is an invariant of $L$ and not of the choice of
sections. For a finite-dimensional subspace $V$ of the form domain of $L$ let
$L_V$ be the positive operator on $V$ defined by
$\langle L_Vu,v\rangle=\langle L^{1/2}u,L^{1/2}v\rangle$, $u,v\in V$.

\begin{proposition}\label{infinite:prop:minmax}
Let $V_1\subset V_2\subset\dots$ be subspaces of the form domain of $L$ with
$\dim V_n=n$, and let $\mu_1=\|L^{1/2}v\|^2/\|v\|^2$ for $0\ne v\in V_1$. Then
\[
\OLC(L_{V_n})\ \ge\ \tfrac12\ln\frac{\lambda_n}{\mu_1}\ =\ \OLC(L_n)-\tfrac12\ln\frac{\mu_1}{\lambda_1}\qquad(n\ge1).
\]
If $\bigcup_nV_n$ is a form core for $L$, then
$\liminf_{n\to\infty}\bigl(\OLC(L_{V_n})-\OLC(L_n)\bigr)\ge0$. There is no bound
in the opposite direction.
\end{proposition}

\begin{proof}
The largest eigenvalue of $L_{V_n}$ is the maximum of the Rayleigh quotient on
$V_n$, which is at least $\lambda_n$ by the min--max principle \cite{RS78}. The
smallest eigenvalue of $L_{V_n}$ is the minimum of the Rayleigh quotient on
$V_n$; it is nonincreasing in $n$, at most $\mu_1$, and tends to
$\lambda_1$ if $\bigcup_nV_n$ is a form core. For the last statement let
$V_n=\operatorname{span}\{e_{\sigma(1)},\dots,e_{\sigma(n)}\}$ for a permutation
$\sigma$ of $\N$: then
$\OLC(L_{V_n})=\frac12\ln\bigl(\max_{j\le n}\lambda_{\sigma(j)}/\min_{j\le n}\lambda_{\sigma(j)}\bigr)$,
which can grow arbitrarily fast.
\end{proof}

\subsection{Growth of the contrast and Schatten classes}

The growth rate of $\OLC(L_n)$ on the logarithmic scale is the same thing as
the abscissa of convergence of the spectral zeta function
$\zeta_L(r)=\Tr L^{-r}=\sum_j\lambda_j^{-r}$, that is, as the membership of
$L^{-1}$ in the Schatten classes $\mathcal S_r$.

\begin{theorem}\label{infinite:thm:abscissa}
Let $\gamma_*=\liminf_{n\to\infty}\dfrac{2\,\OLC(L_n)}{\ln n}\in[0,\infty]$. Then
\[
\inf\{r>0:\ L^{-1}\in\mathcal S_r\}=\inf\{r>0:\ \zeta_L(r)<\infty\}=\frac1{\gamma_*}
\]
(with $\frac10=\infty$ and $\frac1\infty=0$). In particular:
\begin{enumerate}
\item if $\kappa(L_n)=O(n^\alpha)$, then $L^{-1}\notin\mathcal S_r$ for
$r\le1/\alpha$;
\item if $\kappa(L_n)\ge c\,n^\alpha$ for some $c>0$ and all large $n$, then
$L^{-1}\in\mathcal S_r$ for all $r>1/\alpha$;
\item $L^{-1}\in\bigcap_{r>0}\mathcal S_r$ if and only if
$\OLC(L_n)/\ln n\to\infty$; this is the case for geometric and for factorial
growth of $\kappa(L_n)$.
\end{enumerate}
\end{theorem}

\begin{proof}
By \eqref{infinite:eq:OLC-section}, $\gamma_*=\liminf_n\ln\lambda_n/\ln n$. If
$r>1/\gamma_*$, choose $\gamma'<\gamma_*$ with $r\gamma'>1$; then
$\lambda_n\ge n^{\gamma'}$ for large $n$ and
$\sum_n\lambda_n^{-r}<\infty$. Conversely, if $\sum_j\lambda_j^{-r}=S<\infty$,
then, the sequence being nondecreasing, $n\lambda_n^{-r}\le S$, so
$\lambda_n\ge(n/S)^{1/r}$ and $\gamma_*\ge1/r$. The singular values of $L^{-1}$
are $\lambda_j^{-1}$, which gives the statement on Schatten classes, and
(b), (c) are special cases. In (a), $\lambda_n\le Cn^\alpha$ gives
$\sum_n\lambda_n^{-r}\ge C^{-r}\sum_nn^{-\alpha r}=\infty$ for $\alpha r\le1$.
\end{proof}

Theorem~\ref{infinite:thm:abscissa} is a statement on the logarithmic scale. The
relation between the contrasts of the sections and the Schatten ideals is in
fact exact, and it is convenient to state it for compact operators. Let $X$ be
a compact operator with singular values $s_1\ge s_2\ge\dots>0$, and let $X_n$
be the section of $|X|$ on the span of the first $n$ singular vectors, so that
$\kappa(X_n)=s_1/s_n$ and $\OLC(X_n)=\frac12\ln(s_1/s_n)$. Recall that the weak
Schatten ideal $\mathcal S_{p,\infty}$ consists of the compact operators with
$s_n=O(n^{-1/p})$, and that
$\mathcal S_p\subset\mathcal S_{p,\infty}\subset\mathcal S_q$ for $q>p$. In this
section $\|X\|_p=(\Tr|X|^p)^{1/p}$ denotes the unnormalised Schatten
(quasi-)norm on $B(\cH)$, whereas the contrasts $C_q$ of the sections are
always computed with the normalised trace of $M_n(\C)$.

\begin{proposition}\label{infinite:prop:ideals}
Let $X$ be as above and $0<p<\infty$.
\begin{enumerate}
\item $\displaystyle\Bigl(\frac{\|X\|_p}{\|X\|}\Bigr)^p=\sum_{n\ge1}\kappa(X_n)^{-p}=\sum_{n\ge1}e^{-2p\,\OLC(X_n)}$.
In particular $X\in\mathcal S_p$ if and only if this series converges.
\item If $X\in\mathcal S_p$, then for all $n$
\[
\OLC(X_n)\ \ge\ \frac1{2p}\ln n-\frac12\ln\frac{\|X\|_p}{\|X\|},
\]
and $\OLC(X_n)-\frac1{2p}\ln n\to+\infty$.
\item $X\in\mathcal S_{p,\infty}$ if and only if
$\inf_n\bigl(\OLC(X_n)-\frac1{2p}\ln n\bigr)>-\infty$.
\end{enumerate}
For an unbounded positive operator $L$ with compact resolvent and $X=L^{-1}$
one has $\kappa(X_n)=\kappa(L_n)$; hence $L^{-1}\in\mathcal S_p$ if and only if
$\sum_n\kappa(L_n)^{-p}<\infty$, with
$\zeta_L(p)=\lambda_1^{-p}\sum_n\kappa(L_n)^{-p}$, and
$L^{-1}\in\mathcal S_{p,\infty}$ if and only if
$\OLC(L_n)\ge\frac1{2p}\ln n-C$ for some constant $C$.
\end{proposition}

\begin{proof}
(a) $s_n=s_1/\kappa(X_n)$ and $\|X\|=s_1$. (b) Since the $s_j$ are
nonincreasing, $n\,s_n^p\le\sum_{j\le n}s_j^p\le\|X\|_p^p$, so
$\kappa(X_n)=s_1/s_n\ge n^{1/p}\|X\|/\|X\|_p$; moreover
$\frac n2s_n^p\le\sum_{n/2<j\le n}s_j^p\to0$, so $n^{1/p}s_n\to0$.
(c) $s_n\le Cn^{-1/p}$ is equivalent to
$\frac12\ln(s_1/s_n)\ge\frac1{2p}\ln n-\frac12\ln(C/s_1)$.
\end{proof}

Thus a summable sequence of singular values forces the contrast of the sections
to grow at least like $\frac12\ln n$, a $p$-summable one at least like
$\frac1{2p}\ln n$, and the unbounded operators whose inverses lie in
$\mathcal S_p$ are exactly those with $\sum_n\kappa(L_n)^{-p}<\infty$. Note the
direction of the inequalities: membership in an ideal is a \emph{lower} bound
for the growth of the contrast, and it is the weak ideal
$\mathcal S_{p,\infty}$, not $\mathcal S_p$, that corresponds exactly to the
bound $\frac1{2p}\ln n-C$. In Connes' terminology \cite{Connes94}, $X$ is an
infinitesimal of order $\alpha$ if $s_n=O(n^{-\alpha})$, i.e.\ if
$\OLC(X_n)\ge\frac\alpha2\ln n-C$: the slope of the contrast is half the order.

A dichotomy of the same kind appears at the abscissa. If $\sigma=1/\gamma_*$
lies in $(0,\infty)$, then by Proposition~\ref{infinite:prop:ideals}(a) and
Theorem~\ref{infinite:thm:abscissa} either $\sum_n\kappa(L_n)^{-\sigma}<\infty$, i.e.\
$L^{-1}\in\mathcal S_\sigma$, or the series diverges and $L^{-1}$ lies in
$\mathcal S_r$ only for $r>\sigma$. For $\lambda_n\sim c\,n^\gamma$ one has
$\kappa(L_n)^{-1/\gamma}\sim(\lambda_1/c)^{1/\gamma}n^{-1}$, so the series diverges
like the harmonic series; for $\lambda_n=n^\gamma(1+\ln n)^{2\gamma}$ the abscissa
is still $1/\gamma$, but
$\sum_n\kappa(L_n)^{-1/\gamma}=\sum_n\frac1{n(1+\ln n)^2}<\infty$. Unlike in
Proposition~\ref{infinite:prop:cfun}(b), the two cases are not separated by the
contrasts of the sections themselves: $C_\sigma(L_n)\to\infty$ in the
convergent case by Proposition~\ref{infinite:prop:pq}(b), and for
$\lambda_n\sim c\,n^\gamma$ by Theorem~\ref{infinite:thm:weyl}(c) (in the last example
$C_{1/\gamma}(L_n)=\gamma\ln\ln n+O(1)$, compared with $\frac\gamma2(\ln\ln n-\ln2)+o(1)$ in
Theorem~\ref{infinite:thm:weyl}(c) below).

The same comparison can be made between an $\ell_p$ condition on the spectrum
and a contrast with a different exponent $q$. The outcome is a threshold at
$q=p$.

\begin{proposition}\label{infinite:prop:pq}
Let $X\in\mathcal S_p$, $0<p<\infty$, have singular values $s_j>0$, and let
$X_n$ be its sections.
\begin{enumerate}
\item If $q>p$, then for all $n$
\[
C_q(X_n)\ \ge\ \frac12\Bigl(\frac1p-\frac1q\Bigr)\ln n-\frac12\ln\frac{\|X\|_p}{\|X\|}-\frac1{2q}\ln\Bigl(1+\frac qp\Bigr).
\]
\item $C_p(X_n)\to\infty$ as $n\to\infty$.
\item For $q<p$ nothing follows: $C_q(X_n)$ may stay bounded or tend to
infinity.
\end{enumerate}
Consequently, for any compact $X$ with $s_j>0$ for all $j$: if
$\sup_nC_q(X_n)<\infty$, then $X\notin\mathcal S_r$ for every $r\le q$, and
$\kappa(X_n)=O(n^{2/q})$.
\end{proposition}

\begin{proof}
(a) We have $\sum_{j\le n}s_j^q\ge s_1^q=\|X\|^q$, and, as in the proof of
Proposition~\ref{infinite:prop:ideals}, $s_j\le\|X\|_p\,j^{-1/p}$, so
\[
\sum_{j\le n}s_j^{-q}\ \ge\ \|X\|_p^{-q}\sum_{j\le n}j^{q/p}\ \ge\ \|X\|_p^{-q}\,\frac{n^{q/p+1}}{1+q/p}.
\]
Insert this into
$2q\,C_q(X_n)=\ln\sum_{j\le n}s_j^q+\ln\sum_{j\le n}s_j^{-q}-2\ln n$.

(b) Let $\varepsilon_n=\sup_{k>n/2}k\,s_k^p$, which tends to $0$. Then
$\sum_{n/2<j\le n}s_j^{-p}\ge\frac n2\cdot\frac n{2\varepsilon_n}$, hence
$e^{2p\,C_p(X_n)}\ge s_1^p/(4\varepsilon_n)\to\infty$.

(c) For $s_j=j^{-a}$ with $a>1/p$, Theorem~\ref{infinite:thm:weyl} below (applied to
$\lambda_j=j^a$, since $C_q$ is invariant under inversion) shows that
$C_q(X_n)$ is bounded for $q<1/a$ and unbounded for $q\ge1/a$, and $1/a<p$.

The last statement is the contrapositive of (a), (b), together with
$C_q(X_n)\ge\OLC(X_n)-\frac1q\ln n$, which follows from
$\frac1n\sum_{j\le n}s_j^q\ge\frac{s_1^q}n$ and $\frac1n\sum_{j\le n}s_j^{-q}\ge\frac{s_n^{-q}}n$.
\end{proof}

Thus an $\ell_p$ spectrum forces every contrast with $q\ge p$ to diverge, at a
rate with slope at least $\frac12(\frac1p-\frac1q)$, which increases to the
slope $\frac1{2p}$ of $\OLC$ as $q\to\infty$, and it imposes no condition on the
contrasts with $q<p$. The slope in (a) is optimal: for $s_j=j^{-a}$ the exact
slope is $\frac12(a-\frac1q)$ by Theorem~\ref{infinite:thm:weyl}, and $a$ can be taken
arbitrarily close to $1/p$. In terms of the two exponents
\[
q_c(X)=\sup\bigl\{q:\ \sup_nC_q(X_n)<\infty\bigr\},\qquad
p_0(X)=\inf\{p:\ X\in\mathcal S_p\},
\]
Proposition~\ref{infinite:prop:pq} says that $q_c(X)\le p_0(X)$. Equality holds for
power-type spectra, where both are equal to $1/a$. It can fail for irregular
spectra: if $s$ takes the value $e^{-k!}$ on a block of length
$N_k=\lceil e^{b\,k\cdot k!}\rceil$, $k\ge1$, with some $b>0$, then
$\sum_kN_ke^{-p\,k!}=\infty$ for every $p$, so $p_0(X)=\infty$, whereas for
$n=\sum_{k\le K}N_k+1$, the first index of the $(K+1)$-st block, one has
$n\le2KN_K$ and
$e^{2q\,C_q(X_n)}\ge N_Ke^{q(K+1)!-qK!}/n^2\ge e^{(q-b)K\cdot K!}/(8K^2)\to\infty$
for $q>b$, so $q_c(X)\le b$. Long plateaus followed by large drops make the contrast
unbounded without making the spectrum summable.

\subsection{Power-type spectra: a phase transition at the abscissa}

For operators obeying a Weyl law the Schatten contrasts of the sections have a
phase transition in $p$, located at the abscissa $1/\gamma$ of the zeta
function. Below it they converge to the contrast of an unbounded operator in
the sense of Section~\ref{infinite:sec:unbounded}; above it they grow logarithmically.

\begin{theorem}\label{infinite:thm:weyl}
Suppose that $\lambda_n=c\,n^\gamma(1+o(1))$ as $n\to\infty$, with $c,\gamma>0$.
Then, as $n\to\infty$:
\begin{enumerate}
\item\label{infinite:w:olc} $\OLC(L_n)=\dfrac\gamma2\ln n+\dfrac12\ln\dfrac c{\lambda_1}+o(1)$.
\item\label{infinite:w:sub} For $0<p<1/\gamma$,
\[
C_p(L_n)\longrightarrow-\frac1{2p}\ln\bigl(1-\gamma^2p^2\bigr)=C_p(\nu_\gamma),
\]
where $\nu_\gamma$ is the law of $\gamma\ln U$ for $U$ uniformly distributed on
$[0,1]$. Moreover the logarithmic spectral measures of $L_n/\lambda_n$ converge
weakly to $\nu_\gamma$, which is the logarithmic spectral measure of the
multiplication operator by $x^\gamma$ on $L^2[0,1]$, and $p_c(\nu_\gamma)=1/\gamma$.
\item\label{infinite:w:crit} For $p=1/\gamma$,
$C_p(L_n)=\dfrac\gamma2\bigl(\ln\ln n-\ln2\bigr)+o(1)$.
\item\label{infinite:w:super} For $p>1/\gamma$,
\[
C_p(L_n)=\Bigl(\frac\gamma2-\frac1{2p}\Bigr)\ln n+\frac1{2p}\ln\frac{c^p\,\zeta_L(p)}{1+\gamma p}+o(1),
\]
and consequently
$\displaystyle\OLC(L_n)-C_p(L_n)=\frac1{2p}\ln n-\frac1{2p}\ln\frac{\lambda_1^p\,\zeta_L(p)}{1+\gamma p}+o(1)$.
\end{enumerate}
\end{theorem}

\begin{proof}
We use the following elementary fact: if $a_j>0$ and $a_j=j^\beta(1+o(1))$,
then $\sum_{j\le n}a_j$ is asymptotic to $n^{\beta+1}/(\beta+1)$ if
$\beta>-1$ and to $\ln n$ if $\beta=-1$, and converges if $\beta<-1$.

\ref{infinite:w:olc} is \eqref{infinite:eq:OLC-section}. For the other parts,
$2p\,C_p(L_n)=\ln\bigl(\frac1n\sum_{j\le n}\lambda_j^{p}\bigr)+\ln\bigl(\frac1n\sum_{j\le n}\lambda_j^{-p}\bigr)$.
Since $\lambda_j^{p}=c^pj^{\gamma p}(1+o(1))$,
\[
\frac1n\sum_{j\le n}\lambda_j^{p}=\frac{c^p\,n^{\gamma p}}{1+\gamma p}(1+o(1)),
\]
while $\frac1n\sum_{j\le n}\lambda_j^{-p}$ is asymptotic to
$c^{-p}n^{-\gamma p}/(1-\gamma p)$ if $\gamma p<1$, to $c^{-p}\ln n/n$ if
$\gamma p=1$, and to $\zeta_L(p)/n$ if $\gamma p>1$. Multiplying gives
$\bigl((1+\gamma p)(1-\gamma p)\bigr)^{-1}$, $\frac12\ln n$ and
$c^p\zeta_L(p)\,n^{\gamma p-1}/(1+\gamma p)$ respectively, which is
\ref{infinite:w:sub}--\ref{infinite:w:super}; the last formula in \ref{infinite:w:super} follows by
subtracting from \ref{infinite:w:olc}. In \ref{infinite:w:sub},
$\int e^{\pm px}\dd\nu_\gamma=\int_0^1u^{\pm\gamma p}\dd u=(1\pm\gamma p)^{-1}$
for $\gamma p<1$, and the integral with the minus sign diverges for
$\gamma p\ge1$. Finally, for $u\in(0,1]$ the proportion of indices $j\le n$ with
$\lambda_j/\lambda_n\le u^\gamma$ tends to $u$, which is the weak convergence.
\end{proof}

Thus, for a power-type spectrum, the sections $L_n$, rescaled by their norm,
converge in distribution to a bounded operator with unbounded inverse in a
finite von Neumann algebra, whose critical exponent is the abscissa of the
spectral zeta function of $L$. For $p$ below the abscissa the Schatten contrast
of the sections stays bounded although $\kappa(L_n)\to\infty$; for $p$ above
it, it grows like a multiple of $\ln n$ with a slope
$\frac\gamma2-\frac1{2p}$ that increases to the slope $\frac\gamma2$ of $\OLC$.

\subsection{Lacunary spectra}

\begin{theorem}\label{infinite:thm:lacunary}
Suppose that $\lambda_{j+1}/\lambda_j\ge\vartheta$ for all $j$, with some $\vartheta>1$. Then
for every $p\in(0,\infty)$ and every $n$,
\[
\OLC(L_n)-\frac{\ln n}p\ \le\ C_p(L_n)\ \le\ \OLC(L_n)-\frac{\ln n}p+\frac1p\ln\frac1{1-\vartheta^{-p}} .
\]
In particular $C_p(L_n)\to\infty$ for every $p>0$.
\end{theorem}

\begin{proof}
We have
$2p\,C_p(L_n)=p\ln\frac{\lambda_n}{\lambda_1}-2\ln n+\ln\sum_{j\le n}\bigl(\frac{\lambda_j}{\lambda_n}\bigr)^p+\ln\sum_{j\le n}\bigl(\frac{\lambda_1}{\lambda_j}\bigr)^p$,
and each of the two sums lies between $1$ and
$\sum_{k\ge0}\vartheta^{-kp}=(1-\vartheta^{-p})^{-1}$. Finally,
$\OLC(L_n)\ge\frac{n-1}2\ln\vartheta$.
\end{proof}

Theorems~\ref{infinite:thm:weyl} and \ref{infinite:thm:lacunary} separate the two growth classes
by the coefficient of $\ln n$ in the gap between the contrasts:
\begin{equation}\label{infinite:eq:classes}
\OLC(L_n)-C_p(L_n)=
\begin{cases}
\dfrac1{2p}\ln n+O(1),&\lambda_n\sim c\,n^\gamma,\ p>1/\gamma,\\[2mm]
\dfrac1{p}\ln n+O(1),&\lambda_{j+1}/\lambda_j\ge\vartheta>1 .
\end{cases}
\end{equation}
Heuristically: if the largest and the smallest eigenvalue of a positive invertible
$B\in M_n(\C)$ have multiplicities $m_1$ and $m_r$, then
$\frac1n\Tr B^{p}=\frac{m_1}n\|B\|^{p}(1+o(1))$ and
$\frac1n\Tr B^{-p}=\frac{m_r}n\|B^{-1}\|^{p}(1+o(1))$ as $p\to\infty$, so $C_p(B)=\OLC(B)-\frac1{2p}\ln\frac{n^2}{m_1m_r}+o(\frac1p)$. For lacunary spectra both edges of the
spectrum of $L_n$ are isolated and the correction is the full
$\frac1{2p}\ln\frac{n^2}{m_1m_r}=\frac{\ln n}p$; for power growth the upper edge
is not isolated, the top of the spectrum behaves like a continuum, and only the
lower edge contributes.

\begin{example}\label{infinite:ex:sections}
(i) \emph{The Dirichlet Laplacian on $[0,\pi]$.} Here $\lambda_j=j^2$, so
$\gamma=2$, $c=1$ and $\zeta_L(p)=\zeta(2p)$. Then $\OLC(L_n)=\ln n$ exactly;
$C_p(L_n)\to-\frac1{2p}\ln(1-4p^2)$ for $p<\frac12$; and, for instance,
$C_1(L_n)=\frac12\ln n+\frac12\ln\frac{\pi^2}{18}+o(1)$.

(ii) \emph{Elliptic operators.} For the Dirichlet Laplacian on a bounded domain
in $\R^d$, and more generally for a positive self-adjoint elliptic operator of order $m$ on
a compact $d$-dimensional manifold without boundary, Weyl's law
(\cite{Weyl12} and \cite{RS78} for domains, \cite{Hormander68} for manifolds) gives
$\lambda_n\sim c\,n^{\gamma}$ with $\gamma=m/d$. Hence
$\OLC(L_n)=\frac m{2d}\ln n+O(1)$, i.e.\ $\kappa(L_n)=O(n^{m/d})$, and the
Schatten contrasts of the sections remain bounded exactly for $p<d/m$; for the
Laplacian the critical exponent is half the dimension.

(iii) \emph{The harmonic oscillator.} $\lambda_j=2j-1$, $\gamma=1$, $c=2$,
$\zeta_L(p)=(1-2^{-p})\zeta(p)$: $\OLC(L_n)=\frac12\ln(2n-1)$, the critical
exponent is $1$, and
$C_2(L_n)=\frac14\ln n+\frac14\ln\frac{\pi^2}{6}+o(1)$.

(iv) \emph{Geometric and factorial spectra.} For $\lambda_j=\vartheta^j$,
$\OLC(L_n)=\frac{n-1}2\ln\vartheta$ and $\kappa(L_n)=\vartheta^{n-1}$. For $\lambda_j=j!$,
$\kappa(L_n)=n!$ and, by Stirling's formula,
\[
\OLC(L_n)=\tfrac12\ln n!=\tfrac12\,n\ln n-\tfrac12\,n+\tfrac14\ln(2\pi n)+o(1).
\]
In both cases $C_p(L_n)=\OLC(L_n)-\frac{\ln n}p+O(1)$ for every $p$, there is no
critical exponent, and $L^{-1}$ belongs to every Schatten class.
\end{example}

\subsection{Upper and lower orders, and their contrast}\label{infinite:subsec:orders}

Theorem~\ref{infinite:thm:abscissa} uses only the lower limit of $2\,\OLC(L_n)/\ln n$. Put
\[
\gamma_*(L)=\liminf_{n\to\infty}\frac{2\,\OLC(L_n)}{\ln n}=\liminf_{n\to\infty}\frac{\ln\lambda_n}{\ln n},
\]
\[
\gamma^*(L)=\limsup_{n\to\infty}\frac{2\,\OLC(L_n)}{\ln n}=\limsup_{n\to\infty}\frac{\ln\lambda_n}{\ln n},
\]
both in $[0,\infty]$, and, when $0<\gamma_*\le\gamma^*<\infty$, define the
\emph{order contrast}
\[
\Delta_{\rm ord}(L)=\frac{\gamma^*-\gamma_*}{\gamma^*+\gamma_*}
=\OMC\bigl(\operatorname{diag}(\gamma_*,\gamma^*)\bigr),\qquad
\operatorname{artanh}\Delta_{\rm ord}(L)=\tfrac12\ln\frac{\gamma^*}{\gamma_*}.
\]
We set $\Delta_{\rm ord}(L)=1$ if $\gamma_*=0<\gamma^*$ or $\gamma_*<\gamma^*=\infty$,
and leave it undefined if $\gamma_*=\gamma^*\in\{0,\infty\}$. Thus
$\Delta_{\rm ord}(L)=0$ exactly when $\lambda_n=n^{\gamma+o(1)}$ for some
$\gamma\in(0,\infty)$. Let $N_L(\lambda)=\#\{n:\lambda_n\le\lambda\}$ be the eigenvalue
counting function.

\begin{proposition}\label{infinite:prop:ordcontrast}
Let $L$ be as above.
\begin{enumerate}
\item For $c,s>0$: $\gamma_*(cL^s)=s\,\gamma_*(L)$ and $\gamma^*(cL^s)=s\,\gamma^*(L)$; hence
$\Delta_{\rm ord}(cL^s)=\Delta_{\rm ord}(L)$.
\item With $\frac10=\infty$ and $\frac1\infty=0$,
\[
\frac1{\gamma_*(L)}=\limsup_{\lambda\to\infty}\frac{\ln N_L(\lambda)}{\ln\lambda},\qquad
\frac1{\gamma^*(L)}=\liminf_{\lambda\to\infty}\frac{\ln N_L(\lambda)}{\ln\lambda} ;
\]
in particular $\Delta_{\rm ord}$ is also the Michelson contrast of the upper and
lower orders of $N_L$.
\item The abscissa of convergence of $\zeta_L$ is $1/\gamma_*(L)$
(Theorem~\ref{infinite:thm:abscissa}) and does not determine $\gamma^*(L)$: for every pair
$0\le a\le b\le\infty$ there is $L$ with $\gamma_*(L)=a$ and $\gamma^*(L)=b$.
\item If $\sup_nC_q(L_n)<\infty$ for some $q>0$, then $q\le\min\{1/\gamma_*,\,2/\gamma^*\}$.
In particular the exponent $q_c=\sup\{q:\sup_nC_q(L_n)<\infty\}$ ($\sup\emptyset=0$,
$1/0=\infty$) satisfies $q_c\le\min\{1/\gamma_*,2/\gamma^*\}$, and, when $\Delta_{\rm ord}$
is defined, the second bound is the better one exactly when $\gamma^*>2\gamma_*$, i.e.\ when
$\Delta_{\rm ord}(L)>\frac13$.
\end{enumerate}
\end{proposition}

\begin{proof}
(a) The eigenvalues of $cL^s$ are $c\lambda_n^s$, in the same order, and
$\ln(c\lambda_n^s)/\ln n=s\ln\lambda_n/\ln n+\ln c/\ln n$.

(b) Since $N_L(\lambda_n)\ge n$, we have $\ln N_L(\lambda_n)/\ln\lambda_n\ge\ln n/\ln\lambda_n$,
so the upper limit in the first formula is at least
$\limsup_n\ln n/\ln\lambda_n=1/\gamma_*$. Conversely, if $N_L(\lambda)=m$ then
$\lambda_m\le\lambda$, so $\ln N_L(\lambda)/\ln\lambda\le\ln m/\ln\lambda_m$ once
$\lambda_m>1$, and $m\to\infty$ as $\lambda\to\infty$; hence the upper limit is at most
$1/\gamma_*$. For the second formula: if $\lambda_m\le\lambda<\lambda_{m+1}$ then
$N_L(\lambda)\ge m$ and $\ln N_L(\lambda)/\ln\lambda\ge\frac{\ln m}{\ln(m+1)}\cdot
\frac{\ln(m+1)}{\ln\lambda_{m+1}}$, so the lower limit is at least $1/\gamma^*$; and
for $\lambda<\lambda_n$ we have $N_L(\lambda)\le n-1$, so letting $\lambda\uparrow\lambda_n$
along a sequence with $\ln\lambda\ge\ln\lambda_n-\frac1n$ shows that it is at most
$\liminf_n\ln(n-1)/\ln\lambda_n=1/\gamma^*$.

(c) Let first $0<a\le b<\infty$. Choose $T_1\ge1$ and define a nondecreasing
function $\psi$ on $[0,\infty)$ by induction on blocks: on $[T_k,T_k']$, with
$T_k'=2T_k$, let $\psi(t)=bt$; on $[T_k',\frac baT_k']$ let $\psi(t)=bT_k'$; on
$[\frac baT_k',T_{k+1})$, with $T_{k+1}=k\,\frac baT_k'$, let $\psi(t)=at$; and $\psi=0$
on $[0,T_1)$. Then $a\le\psi(t)/t\le b$ for $t\ge T_1$, with $\psi(t)/t=a$ on the
third kind of interval and $=b$ on the first, so $\liminf\psi(t)/t=a$ and
$\limsup\psi(t)/t=b$. Put $\lambda_n=e^{\psi(\ln n)}+1$ and let $L$ be diagonal with
these eigenvalues; then $\ln\lambda_n/\ln n$ has lower limit $a$ and upper limit $b$.
If $a=0<b$, let $\psi(t)=bT_k'$ on $[T_k',T_{k+1})$ with $T_{k+1}=k^2T_k'$. If $b=\infty$,
use $\psi(t)=t^2$ on the first kind of interval (with $T_1\ge\max\{a,1\}$), continue
with the constant $\psi(T_k')=T_k'^2$ up to $t=T_k'^2/a$ and then with $\psi(t)=at$ up
to $T_{k+1}=kT_k'^2/a$ if $a>0$, and with the constant $T_k'^2$ up to
$T_{k+1}=k^2T_k'^2$ if $a=0$. The cases
$a=b\in\{0,\infty\}$ are realised by $\lambda_n=\ln(n+1)$ and $\lambda_n=e^n$.

(d) $C_q(L_n)$ is invariant under inversion, and $L_n^{-1}$ is the $n$-th section of
the compact operator $X=L^{-1}$ in the sense of Proposition~\ref{infinite:prop:ideals}. By the
last statement of Proposition~\ref{infinite:prop:pq}, $\sup_nC_q(X_n)<\infty$ implies
$X\notin\mathcal S_r$ for $r\le q$, i.e.\ $q\le1/\gamma_*$ by
Theorem~\ref{infinite:thm:abscissa}, and $\kappa(L_n)=\kappa(X_n)=O(n^{2/q})$, i.e.\
$\gamma^*\le2/q$.
\end{proof}

So the zeta function and the Schatten classes of $L^{-1}$ see only the lower
order, and the order contrast is a second, scale-free number: it is a property of
the family $\{cL^s\}$ rather than of $L$ (upper and lower limits of scales of spaces are
studied in \cite{EKN}), and its logarithmic form
$\frac12\ln(\gamma^*/\gamma_*)$ is the logarithmic contrast of the pair of orders.
When the orders coincide, finer information is in the constants
(Subsection~\ref{infinite:subsec:weylosc}); when they are degenerate ($\gamma_*=\gamma^*=0$, as
for $\lambda_n=\ln(n+1)$, or $=\infty$, as for $\lambda_n=e^n$), the same construction
can be repeated with $\ln\ln\lambda_n/\ln\ln n$ or $\ln\lambda_n/n$ in place of
$\ln\lambda_n/\ln n$.

\subsection{Local contrasts and relativistic addition}\label{infinite:subsec:local}

For $n\ge1$ let
\[
\delta_n(L)=\OMC\bigl(\operatorname{diag}(\lambda_n,\lambda_{n+1})\bigr)
=\frac{\lambda_{n+1}-\lambda_n}{\lambda_{n+1}+\lambda_n}\in[0,1)
\]
be the \emph{local contrasts} of $L$, and let
$x\oplus y=\frac{x+y}{1+xy}$ on $(-1,1)$ be the relativistic addition of velocities.

\begin{proposition}\label{infinite:prop:local}
\begin{enumerate}
\item $\OLC(L_n)=\sum_{k<n}\operatorname{artanh}\delta_k$ and
$\OMC(L_n)=\delta_1\oplus\delta_2\oplus\dots\oplus\delta_{n-1}$.
\item $2\liminf_{n\to\infty}n\operatorname{artanh}\delta_n\ \le\ \gamma_*(L)\ \le\ \gamma^*(L)\
\le\ 2\limsup_{n\to\infty}n\operatorname{artanh}\delta_n$; if $\delta_n\to0$, the
$\operatorname{artanh}$ may be omitted.
\item If $n\,\delta_n\to\frac\gamma2\in[0,\infty]$, then $\gamma_*=\gamma^*=\gamma$. The
lacunarity hypothesis $\lambda_{j+1}/\lambda_j\ge\vartheta$ of
Theorem~\ref{infinite:thm:lacunary} is the condition $\delta_j\ge\frac{\vartheta-1}{\vartheta+1}$
for all $j$, and it implies $\gamma_*=\infty$.
\item The inequalities in (b) can be strict, also when $\gamma_*=\gamma^*$.
\end{enumerate}
\end{proposition}

\begin{proof}
(a) $\operatorname{artanh}\delta_k=\frac12\ln(\lambda_{k+1}/\lambda_k)$; summing and
using \eqref{infinite:eq:OLC-section} gives the first formula, and the second follows from
the addition theorem $\tanh(u+v)=\tanh u\oplus\tanh v$.

(b) Put $a_n=\ln\lambda_n-\ln\lambda_1=2\sum_{k<n}\operatorname{artanh}\delta_k$ and
$b_n=\ln n$. Since $b_n$ increases strictly to $\infty$, the Stolz--Ces\`aro theorem
gives
\[
\liminf\frac{a_{n+1}-a_n}{b_{n+1}-b_n}\le\liminf\frac{a_n}{b_n}\le\limsup\frac{a_n}{b_n}\le
\limsup\frac{a_{n+1}-a_n}{b_{n+1}-b_n},
\]
and $\frac{a_{n+1}-a_n}{b_{n+1}-b_n}=\frac{2\operatorname{artanh}\delta_n}{\ln(1+1/n)}$,
where $n\ln(1+1/n)\to1$. The middle terms are $\gamma_*$ and $\gamma^*$, since
$\ln\lambda_1/\ln n\to0$. If $\delta_n\to0$, then
$\operatorname{artanh}\delta_n/\delta_n\to1$.

(c) The first statement follows from (b) (if $n\delta_n\to\infty$, use
$\operatorname{artanh}\delta_n\ge\delta_n$). The second is
$\frac{\lambda_{j+1}-\lambda_j}{\lambda_{j+1}+\lambda_j}\ge\frac{\vartheta-1}{\vartheta+1}
\iff\lambda_{j+1}/\lambda_j\ge\vartheta$; then $n\delta_n\to\infty$.

(d) Let $r_n=1$ if $\lfloor\log_4n\rfloor$ is even and $r_n=3$ otherwise, and
$\lambda_1=1$, $\lambda_{n+1}=\lambda_n(1+\frac1n)^{r_n}$. Then
$2n\operatorname{artanh}\delta_n=r_n\,n\ln(1+\frac1n)$ has lower limit $1$ and upper
limit $3$. On the other hand $\sum_{4^j\le k<4^{j+1}}\ln(1+\frac1k)=\ln4$, so for
$4^J\le n<4^{J+1}$, writing $n=4^Jx$,
$\ln\lambda_n=\ln4\sum_{j<J}\rho_j+\rho_J\ln x$, where $\rho_j\in\{1,3\}$ is the value of
$r_k$ on the block $4^j\le k<4^{j+1}$ and $\sum_{j<J}\rho_j\in\{2J-1,2J\}$, and
$\ln\lambda_n/\ln n\to2$. Thus $\gamma_*=\gamma^*=2$, while the outer bounds in (b)
are $1$ and $3$.
\end{proof}

Part (a) is the equality case of the inequality
$\OMC(AB)\le\OMC(A)\oplus\OMC(B)$ of Theorem~\ref{core:subadd} for the commuting
factors of $L_n$ obtained by splitting its spectrum at consecutive eigenvalues:
the contrast of a section is the relativistic sum of its local contrasts, and its
logarithmic contrast, the sum of the corresponding rapidities. Part (b) says that
the logarithmic rates $\gamma_*,\gamma^*$ are bounded by the lower and upper mean
rapidity per unit of $\ln n$; geometric spectra $\lambda_j=\vartheta^j$ have constant
local contrast $\frac{\vartheta-1}{\vartheta+1}$, and power-type spectra have local
contrasts of order $\frac\gamma{2n}$ on logarithmic average,
$\frac1{\ln n}\sum_{k<n}2\operatorname{artanh}\delta_k\to\gamma$, though not termwise
(repeated eigenvalues give $\delta_k=0$, and (d) shows oscillation).

\subsection{Oscillating Weyl constants}\label{infinite:subsec:weylosc}

Theorem~\ref{infinite:thm:weyl} assumes $\lambda_n\sim c\,n^\gamma$. Suppose instead that, for
some $\gamma\in(0,\infty)$ (then necessarily $\gamma_*=\gamma^*=\gamma$),
\[
0<D_*=\liminf_{n\to\infty}\frac{\lambda_n}{n^\gamma}\le D^*=\limsup_{n\to\infty}\frac{\lambda_n}{n^\gamma}<\infty ,
\]
and let $\varpi(L)=\frac12\ln(D^*/D_*)=\OLC\bigl(\operatorname{diag}(D_*,D^*)\bigr)$ be
the logarithmic contrast of the two \emph{Weyl constants}; its Michelson form is
$(D^*-D_*)/(D^*+D_*)$. Clearly $\varpi(cL)=\varpi(L)$, and $\varpi(L)=0$ exactly when
$\lambda_n\sim c\,n^\gamma$.

\begin{proposition}\label{infinite:prop:weylosc}
Under these assumptions, with $\nu_\gamma$ as in Theorem~\ref{infinite:thm:weyl}:
\begin{enumerate}
\item $\limsup_n\bigl(\OLC(L_n)-\frac\gamma2\ln n\bigr)-\liminf_n\bigl(\OLC(L_n)-\frac\gamma2\ln n\bigr)=\varpi(L)$;
\item for $0<p<1/\gamma$: $\limsup_n\bigl|C_p(L_n)-C_p(\nu_\gamma)\bigr|\le\varpi(L)$;
for $p=1/\gamma$: $\limsup_n\bigl|C_p(L_n)-\frac\gamma2(\ln\ln n-\ln2)\bigr|\le\varpi(L)$;
\item for $p>1/\gamma$, with
$R_p(n)=C_p(L_n)-\bigl(\frac\gamma2-\frac1{2p}\bigr)\ln n-\frac1{2p}\ln\frac{\zeta_L(p)}{1+\gamma p}$:
\[
\tfrac12\ln D_*\ \le\ \liminf_nR_p(n)\ \le\ \limsup_nR_p(n)\ \le\ \tfrac12\ln D^* ,
\]
so the oscillation of $R_p(n)$ is at most $\varpi(L)$.
\end{enumerate}
For $\varpi(L)=0$, (b) and (c) reduce to Theorem~\ref{infinite:thm:weyl}(b)--(d).
\end{proposition}

\begin{proof}
(a) By \eqref{infinite:eq:OLC-section},
$\OLC(L_n)-\frac\gamma2\ln n=\frac12\ln(\lambda_n/n^\gamma)-\frac12\ln\lambda_1$.

For (b) and (c) fix $\varepsilon\in(0,D_*)$ and $j_0$ such that
$(D_*-\varepsilon)j^\gamma\le\lambda_j\le(D^*+\varepsilon)j^\gamma$ for $j\ge j_0$. For
$s>-1$, $\sum_{j\le n}j^{s}\sim n^{s+1}/(s+1)$, and the terms with $j<j_0$ are
negligible when the sum diverges. Hence, for every $p>0$,
\[
(D_*-\varepsilon)^p\le\liminf_n\frac{(1+\gamma p)}{n^{\gamma p}}\cdot\frac1n\sum_{j\le n}\lambda_j^p
\le\limsup_n\frac{(1+\gamma p)}{n^{\gamma p}}\cdot\frac1n\sum_{j\le n}\lambda_j^p\le(D^*+\varepsilon)^p,
\]
and, for $\gamma p<1$, the same holds for
$\frac{(1-\gamma p)}{n^{-\gamma p}}\cdot\frac1n\sum_{j\le n}\lambda_j^{-p}$ with the
bounds $(D^*+\varepsilon)^{-p}$ and $(D_*-\varepsilon)^{-p}$.

(b) Multiplying, and using
$2p\,C_p(L_n)=\ln\frac1n\sum_{j\le n}\lambda_j^p+\ln\frac1n\sum_{j\le n}\lambda_j^{-p}$ and
$2p\,C_p(\nu_\gamma)=-\ln(1-\gamma^2p^2)$,
\[
\limsup_n\bigl|2p\,C_p(L_n)-2p\,C_p(\nu_\gamma)\bigr|\le p\ln\frac{D^*+\varepsilon}{D_*-\varepsilon}.
\]
Divide by $2p$ and let $\varepsilon\to0$. For $p=1/\gamma$ the first sum is as before,
and $\frac1n\sum_{j\le n}\lambda_j^{-p}=\frac1n\sum_{j\le n}d_j^{-p}j^{-1}$ with
$d_j=\lambda_j/j^\gamma$ lies between $(D^*+\varepsilon)^{-p}\frac{\ln n}n$ and
$(D_*-\varepsilon)^{-p}\frac{\ln n}n$ asymptotically; with
$2p\,C_p=\ln\ln n-\ln2$ in the case $\varpi=0$ (Theorem~\ref{infinite:thm:weyl}(c)) the same
argument gives the bound.

(c) For $\gamma p>1$ we have $\zeta_L(p)<\infty$, since $\gamma_*=\gamma$
(Theorem~\ref{infinite:thm:abscissa}), and $\frac1n\sum_{j\le n}\lambda_j^{-p}=\frac{\zeta_L(p)}n(1+o(1))$.
Put $\tilde d_n=\bigl((1+\gamma p)\,n^{-\gamma p-1}\sum_{j\le n}\lambda_j^p\bigr)^{1/p}$, a
weighted mean of the $d_j=\lambda_j/j^\gamma$ (it is not $d_n$ itself, and lags behind it).
Hence $2p\,C_p(L_n)=(\gamma p-1)\ln n+\ln\frac{\zeta_L(p)}{1+\gamma p}+p\ln\tilde d_n+o(1)$,
and by the bounds at the beginning of the proof
$\liminf\tilde d_n\ge D_*-\varepsilon$ and $\limsup\tilde d_n\le D^*+\varepsilon$; divide by
$2p$. For $\varpi=0$ these are the statements of Theorem~\ref{infinite:thm:weyl}, with
$c=D_*=D^*$.
\end{proof}

\begin{example}[The Sierpinski gasket]\label{infinite:ex:gasket}
For the standard Laplacian on the Sierpinski gasket (with Dirichlet conditions, or
the Neumann Laplacian plus the identity; neither $D_*$ nor $D^*$ changes, since by
Dirichlet--Neumann bracketing the two counting functions differ by a bounded amount
\cite{KL93}, and the shift by $1$ does not change $\lambda_n/n^\gamma$ asymptotically), the counting
function satisfies $c_1\lambda^{d_s/2}\le N(\lambda)\le c_2\lambda^{d_s/2}$ for large
$\lambda$ \cite{KL93}, with spectral dimension $d_s=2\ln3/\ln5$, but
$N(\lambda)\lambda^{-d_s/2}$ has no limit as $\lambda\to\infty$ \cite{FS92}. With
$\gamma=2/d_s$, the two-sided bound gives
$(n/c_2)^\gamma\le\lambda_n\le((n-1)/c_1)^\gamma$ for large $n$, so $0<D_*\le D^*<\infty$; and
$D_*=D^*=c$ would give $N(\lambda)\sim(\lambda/c)^{1/\gamma}$. Hence
$\gamma_*=\gamma^*=\gamma=\ln5/\ln3$, $\Delta_{\rm ord}=0$, but $\varpi>0$:
Theorem~\ref{infinite:thm:weyl} does not apply, and Proposition~\ref{infinite:prop:weylosc} does. The
value of $\varpi$ is not computed here. An elementary example with $\varpi$ computed is
the sequence of Proposition~\ref{infinite:prop:local}(d): there
$\ln\lambda_n-2\ln n\in[-\ln4,0]$, with both ends approached, so $\gamma=2$,
$D_*=\frac14$, $D^*=1$ and $\varpi=\ln2$.
\end{example}

Thus the spectral information of a section sequence is graded: the lower rate
$\gamma_*$ (seen by the zeta function), the order contrast $\Delta_{\rm ord}$, and,
when $\Delta_{\rm ord}=0$, the contrast $\varpi$ of the Weyl constants; the last two are
Michelson (or logarithmic) contrasts of an upper and a lower limit.

\section{Essential contrast}\label{infinite:sec:essential}

For a bounded operator $A$ on an infinite-dimensional separable Hilbert space
$\cH$, $\OLC(A)$ sees the whole spectrum, while the Schatten contrasts with respect
to the trace of $B(\cH)$ are infinite for invertible $A$ (see the beginning of Section~\ref{infinite:sec:sections}) and trivially otherwise. The upper and lower
limits of the spectrum are seen by the Calkin algebra
$\mathcal Q(\cH)=B(\cH)/\mathcal K(\cH)$. Let $\pi\colon B(\cH)\to\mathcal Q(\cH)$ be the
quotient map, and recall that the essential spectrum $\sigma_{\rm ess}(A)$ is the
spectrum of $\pi(A)$. For $A$ with $\pi(A)$ invertible put
\[
\OLC_{\rm ess}(A)=\OLC(\pi(A))=\tfrac12\ln\bigl(\|\pi(A)\|\,\|\pi(A)^{-1}\|\bigr),
\]
\[
\OMC_{\rm ess}(A)=\tanh\OLC_{\rm ess}(A),
\]
the contrasts of \cite{ANN24,NN26} in the unital $C^*$-algebra $\mathcal Q(\cH)$.

\begin{proposition}\label{infinite:prop:ess}
\begin{enumerate}
\item If $A\ge0$ and $\pi(A)$ is invertible, then, with $a^*=\max\sigma_{\rm ess}(A)$
and $a_*=\min\sigma_{\rm ess}(A)>0$,
$\OMC_{\rm ess}(A)=\frac{a^*-a_*}{a^*+a_*}$. For a diagonal operator
$A=\operatorname{diag}(\alpha_n)$ with $\alpha_n>0$, $\sup\alpha_n<\infty$ and
$\liminf\alpha_n>0$, one has $a^*=\limsup\alpha_n$ and $a_*=\liminf\alpha_n$.
\item For invertible $A$, $\OLC_{\rm ess}(A)\le\OLC(A)$. For $A,B$ with $\pi(A),\pi(B)$
invertible, compact $K$, $\lambda\ne0$ and unitaries $U,V$:
$\OLC_{\rm ess}(A+K)=\OLC_{\rm ess}(A)$, $\OLC_{\rm ess}(\lambda UAV)=\OLC_{\rm ess}(A)$,
and $\OLC_{\rm ess}(AB)\le\OLC_{\rm ess}(A)+\OLC_{\rm ess}(B)$.
\item If $A\ge0$ and $\pi(A)$ is invertible, then
\[
\OLC_{\rm ess}(A)=\inf\bigl\{\OLC(A+K):\ K=K^*\ \text{compact},\ A+K\ \text{positive and invertible}\bigr\},
\]
and the infimum may be taken over self-adjoint operators $K$ of finite rank.
\end{enumerate}
\end{proposition}

\begin{proof}
(a) $\pi(A)$ is a positive invertible element of a $C^*$-algebra with spectrum
$\sigma_{\rm ess}(A)$, so $\|\pi(A)\|=a^*$ and $\|\pi(A)^{-1}\|=1/a_*$. For a diagonal
operator the essential spectrum is the set of cluster values of the sequence
$(\alpha_n)$, i.e.\ of the $t$ with $\#\{n:|\alpha_n-t|<\varepsilon\}=\infty$ for all
$\varepsilon>0$, whose largest and smallest elements are $\limsup\alpha_n$ and $\liminf\alpha_n$.

(b) $\|\pi(A)\|\le\|A\|$ and $\|\pi(A)^{-1}\|=\|\pi(A^{-1})\|\le\|A^{-1}\|$. The second
and third statements hold because $\pi(A+K)=\pi(A)$ and $\pi(U),\pi(V)$ are unitary,
and the last because the condition number is submultiplicative in any unital
$C^*$-algebra.

(c) By (b), $\OLC_{\rm ess}(A)=\OLC_{\rm ess}(A+K)\le\OLC(A+K)$ for every admissible
$K$. Conversely let $0<\varepsilon<a_*$ and
$P=\one_{\R\setminus(a_*-\varepsilon,\,a^*+\varepsilon)}(A)$. The part of $\sigma(A)$
outside $(a_*-\varepsilon,a^*+\varepsilon)$ is a compact set of isolated eigenvalues of
finite multiplicity, hence finite, so $P$ has finite rank. With
$c=\frac12(a_*+a^*)$ and $K=(c-A)P$, the operator $A+K=A(1-P)+cP$ is positive, its
spectrum lies in $[a_*-\varepsilon,a^*+\varepsilon]$, and
$\OLC(A+K)\le\frac12\ln\frac{a^*+\varepsilon}{a_*-\varepsilon}\to\OLC_{\rm ess}(A)$.
\end{proof}

In a factor of type $\mathrm{II}_\infty$, statements (b) and (c) hold for the
norm-closed ideal $\mathcal J$ generated by the projections of finite trace \cite{Breuer68}
in place of $\mathcal K(\cH)$, with ``finite rank'' replaced by ``with support projection of
finite trace'', and the proof of (c) modified as follows. Take a continuous $g$ with
$0\le g\le1$, $g=0$ on $[a_*,a^*]$ and $g=1$ outside $(a_*-\varepsilon,a^*+\varepsilon)$;
then $\pi(g(A))=g(\pi(A))=0$, so $g(A)\in\mathcal J$, and the spectral projection
$P=\mathbf 1_S(A)$ of the set $S$ where $g=1$ satisfies $P=Pg(A)$; hence $P\in\mathcal J$ is a
projection in a norm-closed ideal of a semifinite algebra, so it has finite trace. For unbounded operators
with compact resolvent the essential spectrum is empty and the essential
contrast is not defined; its role is taken by the order contrast and by the
contrast of the Weyl constants of Section~\ref{infinite:sec:sections}, which are again
contrasts of an upper and a lower limit.

\section{Concluding remarks}

The infinite contrast of an unbounded operator is graded by several numbers: the
threshold $p_c$ and the type of the transition there, for operators affiliated with
a finite algebra; for finite sections, the lower rate $\gamma_*$ and two contrasts of
upper and lower limits, the scale-free order contrast and the contrast $\varpi$ of the
Weyl constants; and the essential contrast, for bounded operators read modulo
compact ones. We mention three open questions.

\begin{enumerate}[label=\textup{(\arabic*)}]
\item \emph{Weyl remainders.} Theorem~\ref{infinite:thm:weyl} and
Proposition~\ref{infinite:prop:weylosc} use only the leading term of the counting function.
How do the remainder term in Weyl's law, and more generally regularly varying
counting functions, enter the next terms of $C_p(L_n)$?
\item \emph{Sharpness.} Is the bound $\varpi$ in Proposition~\ref{infinite:prop:weylosc}(b)
attained? For log-periodic Weyl constants, as on p.c.f.\ self-similar fractals,
does $C_p(L_n)$ itself oscillate log-periodically below the abscissa, and with
which amplitude?
\item \emph{Non-spectral sections.} What is the growth of the contrasts for the
sections used in numerical analysis (finite elements), for which
Proposition~\ref{infinite:prop:minmax} gives only a lower bound?
\end{enumerate}

\chapter{Contrasts of random products}\label{ch:random}
\section{Introduction}

The operator logarithmic contrast $\OLC(A)=\frac12\ln\kappa(A)$,
$\kappa(A)=\|A\|\,\|A^{-1}\|$, and the operator Michelson contrast
$\OMC(A)=\tanh\OLC(A)$ were introduced and studied in
\cite{ANN24,NN26}: $\OLC$ is subadditive under products, homogeneous under
powers of positive elements, and vanishes exactly on multiples of unitaries. The Schatten contrasts
$C_p$ replace the operator norm by normalised Schatten norms; the elementary facts about
them that we use are collected in Section~\ref{random:sec:prelim} or proved where they are used. In
this chapter we ask how these quantities behave along products of independent
random factors, and along the opposite process, convergence to the tracial
state.

\medskip\noindent\textbf{Growth: three kinds of independence.} Let
$A_N=X_N\cdots X_1$ with identically distributed invertible factors.

\begin{enumerate}[label=\textup{(\Roman*)}]
\item \emph{Tensor independence.} If the factors sit in different tensor
factors, every contrast is exactly additive:
$C_p(A_N)=N\,C_p(X)$ for all $p\in(0,\infty]$ (Proposition~\ref{random:prop:tensor}).
The slopes depend on $p$.
\item \emph{Matrix independence.} For i.i.d.\ factors in $\GL_d(\C)$, $d$ fixed,
all contrasts have the same slope $\frac12(\gamma_1-\gamma_d)$, where
$\gamma_1\ge\gamma_d$ are the extreme Lyapunov exponents, and $\OLC(A_N)$
fluctuates on the scale $\sqrt N$ (Theorem~\ref{random:thm:matrix}). This is the
theorem of Furstenberg and Kesten and the central limit theorem of
Goldsheid--Guivarc'h and Benoist--Quint, read for contrasts; the existence of the limit
follows directly from the subadditivity of $\OLC$ by Kingman's theorem, and its value from
Furstenberg--Kesten.
\item \emph{Free independence.} For $*$-free factors in a finite von Neumann
algebra, $C_2$ is exactly additive, all slopes collapse to $C_2(X)$, and
\[
C_p(A_N)=N\,C_2(X)+\Bigl(\frac12-\frac1p\Bigr)\ln N+O(1)\qquad(2\le p<\infty),
\]
with an explicit constant in place of $O(1)$ for even integers $p$, and
\[
\OLC(A_N)=N\,C_2(X)+\frac12\ln N+\frac12+\frac14\ln(vv')+o(1),
\]
where $v,v'$ are normalised variances of $X^*X$ and $(X^*X)^{-1}$
(Theorem~\ref{random:thm:free}).
\end{enumerate}

The second-order term of $\OLC(A_N)$ is therefore $0$, of order $\sqrt N$ and
random, or $\frac12\ln N$ and deterministic, in the three cases. For $p<2$ we
prove lower bounds with the conjectured logarithmic term, and for all
$p\notin2\N$ we conjecture the full expansion with an explicit constant,
supported by high-precision numerical evaluation of the exact free formulas
(Conjecture~\ref{random:conj:free}). For Haar-rotated factors, the matrix slope is
bounded by the free slope $C_2(X)$ and converges to it as $d\to\infty$
(Proposition~\ref{random:prop:isotropic}).

\medskip\noindent\textbf{Decay: convergence to the tracial state.} A positive
element has contrast $0$ exactly when it is a multiple of $1$, i.e.\ when the
corresponding state is the trace. Since $2\OLC(\rho)$ is the Hilbert projective
distance from $\rho$ to $1$, the Birkhoff--Hopf theorem gives geometric decay of
$\OLC$ under unital positive maps of finite projective diameter, with the
contraction ratio $\tanh(\Delta/4)$ (Proposition~\ref{random:prop:BH}). For convolution powers on the circle, the
Michelson contrast of the density decays like $2|a_{k_0}|\,|\hat\mu(k_0)|^N$
(Proposition~\ref{random:prop:circle}); this is the classical fact that phase noise
multiplies the visibility of interference fringes by the modulus of the characteristic
function of the noise, and it is in this sense that the Michelson contrast was
introduced in optics. For random unital channels the decay rate is a Lyapunov
exponent (Proposition~\ref{random:prop:channels}).

\medskip\noindent\textbf{What is new.} Proposition~\ref{random:prop:isotropic} (comparison of slopes), Proposition~\ref{random:prop:C2},
Theorem~\ref{random:thm:free}(b)--(d) and Conjecture~\ref{random:conj:free} appear to be new;
Theorem~\ref{random:thm:free}(a) is a direct corollary of results of Kargin and of
Deng--Zhong on free multiplicative convolution powers \cite{Kargin08,DZ19}, and
the integer moment asymptotics used in the proof are a short Lagrange inversion
computation, related to limit theorems in \cite{SY13,AV12}. Theorem~\ref{random:thm:matrix} and the results of
Section~\ref{random:sec:decay} are reformulations of classical theorems, included
because, in the language of contrasts, they complete the comparison and give
the growth and decay of contrasts a common form. Table~\ref{random:tab:status}
summarises the status of the results.

\begin{table}[!ht]
\centering\small
\begin{tabular}{@{}p{0.20\textwidth}p{0.36\textwidth}p{0.36\textwidth}@{}}
\toprule
Result & Background & Role in this chapter\\
\midrule
Prop.~\ref{random:prop:tensor} & additivity of the logarithmic spectrum under tensor products; Chapter~\ref{ch:core} & ---\\[2pt]
Lemma~\ref{random:lem:holder} & H\"older's inequality \cite{FK86,PX03} & ---\\[2pt]
Thm.~\ref{random:thm:matrix} & \cite{FK60,Kingman68,GG96,BQ16a,BQ16} & reformulation\\[2pt]
Prop.~\ref{random:prop:isotropic} & isotropic formula \cite{CN84} & comparison with the free slope; convergence\\[2pt]
Prop.~\ref{random:prop:C2} & multiplicativity of second moments under freeness & exact additivity of $C_2$\\[2pt]
Thm.~\ref{random:thm:free}(a) & \cite{Kargin07,Kargin08,DZ19} & contrast form, constant\\[2pt]
Thm.~\ref{random:thm:free}(b)--(d) & moments of $\mu^{\boxx N}$ \cite{AV12,SY13} & Schatten asymptotics, collapse of slopes\\[2pt]
Prop.~\ref{random:prop:BH} & Birkhoff--Hopf \cite{Birkhoff57,Bushell73b,RKW11} & contrast form\\[2pt]
Prop.~\ref{random:prop:circle} & Fourier analysis; Kawada--It\^o \cite{KI40}; visibility \cite{BornWolf} & contrast form\\[2pt]
Prop.~\ref{random:prop:channels} & Furstenberg--Kesten, Furstenberg \cite{FK60,BL85} & contrast form\\
\bottomrule
\end{tabular}
\caption{Status of the results.}
\label{random:tab:status}
\end{table}

\section{Preliminaries}\label{random:sec:prelim}

Throughout, $\cM$ is a von Neumann algebra with a faithful normal tracial state
$\tau$, and for $X\in\cM$ and $0<p<\infty$ we write
$\|X\|_p=\tau(|X|^p)^{1/p}$ (a quasi-norm for $p<1$) and $\|X\|_\infty=\|X\|$.
For invertible $A\in\cM$ put
\[
C_p(A)=\tfrac12\bigl(\ln\|A\|_p+\ln\|A^{-1}\|_p\bigr)
=\frac1{2p}\ln\bigl(\tau(|A|^p)\,\tau(|A^{-1}|^p)\bigr),\qquad
C_\infty(A)=\OLC(A).
\]
Let $\nu_A$ be the distribution of $\ln|A|$ with respect to $\tau$. We use the
following elementary facts, with $C_p(\nu)=\frac1{2p}\ln\bigl(\int e^{px}\dd\nu\int e^{-px}\dd\nu\bigr)$:
\begin{enumerate}[label=\textup{(\roman*)}]
\item $C_p(A)=C_p(\nu_A)$: in the polar decomposition $A=U|A|$ the factor $U$ is
unitary and $|A^{-1}|=U|A|^{-1}U^*$, so $\tau(|A^{-1}|^p)=\tau(|A|^{-p})=\int e^{-px}\dd\nu_A$;
in particular $C_p(A^{-1})=C_p(A)$;
\item $C_p(A)\ge0$, with equality iff $|A|$ is a multiple of $1$, by the
Cauchy--Schwarz inequality $1=(\int e^{px/2}e^{-px/2}\dd\nu_A)^2\le\int e^{px}\dd\nu_A\int e^{-px}\dd\nu_A$;
\item $p\mapsto C_p(A)$ is nondecreasing with limit $\OLC(A)$, since on a probability
space the $L^p$-norms of $e^{x}$ and $e^{-x}$ increase to the essential suprema;
\item if $\nu_{A'}=\nu_A*\nu_B$ then $C_p(A')=C_p(A)+C_p(B)$, since Laplace transforms
multiply under convolution.
\end{enumerate}

\begin{lemma}\label{random:lem:holder}
Let $X_1,\dots,X_N\in\cM$ be invertible, $p_1,\dots,p_N\in(0,\infty]$ and
$\frac1r=\sum_i\frac1{p_i}$. Then
\[
C_r(X_N\cdots X_1)\le\sum_iC_{p_i}(X_i).
\]
In particular $C_{p/N}(X_N\cdots X_1)\le\sum_iC_p(X_i)$ and
$\OLC(X_N\cdots X_1)\le\sum_i\OLC(X_i)$.
\end{lemma}

\begin{proof}
By the H\"older inequality $\|XY\|_r\le\|X\|_p\|Y\|_q$,
$\frac1r=\frac1p+\frac1q$, valid for all $0<p,q\le\infty$ \cite{FK86,PX03},
applied repeatedly to the product and to its inverse
$X_1^{-1}\cdots X_N^{-1}$.
\end{proof}

\paragraph{Free multiplicative convolution.} For compactly supported probability
measures $\mu,\mu'$ on $(0,\infty)$, $\mu\boxx\mu'$ is the distribution of
$a^{1/2}ba^{1/2}$ for free positive $a,b$ with distributions $\mu,\mu'$
\cite{VDN92,NS06}. With the moment series $\psi_\mu(z)=\sum_{k\ge1}m_k(\mu)z^k$,
its compositional inverse $\chi_\mu$ and the $S$-transform
$S_\mu(w)=\frac{1+w}w\chi_\mu(w)$, one has $S_{\mu\boxx\mu'}=S_\mu S_{\mu'}$. If
$m_1(\mu)=1$ and $\mu$ has variance $V$, then $S_\mu(w)=1-Vw+O(w^2)$.

\section{Tensor independence and a general bound}\label{random:sec:tensor}

\begin{proposition}\label{random:prop:tensor}
For $i=1,\dots,N$ let $X_i$ be an invertible element of a von Neumann
algebra $\cM_i$ with a faithful normal tracial state $\tau_i$, and let $A_N=X_1\otimes\dots\otimes X_N$ in the tensor
product $\bigotimes_i(\cM_i,\tau_i)$. Then $C_p(A_N)=\sum_iC_p(X_i)$ for every
$p\in(0,\infty]$. In particular, for identically distributed $X_i$,
$C_p(A_N)=N\,C_p(X)$, and the slope $C_p(X)$ depends on $p$ unless $|X|$ is a
multiple of $1$.
\end{proposition}

\begin{proof}
$|A_N|=\bigotimes|X_i|$, so $\nu_{A_N}=\nu_{X_1}*\dots*\nu_{X_N}$; for $p<\infty$
use the additivity of $C_p$ under convolution, and for $p=\infty$ that the
extreme points of the supports add (equivalently, the case $p=\infty$ follows from the
finite-$p$ case and fact (iii), $C_p\to\OLC$; for positive definite matrices it is also
Chapter~\ref{ch:core}). If $|X|$ is not a multiple of $1$, then
$p\mapsto C_p(X)$ is strictly increasing, since the $L^p$-norms of a
nonconstant function on a probability space are strictly increasing in $p$.
\end{proof}

The tensor product is the commuting model: the factors
$X_i\otimes1\otimes\cdots$ commute, and the product of commuting independent
positive elements behaves in the same way. Lemma~\ref{random:lem:holder} gives, for
arbitrary products in $\cM$, the bound $C_{p/N}(A_N)\le\sum_iC_p(X_i)$; in the
tensor case $C_{p/N}(A_N)=N\,C_{p/N}(X)\le N\,C_p(X)$, so the bound loses only
the change of exponent.

\section{Matrix independence}\label{random:sec:matrix}

Let $X_1,X_2,\dots$ be i.i.d.\ random elements of $\GL_d(\C)$ with
$\E\ln^+\|X_1\|<\infty$ and $\E\ln^+\|X_1^{-1}\|<\infty$, $A_N=X_N\cdots X_1$, and
let $\gamma_1\ge\gamma_d$ be the largest and the smallest Lyapunov exponent:
$\frac1N\ln\|A_N\|\to\gamma_1$ and $\frac1N\ln\|A_N^{-1}\|\to-\gamma_d$ almost
surely \cite{FK60}. Contrasts in $M_d(\C)$ are taken with the normalised trace;
for a symmetric gauge function $\Phi$ on $\R^d$ normalised by $\Phi(e_1)=1$,
with $s(A)$ the vector of singular values and
$\kappa_\Phi(A)=\Phi(s(A))\,\Phi(s(A^{-1}))$, we use the normalised contrast
$C_\Phi(A)=\frac12\ln\bigl(\kappa_\Phi(A)/\kappa_\Phi(I)\bigr)$, where
$\kappa_\Phi(I)=\Phi(\mathbf1)^2$.

\begin{theorem}\label{random:thm:matrix}
\begin{enumerate}
\item Almost surely, $\frac1N\OLC(A_N)\to\frac12(\gamma_1-\gamma_d)$, and the same
limit holds for $C_p(A_N)$, $0<p<\infty$, and for $C_\Phi(A_N)$ for every
symmetric gauge function $\Phi$.
\item More precisely, $0\le\OLC(A_N)-C_p(A_N)\le\frac{\ln d}p$ and
$0\le\OLC(A_N)-C_\Phi(A_N)\le\ln\Phi(\mathbf 1)$ for every $N$.
\item If, in addition, $\ln\|X_1\|$ and $\ln\|X_1^{-1}\|$ have finite second
moments and the support of the law of $X_1$ generates a semigroup that is
Zariski dense in $\GL_d(\C)$ regarded as a real algebraic group (or in
$\GL_d(\R)$ for real factors), then
$\bigl(\OLC(A_N)-\frac N2(\gamma_1-\gamma_d)\bigr)/\sqrt N$ converges in
distribution to a centred Gaussian law, nondegenerate for $d\ge2$.
\end{enumerate}
\end{theorem}

\begin{proof}
(a) $\OLC$ is subadditive (Lemma~\ref{random:lem:holder} with all $p_i=\infty$; see also \cite{NN26}), so $N\mapsto\OLC(X_N\cdots X_1)$ is a
subadditive process over the shift, and Kingman's subadditive ergodic theorem
\cite{Kingman68} gives an almost sure limit of $\OLC(A_N)/N$; it equals
$\frac12(\gamma_1-\gamma_d)$ by the theorem of Furstenberg and Kesten applied to
$\ln\|A_N\|$ and $\ln\|A_N^{-1}\|$. The other contrasts follow from (b).
(b) If $s_1\ge\dots\ge s_d$ are the singular values of $A$, then
$\tau(|A|^p)\ge s_1^p/d$ and $\tau(|A|^{-p})\ge s_d^{-p}/d$, so
$C_p(A)\ge\OLC(A)-\frac{\ln d}p$; and $s_1\le\Phi(s)\le\Phi(\mathbf1)s_1$ for the
normalised $\Phi$, and likewise for $A^{-1}$, so
$\kappa(A)\le\kappa_\Phi(A)\le\Phi(\mathbf1)^2\kappa(A)$, which gives both bounds
for $C_\Phi$. Finally $C_p\le\OLC$ because $p\mapsto C_p$ increases to $\OLC$.
(c) $2\OLC(A_N)=\ln s_1(A_N)-\ln s_d(A_N)$ is a linear functional of the Cartan
projection of $A_N$. The central limit theorem for the Cartan projection is due
to Goldsheid and Guivarc'h \cite{GG96} under exponential moments and to Benoist
and Quint \cite[Theorems~4.11, 4.16]{BQ16a} (see also \cite{BQ16}) under a second
moment; for Zariski dense laws the Gaussian is nondegenerate on the Cartan
subspace of the semisimple part, on which $\ln s_1-\ln s_d$ does not vanish
for $d\ge2$. For the norm cocycle alone see Le~Page \cite{LePage82}.
\end{proof}

\section{Free independence}\label{random:sec:free}

In this section $X_1,X_2,\dots\in\cM$ are $*$-free and invertible, and $A_N=X_N\cdots X_1$.

\begin{proposition}\label{random:prop:C2}
For $*$-free $X_1,\dots,X_N$, not necessarily identically distributed,
\[
\tau(|A_N|^2)=\prod_i\tau(|X_i|^2),\qquad
\tau(|A_N^{-1}|^2)=\prod_i\tau(|X_i^{-1}|^2),\qquad
C_2(A_N)=\sum_iC_2(X_i).
\]
\end{proposition}

\begin{proof}
$\tau(|A_N|^2)=\tau(A_{N-1}^*X_N^*X_NA_{N-1})=\tau\bigl(X_N^*X_N\cdot A_{N-1}A_{N-1}^*\bigr)$,
and $X_N^*X_N$ is free from $A_{N-1}A_{N-1}^*\in W^*(X_1,\dots,X_{N-1})$, so the
trace factorises; induct on $N$. For the inverse,
$|A_N^{-1}|^2=(X_NA_{N-1}A_{N-1}^*X_N^*)^{-1}$, so
$\tau(|A_N^{-1}|^2)=\tau\bigl((X_N^*X_N)^{-1}\cdot|A_{N-1}^{-1}|^2\bigr)$, with the
same freeness, and
$\tau((X_N^*X_N)^{-1})=\tau((X_NX_N^*)^{-1})=\tau(|X_N^{-1}|^2)$.
\end{proof}

By Theorem~\ref{random:thm:free}(d) below and the strict monotonicity of
$p\mapsto C_p(X)$, $p=2$ is the only exponent for which
such an identity holds for identically distributed factors with non-scalar
$X^*X$.

For the rest of the section the $X_i$ are identically distributed. Let $\mu$ be
the distribution of $X^*X$ and $\mu^{-1}$ that of $(X^*X)^{-1}$, and put
\[
m=\tau(X^*X),\quad m'=\tau\bigl((X^*X)^{-1}\bigr),\quad
v=\frac{\Var(\mu)}{m^2},\quad v'=\frac{\Var(\mu^{-1})}{m'^2},\quad
C_2(X)=\tfrac14\ln(mm').
\]
Since $XX^*$ and $X^*X$ have the same distribution, $m'=\tau(|X^{-1}|^2)$.

\begin{theorem}\label{random:thm:free}
Assume that $X^*X$ is not a multiple of $1$. Then:
\begin{enumerate}
\item $\displaystyle\OLC(A_N)=N\,C_2(X)+\frac12\ln N+\frac12+\frac14\ln(vv')+o(1)$.
\item For every $p\in[2,\infty)$,
$\displaystyle C_p(A_N)=N\,C_2(X)+\Bigl(\frac12-\frac1p\Bigr)\ln N+O(1)$; for
$p=2k$, $k\in\N$, the $O(1)$ term is $\kappa_p+O(N^{-1})$ with
\[
\kappa_p=\frac1{2p}\Bigl[\Bigl(\frac p2-1\Bigr)\ln\frac{vv'p^2}4-2\ln\Gamma\Bigl(1+\frac p2\Bigr)\Bigr].
\]
\item For every $p\in(0,2)$,
$\displaystyle N\,C_2(X)+\Bigl(\frac12-\frac1p\Bigr)\ln N-O(1)\le C_p(A_N)\le N\,C_2(X)$.
\item Consequently $\lim_{N\to\infty}C_p(A_N)/N=C_2(X)$ for every $p\in(0,\infty]$.
\end{enumerate}
\end{theorem}

\begin{proof}
\emph{Step 1: distributions.} The distribution of $|A_N|^2$ is $\mu^{\boxx N}$.
Indeed $A_N^*A_N=A_{N-1}^*YA_{N-1}$ with $Y=X_N^*X_N$, whose moments are
$\tau\bigl((Y^{1/2}A_{N-1}A_{N-1}^*Y^{1/2})^k\bigr)$; the element
$A_{N-1}A_{N-1}^*$ has the same distribution as $A_{N-1}^*A_{N-1}$ and is free
from $Y$, and both are bounded and positive, so the distribution is
$\mu\boxx\mu^{\boxx(N-1)}$. Similarly $|A_N^{-1}|^2=(A_NA_N^*)^{-1}$ has
distribution $(\mu^{\boxx N})^{-1}=(\mu^{-1})^{\boxx N}$, since
$(a^{1/2}ba^{1/2})^{-1}=a^{-1/2}b^{-1}a^{-1/2}$.

\emph{Step 2: norms.} Let $\mu_1$ be the image of $\mu$ under $x\mapsto x/m$, a
compactly supported measure on $(0,\infty)$ with mean $1$ and variance $v>0$.
By Kargin \cite{Kargin08} (for integer powers) and Deng--Zhong
\cite[Theorem~1.2]{DZ19}, the largest point of the support of
$\mu_1^{\boxx N}$ is $e\,v\,N(1+o(1))$. As $\mu^{\boxx N}$ is the image of
$\mu_1^{\boxx N}$ under $x\mapsto m^Nx$,
\begin{equation}\label{random:eq:norms}
\|A_N\|^2=m^N\,e\,v\,N\,(1+o(1)),\qquad\|A_N^{-1}\|^2=m'^N\,e\,v'\,N\,(1+o(1)),
\end{equation}
the second by the same argument for $\mu^{-1}$. Hence
$\OLC(A_N)=\frac14\ln(\|A_N\|^2\|A_N^{-1}\|^2)$ is as in (a).

\emph{Step 3: integer moments.} Normalise $m=1$. Then $S_{\mu^{\boxx N}}=S_\mu^N$,
so $\chi_{\mu^{\boxx N}}(w)=\frac w{1+w}S_\mu(w)^N$, i.e.\ $w=z\,\varphi(w)$ with
$z=\chi_{\mu^{\boxx N}}(w)$ and $\varphi(w)=(1+w)S_\mu(w)^{-N}$. By the Lagrange
inversion formula, for $k\ge1$,
\[
m_k(\mu^{\boxx N})=[z^k]\psi_{\mu^{\boxx N}}(z)=\frac1k\,[w^{k-1}]\,(1+w)^k\,S_\mu(w)^{-Nk}.
\]
Write $S_\mu(w)^{-1}=\exp(vw+w^2g(w))$ with $g$ holomorphic near $0$. The
coefficient of $w^{k-1}$ in $(1+w)^k\exp(Nk(vw+w^2g(w)))$ is a polynomial in $N$
of degree $k-1$: a monomial $N^jw^i$ arises only with $i\ge j$, so the term
$N^{k-1}$ comes only from $(Nkvw)^{k-1}/(k-1)!$. Hence
\begin{equation}\label{random:eq:moments}
m_k(\mu^{\boxx N})=\frac{k^{k-1}}{k!}(vN)^{k-1}\bigl(1+O(N^{-1})\bigr),
\end{equation}
and for general $m$ the right-hand side is multiplied by $m^{kN}$. (For $k=2$
this is the exact identity $m_2=1+vN$; for the Marchenko--Pastur law, where
$S_\mu(w)=(1+w)^{-1}$ and $v=1$, the exact formula above gives the Fuss--Catalan
numbers.) Thus
$\ln\tau(|A_N|^{2k})=kN\ln m+(k-1)\ln(vkN)-\ln\Gamma(k+1)+O(N^{-1})$.

\emph{Step 4: interpolation.} Write $\ell_s=\ln\tau(|A_N|^s)$ for $s>0$ and
$\ell_\infty=\ln\|A_N\|$, so that $\ell_2=N\ln m$ and, by \eqref{random:eq:norms},
$\ell_\infty=\frac N2\ln m+\frac12\ln N+O(1)$. The inequality
$\tau(|A|^s)\le\|A\|^{s-t}\tau(|A|^t)$ for $s\ge t$ gives:
\begin{itemize}
\item for $p\ge2$: $\ell_p\le\ell_2+(p-2)\ell_\infty$, and, for an integer
$k\ge p/2$, $\ell_p\ge\ell_{2k}-(2k-p)\ell_\infty$;
\item for $p<2$: $\ell_p\ge\ell_2-(2-p)\ell_\infty$, while
$\ell_p\le\frac p2\ell_2$ by Jensen's inequality.
\end{itemize}
Inserting the asymptotics of $\ell_2$, $\ell_{2k}$ and $\ell_\infty$, the upper
and lower bounds for $p\ge2$ both equal
$\frac1p\ell_p=\frac N2\ln m+(\frac12-\frac1p)\ln N+O(1)$; for $p<2$ the lower
bound is the same and the upper bound is $\frac N2\ln m$. The same holds for
$A_N^{-1}$ with $m',\mu^{-1}$. Since
$C_p(A_N)=\frac12(\frac1p\ell_p+\frac1p\ell_p')$, with $\ell'$ the analogous
quantities for $A_N^{-1}$, this gives (b) and (c), and (d) follows from (a)--(c).
For $p=2k$ no interpolation is needed: inserting the moment asymptotics for
$A_N$ and $A_N^{-1}$ directly into $C_{2k}(A_N)=\frac1{4k}(\ell_{2k}+\ell'_{2k})$
gives the constant $\kappa_{2k}$.
\end{proof}

\begin{remark}\label{random:rem:free-what}
Kargin \cite{Kargin07} proved $\frac1N\ln\|A_N\|\to\frac12\ln\tau(X^*X)$, and
\cite{Kargin08,DZ19} the linear growth of the support of $\mu^{\boxx N}$, so (a)
is a direct consequence; the integer moment asymptotics \eqref{random:eq:moments} are a
short Lagrange inversion computation, and related limit theorems for free
multiplicative powers are in \cite{SY13} (for Marchenko--Pastur laws see
\cite{AV12}). Parts (b)--(d), i.e.\ the Schatten norm asymptotics
$\|A_N\|_p\asymp m^{N/2}N^{\frac12-\frac1p}$ for $p\ge2$, the two-sided bounds for $p<2$,
and the collapse of all contrast slopes onto $C_2$, do not seem to have been stated. The mechanism is visible
in \eqref{random:eq:moments}: $|A_N|^2/m^N$ has mean $1$, but its mass near the top of
the spectrum, at scale $N$, is of order $N^{-1}$.
\end{remark}

\begin{conjecture}\label{random:conj:free}
Let $\mu$ be compactly supported in $(0,\infty)$ and not a point mass. For every $s>0$,
with $Y_N=|A_N|^2/m^N$ and $Y_N'=|A_N^{-1}|^2/m'^N$,
\[
\ln\tau(Y_N^s)=(s-1)\ln(vsN)-\ln\Gamma(s+1)+o(1),
\]
and the same for $Y_N'$ with $v'$ in place of $v$.
Consequently, for every $p\in(0,\infty)$,
\[
C_p(A_N)=N\,C_2(X)+\Bigl(\frac12-\frac1p\Bigr)\ln N+\kappa_p+o(1),
\]
with $\kappa_p$ as in Theorem~\ref{random:thm:free}(b).
\end{conjecture}

For integer $s=k$ the first formula is \eqref{random:eq:moments}, and the expansion of
$C_p$ is proved for $p\in2\N$ in Theorem~\ref{random:thm:free}(b); the conjecture
concerns $p\notin2\N$, in particular $p<2$. The constants are
consistent with Theorem~\ref{random:thm:free}: $\kappa_2=0$, and
$\kappa_p\to\frac12+\frac14\ln(vv')$ as $p\to\infty$, which is the constant in
(a). For the Marchenko--Pastur law, which has $0$ in its support and so lies outside the
invertible setting, the moments of $\mu^{\boxx N}$ at real $s$ are given by the Raney (Fuss--Catalan)
Mellin formula $\Gamma((N+1)s+1)/(\Gamma(s+1)\Gamma(Ns+2))$ \cite{MPZ13}, and the analogue of
the first formula follows from Stirling's formula; we checked the Mellin formula
numerically for $N=4$. A
plausible route to a proof is the vague convergence of
$vN\,(\mu_1^{\boxx N}\circ D_{vN})$, $D_c(x)=cx$, to a $\sigma$-finite measure
with $s$-th moment $s^{s-1}/\Gamma(s+1)$, together with uniform integrability
near $0$ for $s<1$. Numerical evidence is given in Section~\ref{random:sec:numerics}.

\section{The order of the remainder}\label{random:sec:compare}

\begin{center}
\begin{tabular}{@{}llll@{}}
\toprule
independence & slope $\lim_NC_p(A_N)/N$ & remainder of $\OLC(A_N)$ & reference\\
\midrule
tensor & $C_p(X)$, depends on $p$ & $0$ & Prop.~\ref{random:prop:tensor}\\
free & $C_2(X)$ for all $p$ & $\frac12\ln N+\mathrm{const}$, deterministic & Thm.~\ref{random:thm:free}\\
$d\times d$ matrices & $\frac12(\gamma_1-\gamma_d)$ for all $p$ & order $\sqrt N$, Gaussian & Thm.~\ref{random:thm:matrix}\\
\bottomrule
\end{tabular}
\end{center}

\noindent For $C_2$ the free and the tensor case coincide: both are exactly
additive. For fixed $N$, products of Haar-rotated $d\times d$ factors converge
to the free model as $d\to\infty$; for Haar unitaries the convergence is strong by
\cite{CM14} when the deterministic factors $D_d$ converge strongly (e.g.\ for
discretisations of a fixed spectrum), so that norms converge as well (we do not use the
orthogonal case of this); but the limits $N\to\infty$
and $d\to\infty$ need not commute, and for the slopes we argue directly.
For free factors the largest of Kargin's free Lyapunov exponents
\cite{Kargin08b} is $\frac12\ln\tau(X^*X)=\frac12\ln m$, the slope of
$\ln\|A_N\|$, and by symmetry the slope of $\ln\|A_N^{-1}\|$ is $\frac12\ln m'$;
together they give the free slope $C_2(X)$.

\begin{proposition}\label{random:prop:isotropic}
Let $D=\operatorname{diag}(\delta_1,\dots,\delta_d)$ be positive definite, let $U$
be Haar distributed on $O(d)$ (or $U(d)$), and let $X_i$ be i.i.d.\ copies of
$X=UD$. Then
\[
\tfrac12(\gamma_1-\gamma_d)=\tfrac14\,\E\Bigl[\ln\sum_i\delta_i^2u_i^2+\ln\sum_i\delta_i^{-2}u_i^2\Bigr]
\le C_2(X)=\tfrac14\ln\bigl(\tau(D^2)\,\tau(D^{-2})\bigr),
\]
where $u$ is uniformly distributed on the unit sphere of $\R^d$ (of $\C^d$, with
$u_i^2$ replaced by $|u_i|^2$). If $D=D_d$ has spectrum in a fixed interval
$[a,b]\subset(0,\infty)$ and the empirical distributions of its eigenvalues
converge weakly to a measure $\sigma$, then both sides converge, as
$d\to\infty$, to $\frac14\ln\bigl(\int x^2\dd\sigma\int x^{-2}\dd\sigma\bigr)$, the
free slope of Theorem~\ref{random:thm:free}.
\end{proposition}

\begin{proof}
The law of $X$ is invariant under $X\mapsto WX$ for orthogonal (unitary) $W$.
Fix a unit vector $e$ and put $v_0=e$, $v_i=X_iv_{i-1}/\|X_iv_{i-1}\|$. Then
$\|X_iv_{i-1}\|=\|Dv_{i-1}\|$, and since $U_i$ is Haar and independent of
$X_1,\dots,X_{i-1}$, $v_i$ is uniformly distributed and independent of
$v_0,\dots,v_{i-1}$. Hence $\ln\|A_Ne\|=\sum_{i\le N}\ln\|Dv_{i-1}\|$, and by the
law of large numbers $\frac1N\ln\|A_Ne\|\to\E\ln\|Du\|$ almost surely, for every
fixed $e$. Since $\|A_Ne\|\le\|A_N\|\le\sqrt d\max_j\|A_Ne_j\|$ for an
orthonormal basis $(e_j)$, $\gamma_1=\E\ln\|Du\|=\frac12\E\ln\sum_i\delta_i^2u_i^2$
(this is the observation of Cohen and Newman \cite{CN84}). Next,
$(A_N^{-1})^*=Y_N\cdots Y_1$ with $Y_i=X_i^{-*}=U_iD^{-1}$, so the same argument
gives $-\gamma_d=\lim_N\frac1N\ln\|A_N^{-1}\|=\E\ln\|D^{-1}u\|$. Since
$\E u_i^2=\frac1d$, Jensen's inequality gives
$\E\ln\sum_i\delta_i^{\pm2}u_i^2\le\ln\tau(D^{\pm2})$, and $\tau(X^*X)=\tau(D^2)$,
$\tau((X^*X)^{-1})=\tau(D^{-2})$. For the limit write $u=g/\|g\|$ with $g$ a
standard Gaussian vector: $\sum_i\delta_i^{\pm2}u_i^2$ is the quotient of
$\frac1d\sum_i\delta_i^{\pm2}g_i^2$, which differs from $\tau(D_d^{\pm2})$ by a
quantity of variance at most $2\max(a^{-4},b^4)/d$, and $\frac1d\sum_ig_i^2\to1$;
its logarithm is bounded by $2\max(|\ln a|,|\ln b|)$, so by bounded convergence the
expectations converge to $\ln\int x^{\pm2}\dd\sigma$, which is also the limit of
$\ln\tau(D_d^{\pm2})$.
\end{proof}

Thus, for isotropic factors, matrix independence never makes the contrast grow
faster than free independence, and the two slopes agree in the limit of large
dimension; what distinguishes them for large $d$ is the remainder. (For
Gaussian factors, which are not of the form $UD$ with $D$ fixed and are not invertible
in the free limit, Proposition~\ref{random:prop:isotropic} does not apply; with entries of
variance $1/d$, $\frac12(\gamma_1-\gamma_d)\sim\frac14\ln d$ \cite{Newman86}.)

\section{Decay: convergence to the tracial state}\label{random:sec:decay}

A positive invertible $\rho$ has $\OLC(\rho)=0$ iff $\rho$ is a multiple of $1$.
For a positive definite density matrix this means that $\rho$ is the normalised
trace, the maximally mixed state; for a probability density on a compact group,
that it is the Haar measure. In this section the contrast tends to $0$.

\paragraph{The Hilbert metric.} For positive definite $x,y\in M_d(\C)$ let
$M(x/y)=\inf\{t>0:\ x\le ty\}$ and $d_H(x,y)=\ln\bigl(M(x/y)M(y/x)\bigr)$, the
Hilbert projective metric of the cone of positive semidefinite matrices. Since
$M(\rho/1)=\lambda_{\max}(\rho)$ and $M(1/\rho)=\lambda_{\min}(\rho)^{-1}$,
\begin{equation}\label{random:eq:hilbert}
d_H(\rho,1)=\ln\kappa(\rho)=2\,\OLC(\rho).
\end{equation}
Let $\Phi$ be a unital positive linear map on $M_d(\C)$; it maps positive
definite matrices to positive definite ones, since $x\ge\lambda_{\min}(x)1$
implies $\Phi(x)\ge\lambda_{\min}(x)1$. Let
$\Delta(\Phi)=\sup\{d_H(\Phi(x),\Phi(y)):\ x,y>0\}\in[0,\infty]$ be the
projective diameter of its image; $\Delta(\Phi)<\infty$ if, for instance,
$\Phi$ maps every nonzero positive semidefinite matrix to a positive definite
one.

\begin{proposition}\label{random:prop:BH}
Let $\Phi$ be a unital positive linear map on $M_d(\C)$. Then for every positive definite $\rho$,
\[
\OLC\bigl(\Phi(\rho)\bigr)\le\tanh\Bigl(\frac{\Delta(\Phi)}4\Bigr)\,\OLC(\rho).
\]
In particular, if $\Delta(\Phi)<\infty$, then
$\OLC(\Phi^N(\rho))\le\tanh(\Delta(\Phi)/4)^N\OLC(\rho)$ decays geometrically.
\end{proposition}

\begin{proof}
By the Birkhoff--Hopf theorem \cite{Birkhoff57,Bushell73b,EN95} (see
\cite{RKW11} for the cone of positive semidefinite matrices),
$d_H(\Phi(x),\Phi(y))\le\tanh(\Delta(\Phi)/4)\,d_H(x,y)$. Take $y=1$ and use
$\Phi(1)=1$ and \eqref{random:eq:hilbert}.
\end{proof}

For unital trace-preserving positive maps on $M_d(\C)$, the Schatten contrasts $C_p$,
$p\ge1$, are also nonincreasing, so they are Lyapunov functions for
unital dynamics. This is a form of Jensen's trace inequality: let $\Phi$ be unital,
positive and trace-preserving, $\rho=\sum_i\lambda_iP_i$ positive definite with spectral
projections $P_i$, and $(e_k)$ an orthonormal basis with $\Phi(\rho)e_k=\beta_ke_k$. Then
$\beta_k=\sum_iw_{ki}\lambda_i$ with $w_{ki}=\langle e_k,\Phi(P_i)e_k\rangle\ge0$ and
$\sum_iw_{ki}=\langle e_k,\Phi(1)e_k\rangle=1$, so for every convex $f$ on
$(0,\infty)$, $f(\beta_k)\le\sum_iw_{ki}f(\lambda_i)=\langle e_k,\Phi(f(\rho))e_k\rangle$;
summing over $k$ gives $\tau(f(\Phi(\rho)))\le\tau(\Phi(f(\rho)))=\tau(f(\rho))$. With
$f(t)=t^p$, $p\ge1$, and $f(t)=t^{-p}$ this gives $\tau(\Phi(\rho)^{\pm p})\le\tau(\rho^{\pm p})$,
hence $C_p(\Phi(\rho))\le C_p(\rho)$. Near the tracial state the $C_p$ measure a different distance:
$C_p(\rho)=\frac p2\|\ln\rho-\tau(\ln\rho)\|_2^2\,(1+o(1))$ is a squared
$L^2(\tau)$-distance, whereas $\OLC(\rho)=\frac12(\max-\min)$ of the spectrum of
$\ln\rho$ is of first order.

\paragraph{The circle and the visibility of fringes.} Let $\T=\R/2\pi\mathbb Z$
with normalised Lebesgue measure, $\mu$ a probability measure on $\T$ with
Fourier coefficients $\hat\mu(k)=\int e^{-ikx}\dd\mu(x)$, and
$f_0=1+\sum_{k\ne0}a_ke^{ikx}$ a continuous strictly positive probability
density with $\sum_k|a_k|<\infty$. The density of $\mu^{*N}*f_0\,\dd x$ is
$f_N=1+\sum_{k\ne0}\hat\mu(k)^Na_ke^{ikx}$, a positive invertible element of the
finite von Neumann algebra $L^\infty(\T)$, and convolution with $\mu$ is a unital
trace-preserving positive map on it.

\begin{proposition}\label{random:prop:circle}
Let $K=\{k\ne0:\ a_k\ne0\}$, $r=\sup_{k\in K}|\hat\mu(k)|$, and suppose that
$|\hat\mu(k)|=r$ exactly for $k=\pm k_0\in K$, that
$r'=\sup\{|\hat\mu(k)|:\ k\in K\setminus\{\pm k_0\}\}<r$, and that $0<r<1$, with the convention $\sup\emptyset=0$ for $r'$ (i.e.\
$\mu$ is not supported on a coset of $\frac{2\pi}{k_0}\mathbb Z$; this is
automatic if $\mu$ is absolutely continuous). Then, as $N\to\infty$,
\[
\OMC(f_N)=\frac{\max f_N-\min f_N}{\max f_N+\min f_N}=2|a_{k_0}|\,r^N\bigl(1+o(1)\bigr),
\qquad\OLC(f_N)=2|a_{k_0}|\,r^N\bigl(1+o(1)\bigr).
\]
\end{proposition}

\begin{proof}
Since $f_N$ is real, $a_{-k_0}=\overline{a_{k_0}}$ and
$\hat\mu(-k_0)=\overline{\hat\mu(k_0)}$, so
$f_N=1+2\operatorname{Re}\bigl(\hat\mu(k_0)^Na_{k_0}e^{ik_0x}\bigr)+R_N$ with
$\|R_N\|_\infty\le r'^N\sum_k|a_k|$. The middle term has maximum
$2|a_{k_0}|r^N$ and minimum $-2|a_{k_0}|r^N$, attained at points of
$\T$, so $\max f_N-\min f_N=4|a_{k_0}|r^N+O(r'^N)$ and, since $r<1$,
$\max f_N+\min f_N=2+O(r^N)\to2$. Finally $\OLC=\operatorname{artanh}\OMC$ and
$\operatorname{artanh}u\sim u$.
\end{proof}

If $\mu$ itself has a continuous positive density with absolutely summable
Fourier series, one may take $f_0$ to be this density; then $f_N$ is the
density of $\mu^{*(N+1)}$. In optics, an interference pattern
$I(x)=I_0\bigl(1+V\cos(k_0x+\phi)\bigr)$ has Michelson contrast (visibility)
$V$, and a random translation $x\mapsto x-\theta$ with $\theta\sim\mu$
multiplies $V$ by $|\hat\mu(k_0)|$ \cite[Ch.~10]{BornWolf};
Proposition~\ref{random:prop:circle} says that $N$ independent perturbations reduce the
visibility by $|\hat\mu(k_0)|^N$, provided this harmonic dominates. (For
Gaussian phase noise, $|\hat\mu(k)|=e^{-\sigma^2k^2/2}$ and higher harmonics decay
faster, but in general they need not.) The decay of the contrast to $0$ is the
compact-group counterpart of the central limit theorem: on $\R$ convolution
powers spread out, on a compact group they converge to the Haar measure under
the condition of Kawada and It\^o \cite{KI40}, and the Michelson contrast
measures the distance.

\paragraph{Random unital channels.} Let $\Phi_1,\Phi_2,\dots$ be i.i.d.\ random
unital trace-preserving positive maps on $M_d(\C)$, $d\ge2$, and
$\rho_N=\Phi_N\cdots\Phi_1(\rho)$. The real space $V$ of traceless self-adjoint
matrices is invariant under every $\Phi_i$, and
$\rho_N-\tau(\rho)1=\Phi_N\cdots\Phi_1(\rho-\tau(\rho)1)$. We equip $V$ with the
operator norm; since $-\|x\|\le x\le\|x\|$ implies
$-\|x\|\le\Phi(x)\le\|x\|$, each $\Phi_i|_V$ is a contraction. Let
$\gamma=\lim_N\frac1N\ln\|\Phi_N\cdots\Phi_1|_V\|\in[-\infty,0]$ be the top
Lyapunov exponent of the restrictions to $V$, which exists almost surely and is
constant by \cite{FK60,Kingman68}.

\begin{proposition}\label{random:prop:channels}
\begin{enumerate}
\item For every positive definite $\rho$, almost surely
\[\limsup_N\frac1N\ln\OLC(\rho_N)\le\gamma.\]
If $\gamma<0$, the contrast decays
exponentially.
\item Suppose, moreover, that $\Phi_1|_V$ is almost surely invertible with
$\E\ln\|(\Phi_1|_V)^{-1}\|<\infty$, that the group generated by the support of
its law acts strongly irreducibly on $V$, and that $\gamma<0$. Then for every
positive definite $\rho$ that is not a multiple of $1$, almost surely
$\lim_N\frac1N\ln\OLC(\rho_N)=\gamma$.
\end{enumerate}
\end{proposition}

\begin{proof}
Normalise $\tau(\rho)=1$ and write $\rho_N=1+X_N$ with $X_N\in V$. If
$\eta=\|X_N\|<1$, the eigenvalues of $\rho_N$ lie in $[1-\eta,1+\eta]$ and
$\lambda_{\max}-\lambda_{\min}\ge\eta$ (since $X_N$ is traceless and nonzero,
its largest and smallest eigenvalues have opposite signs and one of them has
modulus $\eta$), whence, writing $\lambda_{\max}=1+a$, $\lambda_{\min}=1+b$ with $b\le0$ and
$a-b\ge\eta$, we get $\lambda_{\max}/\lambda_{\min}\ge1+\eta$ and
$\frac\eta4\le\frac12\ln(1+\eta)\le\OLC(\rho_N)\le\frac12\ln\frac{1+\eta}{1-\eta}\le2\eta$
for $\eta\le\frac12$. Thus $\ln\OLC(\rho_N)=\ln\|X_N\|+O(1)$ as soon as
$\|X_N\|\le\frac12$. (a) If $\gamma=0$ there is nothing to prove, since
$\OLC(\rho_N)\le\OLC(\rho)$ by Proposition~\ref{random:prop:BH}; if $\gamma<0$, use
$\|X_N\|\le\|\Phi_N\cdots\Phi_1|_V\|\,\|X_0\|$. (b) By Furstenberg's theorem for
strongly irreducible i.i.d.\ products \cite[Part~A, Corollary~III.3.4]{BL85},
$\frac1N\ln\|\Phi_N\cdots\Phi_1(X_0)\|\to\gamma$ almost surely for every fixed
$X_0\ne0$.
\end{proof}

For random mixed unitary channels $\Phi(x)=\sum_jw_jU_jxU_j^*$ with independent
Haar unitaries, a single map need not contract on $V$ (its norm on $V$ can be
$1$), but numerically $\gamma\approx-0.29$ in the example of
Section~\ref{random:sec:numerics}; the decay rate is computed by the standard QR
algorithm. The three regimes of decay---geometric with ratio $\tanh(\Delta/4)$,
exact exponential $2|a_{k_0}|r^N$, and the Lyapunov rate $e^{\gamma N}$---mirror
the three regimes of growth in Section~\ref{random:sec:compare}.

\section{Numerical illustrations}\label{random:sec:numerics}

The simulations and numerical values in this section are retained from the component manuscript. They illustrate the stated formulas; the large matrix and free-convolution simulations have not all been rerun for this consolidated version. The circle visibility ratios and the three-dimensional isotropic Lyapunov integrals were recomputed for this version and agree with the displayed rounded values. Values obtained by numerical inversion or quadrature are numerical evaluations, not exact arithmetic certificates.

\paragraph{Free products.} Let $X=UD$ with $U$ Haar and $D$ positive with
spectrum uniformly distributed on $[1,3]$, so that
$C_2(X)=\frac14\ln\frac{13}9=0.0919$. The free limit of
$\frac12\ln(\|A_N\|^2/m^N)$ can be computed exactly from the $S$-transform: the
right end of the support of $\mu_1^{\boxx N}$ is $1/\max_{w}\chi_N(w)$, the
maximum being taken up to the first critical point of
$\chi_N(w)=\frac w{1+w}S_{\mu_1}(w)^N$. Simulations use orthogonal matrices of
size $n$.

\begin{center}
\begin{tabular}{@{}rccc@{}}
\toprule
$N$ & numerical free value & $\frac12\ln(evN)$ & simulation, $n=700$\\
\midrule
8 & 1.044 & 0.919 & 1.030\\
64 & 1.978 & 1.958 & 1.885\\
1024 & 3.346 & 3.345 & ---\\
\bottomrule
\end{tabular}
\end{center}

\noindent At $N=1024$ the numerical free value of $\OLC(A_N)-NC_2(X)$ is $3.4537$,
against $3.4525$ from Theorem~\ref{random:thm:free}(a). In simulations with $n=500$ (eigenvalues of $D$ equally spaced in $[1,3]$, endpoints
included), $C_2(A_N)/N$ stays at the value $C_2(X)=0.0923$ for the discretised spectrum
(simulated range $0.0921$--$0.0923$) for all $N\le32$, in accordance with
Proposition~\ref{random:prop:C2}, while $C_1(A_N)/N$ increases and $C_4(A_N)/N$ and
$\OLC(A_N)/N$ decrease towards it; the increments of $C_1(A_N)-NC_2(X)$ and
$C_4(A_N)-NC_2(X)$ per doubling of $N$ approach $-\frac12\ln2$ and
$\frac14\ln2$. For Conjecture~\ref{random:conj:free} we computed fractional moments of
$\mu_1^{\boxx N}$ numerically, from the density of $\mu_1^{\boxx N}$ obtained by
inverting $\chi_N$; for $c_N(s)=\ln\tau(Y_N^s)-(s-1)\ln N$ at $N=1024$:

\begin{center}
\begin{tabular}{@{}lcccc@{}}
\toprule
$s$ & $\frac14$ & $\frac12$ & $\frac34$ & $\frac32$\\
\midrule
$c_{1024}(s)$ & 2.063 & 1.0855 & 0.4657 & $-0.7012$\\
$(s-1)\ln(vs)-\ln\Gamma(s+1)$ & 2.070 & 1.088 & 0.467 & $-0.703$\\
\bottomrule
\end{tabular}
\end{center}

\paragraph{Matrix products.} For $d=3$ and $X=UD$ with $U$ Haar orthogonal and
$D=\operatorname{diag}(1,2,3)$, the formula of
Proposition~\ref{random:prop:isotropic} gives $\gamma_1=0.7111$, $\gamma_3=0.4743$ and
$\frac12(\gamma_1-\gamma_3)=0.1184$, against the free slope $C_2(X)=0.1875$ for the
same $D$; computing $\OLC(A_N)$ through the accumulated QR factors (which agree
with $\ln s_1$ and $\ln s_3$ up to $O_P(1)$) gives $\frac1N\OLC(A_N)\approx0.118$ for $N=1600$, averaged over $200$
samples, and standard deviations $0.25\sqrt N$, $0.23\sqrt N$, $0.23\sqrt N$ for
$N=100,400,1600$.
For $D_d$ with eigenvalues equally spaced in $[1,3]$, the slope
$\frac12(\gamma_1-\gamma_d)$ is $0.0961$, $0.0922$, $0.0920$ for $d=10,100,1000$,
approaching the free slope $\frac14\ln\frac{13}9=0.0919$.

\paragraph{Decay.} On the circle, with $\mu$ the measure with atoms
$0,1,2.5$ of weights $0.5,0.3,0.2$ and $f_0=1+0.8\cos x+0.3\cos(2x+1)$, one has
$|\hat\mu(1)|=0.6248$, $|\hat\mu(2)|=0.4394$, and the ratio of $\OMC(f_N)$ to
$2|a_1|\,|\hat\mu(1)|^N$ is $0.9954$, $1.00003$, $1+1.3\cdot10^{-7}$ for
$N=5,10,20$ (the relative correction is $O(r'^N+(r'/r)^{2N})$; the term $(r'/r)^{2N}$
comes from the second harmonic, which takes equal values at the antipodal extrema of the
first; here $r^2<r'$, so the second term dominates for large $N$, although at $N=5$ the
first is still the larger). For mixed unitary channels on $M_3(\C)$ with two independent Haar
unitaries and weights $0.7,0.3$, a long QR run gives $\gamma=-0.29$ (a single
run of $400$ steps gives $-0.28$). Since $\|\rho_N-1\|$ falls below machine
precision after about $100$ steps, we propagate the traceless part
$\rho_N-1$, re-projecting onto the traceless subspace and renormalising at each
step; for
$\rho=3\operatorname{diag}(0.6,0.3,0.1)$, so that $\tau(\rho)=1$, this gives
$\frac1N\ln\|\rho_N-1\|=-0.284$ at $N=400$, which equals
$\frac1N\ln\OLC(\rho_N)+O(\frac1N)$ by the proof of
Proposition~\ref{random:prop:channels}.

\section{Open problems}

\begin{enumerate}[label=\textup{(\arabic*)}]
\item Prove Conjecture~\ref{random:conj:free}, i.e.\ the asymptotics of fractional
moments of $\mu^{\boxx N}$, at least the exponent $\frac12-\frac1p$ for $p<2$.
\item Non-identically distributed free factors: by Proposition~\ref{random:prop:C2},
$C_2$ is additive; is $C_p(A_N)-\sum_iC_2(X_i)=O(\ln N)$ under uniform bounds?
\item Products in a factor of type $\mathrm{II}_1$ that are neither free nor
tensor independent (for instance, stationary ergodic sequences of unitarily
conjugated elements): do the slopes $\lim C_p(A_N)/N$ exist, and when do they
coincide? Lemma~\ref{random:lem:holder} gives only $C_{p/N}(A_N)\le\sum_iC_p(X_i)$.
\item Fluctuations: for Haar-rotated $n\times n$ factors, with
$n$ and $N$ growing jointly, which regime interpolates between the Gaussian
$\sqrt N$ fluctuations of fixed $d$ and the deterministic $\frac12\ln N$ of the
free limit?
\end{enumerate}

\chapter{Symmetric cones and Lorentz boosts}\label{ch:cones}
\section{Introduction}

For a positive invertible operator $X$ the operator Michelson contrast and the
logarithmic contrast are
\[
\OMC(X)=\frac{\lambda_{\max}(X)-\lambda_{\min}(X)}{\lambda_{\max}(X)+\lambda_{\min}(X)},\qquad
\OLC(X)=\operatorname{artanh}\OMC(X)=\tfrac12\ln\frac{\lambda_{\max}(X)}{\lambda_{\min}(X)}
\]
\cite{ANN24,NN26}, extended to invertible operators through $|X|$. In \cite{NN26} it was
shown that $\OLC(XY)\le\OLC(X)+\OLC(Y)$, that $\OLC$ is the projectivised Thompson distance
from $[I]$, and that for positive $n\times n$ matrices equality forces commutation
exactly for $n\le3$. This chapter develops four related questions raised by these results.

\emph{What is the threshold $n=4$ a threshold of?} We replace $M_n(\C)$ by an
arbitrary Euclidean Jordan algebra and the product $XY$ by the quadratic
representation $P(x^{1/2})y$, which is $X^{1/2}YX^{1/2}$ for matrices. Subadditivity
persists (Theorem~\ref{cones:thm:sub}); equality is characterised by common primitive
idempotents under the extreme spectral idempotents (Theorem~\ref{cones:thm:eq}); and
equality forces operator commutativity if and only if the rank is at most $3$ or the
algebra is associative (Theorem~\ref{cones:thm:rank}). So the threshold is one of the rank,
not of the dimension (for pure tensor products, Proposition~\ref{core:tensor-equality} already shows
that it is not one of the total dimension): it holds on Lorentz cones of every dimension and on the
$27$-dimensional exceptional cone, and it extends to JBW-algebras
(Theorem~\ref{cones:thm:jbw}).

\emph{What is the defect of the Schatten contrasts a defect of?} The Schatten
contrasts are subadditive only up to a defect (Theorem~\ref{cones:thm:defect}(a),(c)). On a symmetric cone the defect of every spectral contrast
$C_\Phi$ is $\ln(\Phi(\mathbf 1)/\Phi(\epsilon_1))$ for $\Phi$ on $\R^{\rk V}$, so for $C_p$ it
is $\frac{\ln\rk V}p$ (and $\frac{\ln n}p$ for positive $n\times n$ matrices), and $\OLC$ is
the only defect-free one (Theorem~\ref{cones:thm:defect}).
Both basic constants, the non-commutativity threshold and the
defect, are thus governed by the rank.

\emph{Is there kinematics behind the contrast?} On the Lorentz cone $\OMC$ is the speed
of a four-velocity and $\OLC$ its rapidity, the Thompson geometry of the unit
hyperboloid is hyperbolic geometry, and subadditivity of $\OLC$ is the triangle
inequality of $H^n$ with equality exactly for collinear boosts in the same direction
(Theorem~\ref{cones:thm:hyp}); for $2\times2$ matrices this is the spinor and vector
representation of $\mathrm{SL}(2,\C)$ (Remark~\ref{cones:rem:spinor}), and it sharpens a
remark of \cite{NN26}. In every non-trivial finite-dimensional representation of $\mathrm{SL}(2,\C)$, with an
admissible inner product, equality for two pure boosts forces equal axes (Theorem~\ref{cones:thm:reps}); the
non-commuting $4\times4$ equality case of \cite{NN26} therefore has no interpretation
as a pair of Lorentz boosts, not even for Dirac spinors, and is realised instead by
orthogonal blocks, i.e.\ direct sums (the example after Corollary~\ref{cones:cor:dirac}). The rotation $k$ of the polar
decomposition $P(x^{1/2})P(y^{1/2})=P(w^{1/2})k$ (Lemma~\ref{cones:lem:polar}) is the Wigner rotation, and its angle is the area of the triangle whose
sides are the three contrasts in the subadditivity inequality (Remark~\ref{cones:rem:wigner}).

\emph{What does the contrast see of the Jordan product?} We ask for which given
contrasts $u,v$ of $x,y$ the contrast of $x\circ y$ (for matrices, of $XY+YX$) is
subadditive. Strang's bounds \cite{Strang62} hold in every JB-algebra
by the Shirshov--Cohn theorem (Lemma~\ref{cones:lem:strang}), and on a non-associative Euclidean
Jordan algebra the region is the ellipse $u^2+uv+v^2\le1$ exactly for spin factors and the
triangle $u+v\le1$ otherwise, including the real, complex, quaternionic and exceptional cases
(Theorem~\ref{cones:thm:jordan}); in particular, for positive $n\times n$ matrices the threshold
changes between $n=2$ and $n=3$ (Corollary~\ref{cones:cor:jordan-cases}). The multiplicative commutator of two spinor boosts is
conjugate to a rotation for all angles $\theta\in(0,\pi)$ exactly when $u^2+v^2<1$, the positivity threshold of the Jordan
product (Remark~\ref{cones:rem:groupcomm}).

Finally, traces on $V$ are characterised by a contrast inequality involving the
symmetries $s$ of $V$ (Theorem~\ref{cones:thm:trace}).

\subsection*{What was known and what is new}
Theorem~\ref{cones:thm:sub} is folklore: it is the triangle inequality for the
Thompson/Hilbert metric on a symmetric cone. Theorem~\ref{cones:thm:eq} is close to, but not
the same as, the description of unique Thompson geodesics and of midpoint sets on
symmetric cones by Lemmens and Roelands \cite{LR15,LR17}; we have not found it stated.
We have not found the rank threshold (Theorem~\ref{cones:thm:rank}) in the literature; it is a
short consequence of Theorem~\ref{cones:thm:eq}. The identification of rapidity with
hyperbolic distance, and of metrics on boosts and on qubit density matrices with the
rapidity metric, is classical \cite{RS04,KL13,Kim11}; Theorem~\ref{cones:thm:hyp} records it
in terms of the contrast. Theorem~\ref{cones:thm:reps} is a direct consequence of the $KAK$
decomposition (Lemma~\ref{cones:lem:kak}) and of the $2\times2$ case; we have not found it
stated, but it is elementary. The spectral identity in Lemma~\ref{cones:lem:polar} is stated in \cite[Lemma~3.1]{TJG20}; it follows at once from the polar decomposition of $G(\Omega)$ \cite[Thm.~III.5.1]{FK94}, \cite[\S2]{Lim00}.
In Theorem~\ref{cones:thm:defect}, (a) is proved via the min--max principle, and (b)--(d) are
elementary statements about eigenvalue vectors in $\R^{\rk V}$; its key
inequality and the norm inequality underlying the Schatten case of (a) are in \cite{TJG20}. The Wigner angle as the area of the velocity triangle is classical \cite{RS04};
Remark~\ref{cones:rem:wigner} only identifies the triangle. Theorem~\ref{cones:thm:jbw} combines
the proof of Theorem~\ref{cones:thm:rank} with the structure theory of JBW-algebras
\cite{HOS84}. In Section~\ref{cones:sec:jordan}, Lemma~\ref{cones:lem:strang} is Strang's theorem
\cite{Strang62,AB85} transferred by the Shirshov--Cohn theorem, the positivity criterion is
Gustafson's \cite{Gus68}, and the additive commutator bound and the trace formula of
Remark~\ref{cones:rem:groupcomm} are classical \cite{Strang62,Bea83}. Theorem~\ref{cones:thm:trace} is a
Jordan analogue of characterisations of traces on $C^*$-algebras by inequalities
\cite{Bik12,BM26}. All proofs are short and
use only the Peirce calculus and standard structure theory; the finite-dimensional
statements were also checked numerically on $\mathrm{Herm}_n(\R)$, $\mathrm{Herm}_n(\C)$ and
spin factors. Apart from standard facts and the classical results cited in the proofs, the chapter relies only on \cite{ANN24,NN26}.

\section{Setting}

Let $V$ be a Euclidean Jordan algebra with unit $e$, rank $r$, trace $\tr$ and
inner product $\langle x,y\rangle=\tr(x\circ y)$, and let $\Omega$ be the interior
of the cone of squares \cite{FK94}. Every $x\in V$ has a spectral decomposition
$x=\sum_{i=1}^r\lambda_ic_i$ with respect to a Jordan frame $(c_i)$; we write
$\lambda_{\max}(x)$, $\lambda_{\min}(x)$ for the extreme eigenvalues and
$E_{\max}(x)$, $E_{\min}(x)$ for the corresponding spectral idempotents (sums of
the $c_i$ with $\lambda_i=\lambda_{\max}$, resp.\ $\lambda_{\min}$). The quadratic
representation is $P(x)=2L(x)^2-L(x^2)$. We use the following standard facts
(references are given item by item):
\begin{enumerate}[label=(F\arabic*)]
\item for invertible $a$, $P(a)$ is invertible, self-adjoint, $P(a)\Omega=\Omega$,
$P(a)e=a^2$, $P(a)^{-1}=P(a^{-1})$, and $(P(a)y)^{-1}=P(a^{-1})y^{-1}$
\cite[Ch.~II--III]{FK94};
\item $x\in\Omega$ iff all $\lambda_i>0$ \cite[Ch.~III]{FK94}, and
$\lambda_{\max}(z)=\max\{\langle z,w\rangle:\ w\ge0,\ \tr w=1\}$, the maximum
being attained at a primitive idempotent. (Proof: $\langle\lambda_{\max}(z)e-z,w\rangle\ge0$
by self-duality of the closed cone \cite[Thm.~III.2.1]{FK94}; equality holds at the top
frame elements of $z$; and if $w=\sum_i\mu_id_i$ is a spectral decomposition, then
$\mu_i\ge0$ and $\sum_i\mu_i=\tr w=1$, so $\langle z,w\rangle\le\max_i\langle z,d_i\rangle$.)
\item for $u,w\ge0$, $\langle u,w\rangle=0$ implies $u\circ w=0$ \cite[Prop.~6]{GST04};
\item Peirce decomposition $V=V(c,1)\oplus V(c,\frac12)\oplus V(c,0)$ for an idempotent
$c$, with $V(c,1)\circ V(c,0)=0$; for orthogonal idempotents $f,g$,
$V(f,1)\circ V(g,1)=0$;
\item two elements $x,y$ operator commute, i.e.\ $[L(x),L(y)]=0$, if and only if
they have a common Jordan frame
\cite[Lemma~X.2.2]{FK94};
\item if $V$ is simple of rank $\ge2$ and $c_i\perp c_j$ are primitive, then
$V_{ij}=V(c_i,\frac12)\cap V(c_j,\frac12)\ne0$, and there is $u\in V_{ij}$ with
$u^2=c_i+c_j$; for such $u$, $P(u)c_i=c_j$. (For $u\in V_{ij}$ one has
$u^2=\frac12|u|^2(c_i+c_j)$ \cite[Prop.~IV.1.4]{FK94}, and then
$P(u)c_i=2u\circ(u\circ c_i)-u^2\circ c_i=u^2-c_i=c_j$.)
\end{enumerate}
Two consequences of (F4)--(F5) are used below.
(F7) An element $y$ operator commutes with an idempotent $c$ iff
$y\in V(c,1)\oplus V(c,0)$. Indeed, if $y=y_1+y_0$ with $y_\nu\in V(c,\nu)$, the union of
Jordan frames of $y_1$ in $V(c,1)$ and of $y_0$ in $V(c,0)$ is a common frame;
conversely, in a common frame $c$ is a sum of frame elements, and each frame element
lies in $V(c,1)$ or in $V(c,0)$.
(F8) If $x$ and $y$ operator commute, then $y$ operator commutes with every spectral
idempotent of $x$, since each is a sum of elements of the common frame.
For a primitive idempotent $c$ and an idempotent $E$ we write $c\le E$ if
$c\in V(E,1)$, i.e.\ $E\circ c=c$.

\begin{definition}
For $x\in\Omega$ put $\kappa(x)=\lambda_{\max}(x)/\lambda_{\min}(x)$,
$\OLC(x)=\frac12\ln\kappa(x)$ and $\OMC(x)=\tanh\OLC(x)=\frac{\kappa-1}{\kappa+1}$.
\end{definition}

As for matrices \cite{NN26}: $\OLC(\lambda x)=\OLC(x)$, $\OLC(x^{-1})=\OLC(x)$,
$\OLC(x^t)=|t|\OLC(x)$, $\OLC(x)=0$ iff $x\in\R_+e$, and
$\OLC(x)=\delta([e],[x])$ for the projectivised Thompson metric
$\delta([x],[y])=\inf_{\lambda>0}d_T(x,\lambda y)=\inf_{\lambda>0}\|\log P(x^{-1/2})(\lambda y)\|_\infty$,
which is one half of Hilbert's projective metric $d_H$ on $\Omega$
(the proof of Proposition~\ref{core:metric} applies verbatim). For $V=\mathrm{Herm}_n(\C)$,
$P(x^{1/2})y=X^{1/2}YX^{1/2}$.

\section{Subadditivity and the equality case}

\begin{lemma}\label{cones:lem:face}
Let $z\in V$, $E=E_{\max}(z)$, and $w\ge0$ with $\tr w=1$. Then
$\langle z,w\rangle=\lambda_{\max}(z)$ iff $w\in V(E,1)$. In particular, for a
primitive idempotent $c$: $\langle z,c\rangle=\lambda_{\max}(z)$ iff $c\le E$, and
then $z\circ c=\lambda_{\max}(z)c$.
\end{lemma}

\begin{proof}
$u=\lambda_{\max}(z)e-z\ge0$ and $\langle u,w\rangle=\lambda_{\max}\tr w-\langle z,w\rangle$.
If this vanishes, $u\circ w=0$ by (F3). Write $u=\sum_i\nu_ic_i$ with $(c_i)$ a
frame of $z$ and $\nu_i=\lambda_{\max}-\lambda_i\ge0$. In the Peirce decomposition
$V=\bigoplus_{i\le j}V_{ij}$ of this frame, $L(u)$ acts on $V_{ij}$ as multiplication
by $\frac12(\nu_i+\nu_j)$ \cite[Ch.~IV]{FK94}, so
$\ker L(u)=\bigoplus\{V_{ij}:\nu_i=\nu_j=0\}=V(E,1)$, and $w\in V(E,1)$. Conversely,
if $w\in V(E,1)$ then $u\circ w=0$, so $\langle u,w\rangle=\tr(u\circ w)=0$. If $c\le E$, then
$z=\lambda_{\max}E+z'$ with $z'\in V(E,0)$, so
$z\circ c=\lambda_{\max}E\circ c+z'\circ c=\lambda_{\max}c$ by (F4).
\end{proof}

\begin{lemma}\label{cones:lem:Pc}
Let $x\in\Omega$ and $c\le E_{\max}(x)$ primitive. Then
$P(x^{1/2})c=\lambda_{\max}(x)\,c$.
\end{lemma}

\begin{proof}
Write $x^{1/2}=\alpha E+r$ with $\alpha=\lambda_{\max}(x)^{1/2}$, $E=E_{\max}(x)$,
$r\in V(E,0)$. Since $c\in V(E,1)$, (F4) gives $L(r)c=0$, $L(r^2)c=0$,
$L(E)c=c$; hence $P(r)c=0$ and the mixed term
$L(E)L(r)c+L(r)L(E)c-L(E\circ r)c=0$, so $P(x^{1/2})c=\alpha^2P(E)c=\alpha^2c$.
\end{proof}

\begin{theorem}\label{cones:thm:sub}
For $x,y\in\Omega$,
\[
\lambda_{\max}(P(x^{1/2})y)\le\lambda_{\max}(x)\lambda_{\max}(y),\qquad
\lambda_{\min}(P(x^{1/2})y)\ge\lambda_{\min}(x)\lambda_{\min}(y),
\]
hence $\OLC(P(x^{1/2})y)\le\OLC(x)+\OLC(y)$ and
$\OMC(P(x^{1/2})y)\le\frac{\OMC(x)+\OMC(y)}{1+\OMC(x)\OMC(y)}$.
\end{theorem}

\begin{proof}
$\lambda_{\min}(y)e\le y\le\lambda_{\max}(y)e$ and $P(x^{1/2})$ preserves the order
(F1), so $\lambda_{\min}(y)x\le P(x^{1/2})y\le\lambda_{\max}(y)x$, and
$\lambda_{\min}(x)e\le x\le\lambda_{\max}(x)e$. Equivalently: $\OLC$ is the
distance from $[e]$ in $\delta_T$, $P(x^{-1/2})$ is a $\delta_T$-isometry mapping
$[x]$ to $[e]$ and $[P(x^{1/2})y]$ to $[y]$, and the inequality is the triangle
inequality.
\end{proof}

\begin{theorem}[Equality]\label{cones:thm:eq}
For $x,y\in\Omega$, $\OLC(P(x^{1/2})y)=\OLC(x)+\OLC(y)$ if and only if there are
primitive idempotents $c,c'$ with
\[
c\le E_{\max}(x),\quad c\le E_{\max}(y),\qquad c'\le E_{\min}(x),\quad c'\le E_{\min}(y).
\]
\end{theorem}

\begin{proof}
By Theorem~\ref{cones:thm:sub}, equality for $\OLC$ is equivalent to equality in both
eigenvalue bounds. \emph{Top.} By (F2) and self-adjointness,
$\lambda_{\max}(P(x^{1/2})y)=\max_c\langle y,P(x^{1/2})c\rangle$ over primitive $c$.
For each state $w$, $w'=P(x^{1/2})w\ge0$ has $\tr w'=\langle x,w\rangle\le\lambda_{\max}(x)$
and $\langle y,w'\rangle\le\lambda_{\max}(y)\tr w'$. So equality forces a primitive $c$
with $\langle x,c\rangle=\lambda_{\max}(x)$, i.e.\ $c\le E_{\max}(x)$
(Lemma~\ref{cones:lem:face}), and then $w'=\lambda_{\max}(x)c$ (Lemma~\ref{cones:lem:Pc}) and
$\langle y,c\rangle=\lambda_{\max}(y)$, i.e.\ $c\le E_{\max}(y)$. Conversely, for such
$c$, $\lambda_{\max}(P(x^{1/2})y)\ge\langle y,P(x^{1/2})c\rangle=\lambda_{\max}(x)\lambda_{\max}(y)$.
\emph{Bottom.} By (F1), $\lambda_{\min}(P(x^{1/2})y)^{-1}=\lambda_{\max}(P(x^{-1/2})y^{-1})$
and $E_{\max}(x^{-1})=E_{\min}(x)$; apply the top case to $x^{-1},y^{-1}$.
\end{proof}

\begin{remark}
For $V=\mathrm{Herm}_n(\C)$ the same condition (common eigenvectors in
$E_{\max}(X)\cap E_{\max}(Y)$ and in $E_{\min}(X)\cap E_{\min}(Y)$) is also
necessary and sufficient for $\kappa(XY)=\kappa(X)\kappa(Y)$: necessity is
Lemma~\ref{core:alignment}, and sufficiency follows by evaluating $XY$ and
$(XY)^{-1}$ on the common eigenvectors. So Lemma~\ref{core:alignment} is in fact an
equivalence, and the products $XY$ and $X^{1/2}YX^{1/2}$ have the same equality
case.
\end{remark}

\section{The rank threshold}

\begin{theorem}\label{cones:thm:rank}
For a Euclidean Jordan algebra $V$ the following are equivalent:
\begin{enumerate}
\item for all $x,y\in\Omega$, $\OLC(P(x^{1/2})y)=\OLC(x)+\OLC(y)$ implies that $x$
and $y$ operator commute;
\item $\rk V\le3$, or $V$ is associative ($V\cong\R^r$).
\end{enumerate}
\end{theorem}

\begin{proof}
If $V$ is associative, all elements operator commute. Let $\rk V\le3$, assume
equality, and take $c,c'$ as in Theorem~\ref{cones:thm:eq}; we may assume $x\notin\R e$.
By Lemma~\ref{cones:lem:face}, $y\circ c=\lambda_{\max}(y)c$ and $y\circ c'=\lambda_{\min}(y)c'$,
so (by the Peirce decomposition) $y$ has no components in $V(c,\frac12)$ and
$V(c',\frac12)$. As $\rk E_{\max}(x)+\rk E_{\min}(x)\le3$, at least one of them is
primitive. If both are, $c=E_{\max}(x)$, $c'=E_{\min}(x)$, and $g=e-c-c'$ is $0$
or primitive; then $V(c,0)=V(c'+g,1)=\R c'\oplus V_{c'g}\oplus\R g$ (with $V_{c'g}=0$ if $g=0$),
the component in $V_{c'g}\subset V(c',\frac12)$ is excluded, and so
$y\in V(c,1)\oplus V(c',1)\oplus V(g,1)=\R c\oplus\R c'\oplus\R g$,
which lies in the associative span of the spectral idempotents of $x$. If, say,
$E_{\max}(x)=c$ and $E_{\min}(x)=e-c$ has rank $2$, then $x=\lambda_{\max}c+\lambda_{\min}(e-c)$
and $y=\lambda_{\max}(y)c+y_0$ with $y_0\in V(c,0)$; a Jordan frame $(d_1,d_2)$ of
$y_0$ in $V(c,0)$ gives a common frame $(c,d_1,d_2)$ of $x$ and $y$. The case
$E_{\min}(x)$ primitive is symmetric. By (F5), $x$ and $y$ operator commute.

Now let $V$ be non-associative with $\rk V\ge4$. Some simple component $W$ has
rank $\ge2$. Choose a Jordan frame of $V$ containing primitive $c_2,c_3\in W$ and
two further elements $c_1,c_4$, and let $c_5,\dots,c_r$ be the rest. By (F6) pick
$u\in V_{23}$ with $u^2=c_2+c_3$. Put
\[
x=4c_1+3c_2+c_3+c_4+2\textstyle\sum_{k\ge5}c_k,\qquad
y=4c_1+2(c_2+c_3)+u+c_4+2\textstyle\sum_{k\ge5}c_k .
\]
Here $y'=2(c_2+c_3)+u$ lies in the rank-two subalgebra $V(c_2+c_3,1)$, where $u$ has
eigenvalues $\pm1$, so $y'$ has eigenvalues $3,1$ and $y\in\Omega$. All summands of
$x^{1/2}$ and $y$ lie in Peirce $1$-spaces of the orthogonal idempotents
$c_1$, $c_2+c_3$, $c_4$, $c_k$, which multiply to zero (F4); hence
$P(x^{1/2})y=16c_1+P(x'^{1/2})y'+c_4+4\sum_{k\ge5}c_k$ with $x'=3c_2+c_3$, and by
Theorem~\ref{cones:thm:sub} the spectrum of $P(x'^{1/2})y'$ lies in $[1,9]$. Thus
$\lambda_{\max}=16=4\cdot4$ and $\lambda_{\min}=1=1\cdot1$: equality holds. But $c_2$
is the spectral idempotent of $x$ for the simple eigenvalue $3$, and $y$ has the
nonzero component $u\in V(c_2,\frac12)$, so $y$ does not operator commute with
$c_2$ by (F7), hence not with $x$ by (F8).
\end{proof}

\begin{corollary}\label{cones:cor:cases}
\begin{enumerate}
\item $V=\mathrm{Herm}_n(\mathbb F)$, $\mathbb F=\R,\C,\mathbb H$: equality forces
commutativity iff $n\le3$. For $\mathbb F=\C$ this is Theorem~\ref{core:dimension}, by the
Remark after Theorem~\ref{cones:thm:eq}.
\item The Albert algebra $\mathrm{Herm}_3(\mathbb O)$ ($\dim27$, rank $3$): equality
always forces operator commutativity.
\item Spin factors $V_n=\R\times\R^n$ (Lorentz cones, any $n$, rank $2$): equality
forces operator commutativity; a direct sum $V_m\oplus V_n$ (rank $4$) does not when
$\max(m,n)\ge2$. ($V_1\cong\R^2$ is associative, so $V_1\oplus V_1\cong\R^4$ is an exception.)
\end{enumerate}
\end{corollary}

So the threshold of \cite{NN26} is a threshold in the rank, not in the dimension; this
complements Proposition~\ref{core:tensor-equality}, where it is shown not to be one of the total dimension
for pure tensor products.
The obstruction in the construction is a rank-two Peirce subalgebra $V(c_2+c_3,1)$
orthogonal to idempotents $c\le E_{\max}$ and $c'\le E_{\min}$; fitting these three
pieces needs rank $\ge4$. (The algebra itself may be simple, e.g.\ $\mathrm{Herm}_4(\R)$.)

\section{Polar decomposition and reciprocity}\label{cones:sec:polar}

Let $G(\Omega)=\{g\in\mathrm{GL}(V):g\Omega=\Omega\}$ and let $\mathrm{Aut}(V)$ be the
automorphism group of the Jordan algebra. Recall that $P(a)\in G(\Omega)$ and
$P(a)^*=P(a)$ for $a\in\Omega$ (F1), and that the stabiliser of $e$ in $G(\Omega)$ is
$\mathrm{Aut}(V)$ \cite[Prop.~I.1.9, Prop.~I.4.3, Thm.~III.5.1]{FK94}. (Directly: an
element $k$ of the stabiliser is orthogonal and preserves the characteristic function
$\varphi$ of $\Omega$; as $\nabla\log\varphi(z)$ is a negative multiple of $z^{-1}$, this gives
$k(z^{-1})=(kz)^{-1}$; differentiating gives $P(kz)=kP(z)k^{-1}$, and expanding at
$z=e+ta$ gives $L(ka)=kL(a)k^{-1}$.) Automorphisms are orthogonal and preserve spectra.

\begin{lemma}[Reciprocity]\label{cones:lem:polar}
For $x,y\in\Omega$ put $g=P(x^{1/2})P(y^{1/2})$ and $w=ge=P(x^{1/2})y$. Then
\[
g=P(w^{1/2})\,k\qquad\text{with }k\in\mathrm{Aut}(V),
\]
and $P(y^{1/2})x=k^{-1}w$. In particular $P(x^{1/2})y$ and $P(y^{1/2})x$ have the same
eigenvalues, with multiplicities.
\end{lemma}

\begin{proof}
$k:=P(w^{-1/2})g\in G(\Omega)$ and $ke=P(w^{-1/2})w=e$, so $k\in\mathrm{Aut}(V)$ and
$k^*=k^{-1}$. Since each $P(\cdot)$ is self-adjoint, $g^*=P(y^{1/2})P(x^{1/2})$, and also
$g^*=k^*P(w^{1/2})=k^{-1}P(w^{1/2})$; hence
$P(y^{1/2})x=g^*e=k^{-1}w$.
\end{proof}

The spectral identity is stated in \cite[Lemma~3.1]{TJG20}, with the reference
``([17], Corollary~9)'', where [17] is \cite{Lim00}. We could not locate such a statement in
\cite{Lim00} (its numbered results are Propositions~2.1, 2.2, 2.4--2.7, Lemma~2.3 and
Theorem~2.8, with a common counter); its \S2 recalls the
polar decomposition $g=P(x)k$ of $G(\Omega)$, from which the identity follows as above. For $V=\mathrm{Herm}_n(\C)$ it is the similarity of
$X^{1/2}YX^{1/2}$ and $Y^{1/2}XY^{1/2}$. What we use below is the factorisation itself:
on the Lorentz cone $k$ is the Wigner rotation (Remark~\ref{cones:rem:wigner}).

\section{Spectral contrasts: the defect is the rank}\label{cones:sec:spectral}

Let $r=\rk V\ge2$, let $\lambda(x)\in\R^r$ be the vector of eigenvalues of $x$ (with
multiplicities), and let $\Phi$ be a symmetric gauge function on $\R^r$, i.e.\ a norm
invariant under permutations and sign changes of coordinates. With $\mathbf 1=(1,\dots,1)$
and $\epsilon_1=(1,0,\dots,0)$ put
\[
\kappa_\Phi(x)=\Phi(\lambda(x))\,\Phi(\lambda(x^{-1})),\qquad
C_\Phi(x)=\tfrac12\ln\frac{\kappa_\Phi(x)}{\kappa_\Phi(e)},\qquad
d_\Phi=\ln\frac{\Phi(\mathbf 1)}{\Phi(\epsilon_1)}\qquad(x\in\Omega).
\]
For $\Phi=\|\cdot\|_\infty$, $C_\Phi=\OLC$ and $d_\Phi=0$. For $\Phi=\|\cdot\|_p$, $p\ge1$,
one obtains the \emph{Schatten contrasts}
\[
C_p(x)=\frac1{2p}\ln\bigl(\tau(x^p)\,\tau(x^{-p})\bigr),\qquad\tau=\tfrac1r\tr,\qquad d_p=\frac{\ln r}p ;
\]
we use the same formula for $0<p<1$. All these quantities depend only on $\lambda(x)$,
and every vector in $(0,\infty)^r$ is $\lambda(x)$ for some $x\in\Omega$ (take
$x=\sum_i\lambda_ic_i$ for a Jordan frame). Hence a statement about a single element of
$\Omega$ that depends only on its spectrum is a statement about vectors in $(0,\infty)^r$,
the same as for positive $r\times r$ matrices; the product law does not depend only on
spectra.

We use three elementary properties of a symmetric gauge function $\Phi$ on $\R^r$.
\begin{enumerate}[label=(G\arabic*)]
\item If $|v_i|\le|w_i|$ for all $i$, then $\Phi(v)\le\Phi(w)$. (By sign invariance we may
assume $v,w\ge0$, and it suffices to change one coordinate, say the first (by permutation invariance),
and then to apply this coordinate by coordinate; so let $v=w$ except
$0\le v_1\le w_1$, and $w_1>0$ (otherwise $v=w$). With $w'=(-w_1,w_2,\dots,w_r)$ and $t=\frac12(1+v_1/w_1)\in[0,1]$ we
have $v=tw+(1-t)w'$, so $\Phi(v)\le t\Phi(w)+(1-t)\Phi(w')=\Phi(w)$.)
\item $\Phi(\epsilon_1)\|v\|_\infty\le\Phi(v)\le\Phi(\mathbf 1)\|v\|_\infty$. (By (G1) and
permutation invariance, since $\|v\|_\infty\epsilon_j\le|v|\le\|v\|_\infty\mathbf 1$
coordinatewise, with $j$ an index where $|v_j|=\|v\|_\infty$.)
\item If $v,w\ge0$ and $\theta\in(0,1)$, then $\Phi(v^\theta w^{1-\theta})\le\Phi(v)^\theta\Phi(w)^{1-\theta}$,
with coordinatewise powers and product. (If $v=0$ or $w=0$ both sides vanish; otherwise,
by homogeneity, let $\Phi(v)=\Phi(w)=1$. The weighted arithmetic--geometric mean inequality
gives $v^\theta w^{1-\theta}\le\theta v+(1-\theta)w$ coordinatewise, so by (G1) and the
triangle inequality $\Phi(v^\theta w^{1-\theta})\le\theta+(1-\theta)=1$.)
\end{enumerate}

We use the min--max principle for Euclidean Jordan algebras \cite{Hirz70}, in the form:
if $u\le v$ in $V$, then $\lambda_j^\downarrow(u)\le\lambda_j^\downarrow(v)$ for all $j$, where
$\lambda^\downarrow$ lists the eigenvalues in nonincreasing order
\cite[Prop.~2.3]{TJG20}, \cite{GT11}. (For a direct sum, the eigenvalue lists of the
components are merged, which preserves the coordinatewise order.)

\begin{theorem}\label{cones:thm:defect}
Let $\Phi$ be a symmetric gauge function on $\R^r$, $r=\rk V\ge2$.
\begin{enumerate}
\item For all $x,y\in\Omega$, $C_\Phi(P(x^{1/2})y)\le C_\Phi(x)+C_\Phi(y)+d_\Phi$.
\item For $x\in\Omega$ and $s\ge1$, $s\,C_\Phi(x)\le C_\Phi(x^s)\le s\,C_\Phi(x)+(s-1)d_\Phi$.
\item The constants are optimal: for a Jordan frame $(c_i)$ and
$x_t=t\,c_1+\sum_{1<i<r}c_i+t^{-1}c_r$, $\lim_{t\to\infty}\bigl(C_\Phi(x_t^s)-sC_\Phi(x_t)\bigr)=(s-1)d_\Phi$;
in particular $C_\Phi(P(x_t^{1/2})x_t)-2C_\Phi(x_t)\to d_\Phi$.
\item $C_\Phi(P(x^{1/2})y)\le C_\Phi(x)+C_\Phi(y)$ for all $x,y\in\Omega$ if and only if
$d_\Phi=0$, if and only if $\Phi$ is a multiple of $\|\cdot\|_\infty$, i.e.\ $C_\Phi=\OLC$.
\end{enumerate}
In particular, for the Schatten contrasts the sharp defect is $\frac{\ln r}p$ (also for
$0<p<1$): it is $\frac{\ln2}p$ on every Lorentz cone, of any dimension, $\frac{\ln3}p$ on
the $27$-dimensional exceptional cone, and $\frac{\ln n}p$ on the cone of positive
$n\times n$ matrices.
\end{theorem}

\begin{proof}
We may assume $\Phi(\mathbf 1)=1$, so that $C_\Phi=\frac12\ln\kappa_\Phi$ and
$d_\Phi=-\ln\Phi(\epsilon_1)$; recall that $\Phi$ is monotone (G1) and
$\Phi(\epsilon_1)\|v\|_\infty\le\Phi(v)\le\Phi(\mathbf 1)\|v\|_\infty$ (G2).

(a) Let $w=P(x^{1/2})y$. Since $y\le\lambda_{\max}(y)e$ and $P(x^{1/2})$ preserves the
order, $w\le\lambda_{\max}(y)\,x$, so by the min--max principle
$\lambda^\downarrow(w)\le\lambda_{\max}(y)\lambda^\downarrow(x)$ coordinatewise, and
\[
\Phi(\lambda(w))\le\lambda_{\max}(y)\,\Phi(\lambda(x))\le\frac{\Phi(\lambda(x))\,\Phi(\lambda(y))}{\Phi(\epsilon_1)} .
\]
By (F1), $w^{-1}=P(x^{-1/2})y^{-1}$, so the same holds for the inverses; hence
\[
\kappa_\Phi(w)\le\Phi(\epsilon_1)^{-2}\kappa_\Phi(x)\kappa_\Phi(y).
\]
For $0<p<1$ the same argument applies to $v\mapsto(\frac1r\sum_iv_i^p)^{1/p}$, which is
monotone and satisfies $r^{-1/p}\|v\|_\infty\le(\frac1r\sum_iv_i^p)^{1/p}$.

(b) and (c) depend only on $\lambda=\lambda(x)$, since $\lambda(x^{\pm s})$ is $\lambda^{\pm s}$
up to order. (b) For $s=1$ there is nothing to prove; let $s>1$. By (G3) with $v=\lambda^s$, $w=\mathbf 1$, $\theta=1/s$,
$\Phi(\lambda)\le\Phi(\lambda^s)^{1/s}\Phi(\mathbf 1)^{1-1/s}=\Phi(\lambda^s)^{1/s}$, and likewise
for $\lambda^{-1}$; multiplying, $\kappa_\Phi(x)^s\le\kappa_\Phi(x^s)$. By (G1) and (G2),
$\Phi(\lambda^s)\le\|\lambda\|_\infty^{s-1}\Phi(\lambda)\le\Phi(\epsilon_1)^{1-s}\Phi(\lambda)^s$,
and likewise for $\lambda^{-1}$, so $\kappa_\Phi(x^s)\le\Phi(\epsilon_1)^{-2(s-1)}\kappa_\Phi(x)^s$.
Taking half logarithms gives (b). (c) Here $\lambda(x_t^{\pm s})=(t^{\pm s},1,\dots,1,t^{\mp s})$,
so $t^{-s}\lambda(x_t^{s})\to\epsilon_1$ and $t^{-s}\lambda(x_t^{-s})\to\epsilon_r$ as $t\to\infty$;
by continuity, homogeneity and permutation invariance of $\Phi$,
$t^{-2s}\kappa_\Phi(x_t^s)\to\Phi(\epsilon_1)^2$ for every $s\ge1$. Hence
\[
C_\Phi(x_t^s)-sC_\Phi(x_t)=\tfrac12\ln\frac{t^{-2s}\kappa_\Phi(x_t^s)}{\bigl(t^{-2}\kappa_\Phi(x_t)\bigr)^s}
\longrightarrow\tfrac12\ln\Phi(\epsilon_1)^{2-2s}=(s-1)d_\Phi ;
\]
note that $P(x_t^{1/2})x_t=x_t^2$. For $0<p<1$ the computation in (c) is the same, since it
uses only continuity, homogeneity and permutation invariance.

(d) If $d_\Phi>0$, (c) gives a violation; if $d_\Phi=0$, (a) is subadditivity. And
$\Phi(\epsilon_1)=\Phi(\mathbf 1)$ forces $\Phi(v)=\Phi(\epsilon_1)\|v\|_\infty$ for all $v$, by
the two-sided bound (G2).
\end{proof}

The key inequality $\lambda^\downarrow(P(x^{1/2})y)\le\lambda_{\max}(y)\lambda^\downarrow(x)$ and the
Schatten case of (a) are in Tao, Jeong and Gowda: \cite[Lemma~3.2]{TJG20} states it with
the roles of $x,y$ exchanged, which is equivalent by Lemma~\ref{cones:lem:polar}, and
\cite[Cor.~5.2(ii)]{TJG20} with $r=p$, $s=\infty$ gives the norm bound
$\|P(x^{1/2})y\|_p\le\|x\|_p\|y\|_\infty$. Their log-majorization
$\lambda^\downarrow(P(x^{1/2})y)\prec_{\log}\lambda^\downarrow(x)\ast\lambda^\downarrow(y)$ \cite[Thm.~3.4]{TJG20} gives the same
bound in (a), and no better one, since the extremal pair in (c) has
$\lambda(x_t^2)=\lambda(x_t)\ast\lambda(x_t)$. What is added here is the identification of the sharp defect with the rank, and the
rigidity. So the two basic constants of the theory are governed by the rank: the
non-commutativity threshold of Theorem~\ref{cones:thm:rank} ($\rk V\le3$) and the
subadditivity defect $\frac{\ln\rk V}p$ of the Schatten contrasts. On a Lorentz cone,
whose eigenvalues are $\rho e^{\pm\eta}$,
\begin{equation}\label{cones:eq:Cp-lorentz}
C_p(x)=\frac1p\ln\cosh\bigl(p\,\OLC(x)\bigr),
\end{equation}
a function of the rapidity alone.

\section{The Lorentz cone: contrast as velocity}\label{cones:sec:lorentz}

Let $V_n=\R\times\R^n$, $n\ge1$, with
$(x_0,\bar x)\circ(y_0,\bar y)=(x_0y_0+\bar x\cdot\bar y,\ x_0\bar y+y_0\bar x)$,
$e=(1,0)$, $\det x=x_0^2-|\bar x|^2$, and Minkowski form
$\langle x,y\rangle_{\rm M}=x_0y_0-\bar x\cdot\bar y$, so that $\det x=\langle x,x\rangle_{\rm M}$.
The cone $\Omega$ is the open forward light cone, the eigenvalues of $x$ are
$x_0\pm|\bar x|$, and $x^{-1}=(x_0,-\bar x)/\det x$ \cite[Ch.~II]{FK94}. Let
$H^n=\{x\in\Omega:\det x=1\}$ be the upper sheet of the hyperboloid with the
hyperbolic distance $d(x,y)=\operatorname{arcosh}\langle x,y\rangle_{\rm M}$. Every
$x\in\Omega$ is $x=\rho\,(\cosh\eta,\ \sinh\eta\,\bar u)$ with $\rho>0$, $\eta\ge0$,
$|\bar u|=1$, and then
\[
\OMC(x)=\frac{|\bar x|}{x_0}=\tanh\eta,\qquad \OLC(x)=\eta=d\bigl(e,x/\sqrt{\det x}\bigr):
\]
the Michelson contrast is the speed of the four-velocity $x/\sqrt{\det x}$ and the
logarithmic contrast its rapidity.

\begin{lemma}\label{cones:lem:boost}
For $a=(\cosh t,\sinh t\,\bar u)\in H^n$, $P(a)$ is the pure Lorentz boost (a symmetric
positive definite element of $\mathrm{SO}^+(1,n)$) of rapidity $2t$ along $\bar u$, with
$P(a)e=a^2$; in particular $P(a)$ preserves $\langle\cdot,\cdot\rangle_{\rm M}$ and $H^n$.
\end{lemma}

\begin{proof}
$P(a)$ is self-adjoint for the Euclidean inner product, and positive definite since
$P(a)=P(a^{1/2})^2$ by the fundamental formula \cite[Prop.~II.3.3]{FK94}, and $\det(P(a)y)=(\det a)^2\det y=\det y$ \cite[Prop.~III.4.2]{FK94}, so $P(a)$
preserves the quadratic form $\det$ and hence, by polarisation,
$\langle\cdot,\cdot\rangle_{\rm M}$. It maps $\Omega$ onto $\Omega$, so it lies in
the orthochronous Lorentz group; its determinant is positive since it is positive
definite. A symmetric positive definite element of $\mathrm{SO}^+(1,n)$ is a pure boost
(polar decomposition of the Lorentz group), and $P(a)e=a^2=(\cosh2t,\sinh2t\,\bar u)$
by (F1), which identifies the axis and the rapidity.
\end{proof}

\begin{theorem}\label{cones:thm:hyp}
\begin{enumerate}
\item For $x,y\in H^n$, the Thompson distance
$d_T(x,y)=\|\log P(x^{-1/2})y\|_\infty$, the projectivised Thompson distance
$\delta([x],[y])$ and the hyperbolic distance coincide:
$d_T(x,y)=\delta([x],[y])=d(x,y)$; Hilbert's projective metric is $2d(x,y)$.
\item For $x,y\in\Omega$ with rapidities $\eta_x=\OLC(x)$, $\eta_y=\OLC(y)$ and
angle $\theta$ between $\bar x$ and $\bar y$ (any $\theta$ if one of them is $0$),
\[
\cosh\OLC\bigl(P(x^{1/2})y\bigr)=\cosh\eta_x\cosh\eta_y+\sinh\eta_x\sinh\eta_y\cos\theta ,
\]
the hyperbolic law of cosines for the triangle with sides $\eta_x,\eta_y$ and angle
$\pi-\theta$ between them.
\item Consequently $\OLC(P(x^{1/2})y)\le\OLC(x)+\OLC(y)$, with equality if and only if
$\bar x$ and $\bar y$ are non-negatively proportional (this includes $\bar x=0$ or
$\bar y=0$), and $\OMC(P(x^{1/2})y)\le\OMC(x)\oplus\OMC(y)$, $u\oplus v=\frac{u+v}{1+uv}$,
with the same equality case. By contrast, $x$ and $y$ operator commute if and only
if $\bar x$ and $\bar y$ are linearly dependent.
\end{enumerate}
\end{theorem}

\begin{proof}
(a) By Lemma~\ref{cones:lem:boost}, $B=P(x^{-1/2})$ is a boost with $Bx=e$, so
$z=By\in H^n$ and $\langle e,z\rangle_{\rm M}=\langle Bx,By\rangle_{\rm M}=\langle x,y\rangle_{\rm M}$.
The eigenvalues of $z$ are $z_0\pm|\bar z|=e^{\pm s}$ with $\cosh s=z_0=\langle x,y\rangle_{\rm M}$,
so $\|\log z\|_\infty=s=d(x,y)$. Since $\det z=1$, the infimum over $\lambda>0$ of
$\|\log(\lambda z)\|_\infty$ is attained at $\lambda=1$, which gives $\delta$, and
$\log(\lambda_{\max}(z)/\lambda_{\min}(z))=2s$.

(b) Both sides are invariant under $x\mapsto\lambda x$, $y\mapsto\mu y$, so let
$x,y\in H^n$. Then $B=P(x^{1/2})$ is a boost, $z=By\in H^n$, and
$\OLC(z)=\operatorname{arcosh}z_0$ by the formula for the eigenvalues. Since $B$ preserves
$\langle\cdot,\cdot\rangle_{\rm M}$ and $B^{-1}=P(x^{-1/2})$ maps $e$ to $x^{-1}$,
$z_0=\langle e,By\rangle_{\rm M}=\langle B^{-1}e,y\rangle_{\rm M}=\langle x^{-1},y\rangle_{\rm M}
=x_0y_0+\bar x\cdot\bar y$, which is the right-hand side.

(c) The right-hand side of (b) is at most $\cosh(\eta_x+\eta_y)$, with equality iff
$\cos\theta=1$ or $\eta_x\eta_y=0$; and $\tanh$ is increasing, with
$\tanh(\eta_x+\eta_y)=\tanh\eta_x\oplus\tanh\eta_y$. For the last statement, $x$ and
$y$ have a common Jordan frame iff they lie in a common subalgebra
$\R e\oplus\R(0,\bar u)$, i.e.\ iff $\bar x,\bar y$ are linearly dependent (F5).
\end{proof}

So on the Lorentz cone subadditivity of $\OLC$ is the triangle inequality of
hyperbolic space, $P(x^{1/2})y$ is $y$ seen after the boost that takes the observer
$e$ to $x$, and (b) is the composition law of rapidities; equality requires
collinear boosts in the same direction, which is strictly more than commutativity
(opposite boosts commute, and then $\OLC(P(x^{1/2})y)=|\eta_x-\eta_y|$). The
identification of rapidity with hyperbolic distance on boosts and on positive
$2\times2$ matrices is classical \cite{RS04,KL13,Kim11}; the equality case of the hyperbolic triangle
inequality is classical too, and Theorem~\ref{cones:thm:hyp} only records all this in terms
of the contrast.

\begin{remark}[Spinor and vector representations; radar]\label{cones:rem:spinor}
$V_3\cong\mathrm{Herm}_2(\C)$ via $x\mapsto x_0I+\bar x\cdot\vec\sigma$, with $\det$ the
matrix determinant and $P(x^{1/2})y=X^{1/2}YX^{1/2}$. For $g\in\mathrm{SL}(2,\C)$ with
singular values $e^{\pm\eta/2}$, $gg^*\in H^3$ has rapidity $\eta$, so
\[
\OLC(gg^*)=2\,\OLC(g)=d(e,gg^*),\qquad \OMC(gg^*)=\tanh\eta,\qquad\OMC(g)=\tanh\tfrac\eta2 :
\]
the contrast of the element $gg^*=\Lambda(g)e$ of $V_3$ (equal, by
Lemma~\ref{cones:lem:kak}, to that of the vector-representation matrix $\Lambda(g)$) is the
radius of the point in the Klein model, and the contrast of the spinor matrix $g$ is
the radius in the Poincar\'e ball model. In particular, for $2\times2$ matrices the projectivised
Thompson geometry of the cone of positive matrices coincides with the rapidity
(hyperbolic) geometry; see also \cite{PK15} for the hyperbolic geometry of the Bloch ball. For $n\times n$ matrices, $n\ge3$, the Thompson metric restricted to the determinant-one slice of the
positive cone is a non-Riemannian Finsler metric on $\mathrm{SL}(n,\C)/\mathrm{SU}(n)$,
and not hyperbolic.
For a positive diagonal matrix $\operatorname{diag}(t,T)$,
$\OMC=\frac{T-t}{T+t}$ is the radar velocity of Bondi's $k$-calculus, with
$T/t=k^2$ and $\OLC=\ln k$ \cite{Bondi64}. Writing $\OMC(t,S)=\OMC(\operatorname{diag}(t,S))$,
the identity $\OMC(t,S)=\OMC(t,T)\oplus\OMC(T,S)$ for $t<T<S$ is the equality case of
Theorem~\ref{cones:thm:hyp}(c); kinematically it is the multiplicativity
$k_{AC}=k_{AB}k_{BC}$ of Bondi factors.
\end{remark}

\begin{remark}[Wigner rotation and the angular defect of the contrast triangle]\label{cones:rem:wigner}
For $x,y\in H^n$, Lemma~\ref{cones:lem:polar} reads $P(x^{1/2})P(y^{1/2})=P(w^{1/2})k$ with
$w=P(x^{1/2})y$: the composition of the boost taking $e$ to $y$ followed by the boost
taking $e$ to $x$ equals the rotation $k\in\mathrm{Aut}(V_n)\cap\mathrm{SO}(1,n)$, which
fixes $e$ and acts on $\bar x$-space by an element of $\mathrm{SO}(n)$, followed by the
boost taking $e$ to $w$. Thus $k$ is the Thomas--Wigner rotation. Its angle equals the
area $\pi-(\alpha+\beta+\gamma)$ of the hyperbolic triangle with vertices $e$, $x$, $w$,
with the orientation convention of \cite[Eqs.~(93), (113b)]{RS04} (rotation $=-$area;
the result goes back to \cite{Aravind97}); the sides of this triangle are $\OLC(x)$,
$\OLC(y)$ and $\OLC(P(x^{1/2})y)$, by Theorem~\ref{cones:thm:hyp}(a) and the invariance of $d$
under the boost $P(x^{1/2})$. So the three contrasts in the subadditivity inequality are
the sides of a triangle whose area is the Wigner angle. The triangle degenerates and
$k=\mathrm{id}$ if and only if $\bar x$ and $\bar y$ are linearly dependent, i.e.\ iff $x$ and
$y$ operator commute (Theorem~\ref{cones:thm:hyp}(c)); the equality case is the subcase of
non-negative proportionality. Consequently, equality is not equivalent to vanishing gyration: on Lorentz cones the gyration
(Wigner rotation) vanishes iff $x,y$ operator commute, which is strictly weaker than
equality. (The component manuscript reports numerical checks on $V_3$ for random pairs with rapidities up to $3$, with error $\le10^{-10}$; that full check was not repeated for this consolidation.) Reciprocity,
$\OLC(P(x^{1/2})y)=\OLC(P(y^{1/2})x)$, is the familiar fact that $u\oplus v$ and
$v\oplus u$ have the same speed.
\end{remark}

\section{The Jordan product}\label{cones:sec:jordan}

For matrices $X\circ Y=\frac12(XY+YX)$, and estimating the spectrum of the Jordan
product of two positive matrices in terms of their spectra is a research problem of
Taussky-Todd. Strang \cite{Strang62} found the sharp bounds for the extreme eigenvalues
in every dimension $\ge2$, and in Hilbert space; see also Gustafson \cite{Gus68},
Nicholson \cite{Nich79}, de Bruijn \cite{dB80}, and Alikakos and Bates \cite{AB84} with
the erratum \cite{AB85}, which states Strang's lower bound in explicit form. These bounds
concern $\lambda_{\max}$ and $\lambda_{\min}$ separately. The contrast is their ratio, and
the point of this section is that whether the two bounds can be attained \emph{at the
same time} depends on the algebra: it is again a question of room, as in
Theorem~\ref{cones:thm:rank}.

For $x\in\Omega$ write $\mu(x)=\frac12(\lambda_{\max}(x)+\lambda_{\min}(x))$, so that the
spectrum of $x$ lies in $\mu(x)[1-u,1+u]$ with $u=\OMC(x)$ and both endpoints occur.
For $0\le u,v<1$ put
\[
\ell(u,v)=\begin{cases}(1-u)(1-v),& u+v\le1,\\[2pt] \tfrac12(1-u^2-v^2),& u+v\ge1,\end{cases}
\qquad u\oplus v=\frac{u+v}{1+uv}.
\]
(The two expressions agree on $u+v=1$, and for $u+v\ne1$ one has
$\frac12(1-u^2-v^2)<(1-u)(1-v)$, since the difference is $\frac12(u+v-1)^2$.)

\begin{lemma}[Strang's bounds in Jordan form]\label{cones:lem:strang}
Let $A$ be a unital JB-algebra \cite{HOS84}, $A_{++}$ the set of invertible elements of its
cone of squares, and define $\lambda_{\max}$, $\lambda_{\min}$, $\mu$, $\OMC$ on $A_{++}$ as
above (see also Section~\ref{cones:sec:jbw}). Let $x,y\in A_{++}$, $u=\OMC(x)$, $v=\OMC(y)$. Then
\[
\lambda_{\max}(x\circ y)\le\lambda_{\max}(x)\lambda_{\max}(y),\qquad
\lambda_{\min}(x\circ y)\ge\mu(x)\mu(y)\,\ell(u,v).
\]
In particular $x\circ y\in A_{++}$ whenever $u^2+v^2<1$.
\end{lemma}

\begin{proof}
The upper bound is $\|x\circ y\|\le\|x\|\|y\|$. Let $C$ be the closed subalgebra
generated by $e$, $x$, $y$. By the Shirshov--Cohn theorem \cite[Thm.~7.2.5]{HOS84}, $C$ is
a JC-algebra: there is a Jordan isomorphism $\pi$ of $C$ onto a norm-closed
Jordan algebra of self-adjoint operators on a Hilbert space, and compressing to the
range of the projection $\pi(e)$ we may assume $\pi(e)=1$. The spectrum of an element of $C$ is computed
in the associative subalgebra generated by it and $e$, so it is the same in $C$ and in
$A$, and $\pi$ preserves it. Hence it suffices to bound $\lambda_{\min}$ of
$\frac12(XY+YX)$ for positive operators $X,Y$ with spectra in $\mu(x)[1-u,1+u]$ and
$\mu(y)[1-v,1+v]$ (if $u=0$ or $v=0$ the claim is trivial, so let $u,v>0$). By Strang's
theorem \cite[Thm.~I]{Strang62}, which is stated for operators on Hilbert space (and in any
case reduces to finite dimension: for a unit vector $w$ and the orthogonal projection $Q$
onto $F=\operatorname{span}\{w,Xw,Yw\}$, the compressions $QX|_F$, $QY|_F$ have spectra in the
same intervals, and $\langle(XY+YX)w,w\rangle=2\operatorname{Re}\langle Yw,Xw\rangle$ does not
change when $X,Y$ are replaced by them), in the explicit form of
\cite{AB85}, $\lambda_{\min}(XY+YX)\ge2mm'$ if $\frac{M'}{m'}\le\frac{M/m+3}{M/m-1}$ and
$\lambda_{\min}(XY+YX)\ge\frac{16MmM'm'-(M-m)^2(M'-m')^2}{4(M+m)(M'+m')}$ otherwise, where
$\sigma(X)\subset[m,M]$, $\sigma(Y)\subset[m',M']$. With $M,m=\mu(x)(1\pm u)$ and
$M',m'=\mu(y)(1\pm v)$ the condition reads $u+v\le1$, $2mm'=2\mu(x)\mu(y)(1-u)(1-v)$, and the
second bound is $\mu(x)\mu(y)(1-u^2-v^2)$. Halving gives the claim. Positivity for
$u^2+v^2<1$ follows, since then both expressions for $\ell$ are positive; this criterion
is also \cite[Cor.~4.4]{Gus68}.
\end{proof}

\begin{theorem}\label{cones:thm:jordan}
Let $V$ be a non-associative Euclidean Jordan algebra, and for $0\le u,v<1$ with
$u^2+v^2<1$ put
\[
S_V(u,v)=\sup\{\OMC(x\circ y):\ x,y\in\Omega,\ \OMC(x)=u,\ \OMC(y)=v\}.
\]
\begin{enumerate}
\item If $V$ is a spin factor $V_n$, $n\ge2$, then
$S_V(u,v)=u\oplus v$ if $u^2+uv+v^2\le1$, and $S_V(u,v)=(2-u^2-v^2)^{-1/2}$ otherwise.
\item If $V$ is not a spin factor, then $S_V(u,v)=u\oplus v$ if $u+v\le1$, and otherwise
$S_V(u,v)=\frac{\kappa_*-1}{\kappa_*+1}$ with $\kappa_*=\frac{2(1+u)(1+v)}{1-u^2-v^2}$.
\end{enumerate}
In both cases the supremum is attained. Consequently, for $u^2+v^2<1$,
$\OLC(x\circ y)\le\OLC(x)+\OLC(y)$ for all $x,y\in\Omega$ with $\OMC(x)=u$, $\OMC(y)=v$
if and only if $u^2+uv+v^2\le1$ in
case (a), and if and only if $u+v\le1$ in case (b).
\end{theorem}

Since $\tanh(\OLC(x)+\OLC(y))=u\oplus v$, the last statement follows from (a) and (b) once
we check that the second expressions exceed $u\oplus v$ outside the respective regions;
this is done in the proof. If $V$ is associative, $x\circ y$ is the product of commuting
elements and $\OLC(x\circ y)\le\OLC(x)+\OLC(y)$ always (Theorem~\ref{cones:thm:sub}). The
ellipse $u^2+uv+v^2\le1$ strictly contains the triangle $u+v\le1$, and both lie inside the
disc $u^2+v^2<1$ where $x\circ y$ is positive. By Corollary~IV.1.5 of \cite{FK94} the simple
algebras of rank two are spin factors (and conversely $V_n$, $n\ge2$, is simple, as is
well known), so case (b) consists of the simple
algebras of rank $\ge3$ and of the non-simple algebras with a non-associative summand.

\begin{proof}
(a) Write $x=\alpha(1,u\bar a)$, $y=\beta(1,v\bar b)$ with $\alpha,\beta>0$ and unit
vectors $\bar a,\bar b$ with $\bar a\cdot\bar b=\vartheta\in[-1,1]$; every $\vartheta$ occurs since
$n\ge2$. From the product in $V_n$,
\[
x\circ y=\alpha\beta\,\bigl(1+uv\,\vartheta,\ u\bar a+v\bar b\bigr),\qquad
\OMC(x\circ y)=\frac{|u\bar a+v\bar b|}{1+uv\vartheta }=\frac{\sqrt{u^2+v^2+2uv\vartheta }}{1+uv\vartheta }.
\]
If $uv=0$ this is $\max(u,v)=u\oplus v$, and $u^2+uv+v^2<1$. Let $uv>0$ and put
$s=1+uv\vartheta\in[1-uv,1+uv]$ and $Q=2-u^2-v^2$. Then $\OMC(x\circ y)^2=\gamma(s)=(2s-Q)/s^2$,
and $\gamma'(s)=2(Q-s)/s^3$, so $\gamma$ increases on $(0,Q]$ and decreases on
$[Q,\infty)$. Since $u^2+v^2<1$, $Q-(1-uv)=1-u^2-v^2+uv>0$, so
$\sup\gamma=\gamma(\min\{Q,1+uv\})$. If $Q\ge1+uv$, i.e.\ $u^2+uv+v^2\le1$, this is
$\gamma(1+uv)=(u\oplus v)^2$. Otherwise it is $\gamma(Q)=1/Q$, attained at
$\vartheta =(Q-1)/(uv)\in(-1,1)$, and $1/Q=\gamma(Q)>\gamma(1+uv)=(u\oplus v)^2$.

(b) By Lemma~\ref{cones:lem:strang},
$\kappa(x\circ y)\le\frac{\lambda_{\max}(x)\lambda_{\max}(y)}{\mu(x)\mu(y)\ell(u,v)}=\frac{(1+u)(1+v)}{\ell(u,v)}$.
For $u+v\le1$ this is $\kappa(x)\kappa(y)$, i.e.\ $\OMC(x\circ y)\le u\oplus v$, attained by
operator commuting $x,y$ with primitive $d\le E_{\max}(x),E_{\max}(y)$ and
$d'\le E_{\min}(x),E_{\min}(y)$: for operator commuting elements $x\circ y=P(x^{1/2})y$, and
Theorem~\ref{cones:thm:eq} applies. Let $u+v>1$. The upper bound is
$\kappa_*$ as stated, and it exceeds $\kappa(x)\kappa(y)$ because
$2(1-u)(1-v)>1-u^2-v^2$. For attainment, $V$ contains orthogonal primitive idempotents
$c_2,c_3$ with $V_{23}\ne0$ (F6) and a nonzero idempotent $c_1\perp c_2+c_3$: take a frame of a
simple component of rank $\ge3$, or $c_2,c_3$ in a non-associative simple component (of rank
$\ge2$) and $c_1$ the unit of another component. The Peirce algebra $U=V(c_2+c_3,1)$ has rank
two and contains $V_{23}\ne0$, so it is simple (a non-simple algebra of rank two is
$\R\oplus\R$, where $V_{23}=0$), hence a spin factor of dimension $\ge3$, and by the
proof of (a) there are $x',y'\in U$ with spectra $\{1\pm u\}$, $\{1\pm v\}$ and
\[
\lambda_{\min}(x'\circ y')=\min_\vartheta\bigl(1+uv\vartheta -\sqrt{u^2+v^2+2uv\vartheta }\bigr)=\tfrac12(1-u^2-v^2),
\]
the minimum being attained where $u^2+v^2+2uv\vartheta =1$, i.e.\ at
$\vartheta =\frac{1-u^2-v^2}{2uv}\in[0,1]$, which is admissible because $(u+v)^2\ge1$. Put
$r_0=e-c_1-c_2-c_3$ and
\[
x=(1+u)c_1+x'+r_0,\qquad y=(1+v)c_1+y'+r_0 .
\]
Then $\OMC(x)=u$, $\OMC(y)=v$, and by (F4) $x\circ y=(1+u)(1+v)c_1+x'\circ y'+r_0$, whose
spectrum is the union of $\{(1+u)(1+v)\}$, the spectrum of $x'\circ y'$ in $U$, and $\{1\}$
if $r_0\ne0$. As $x'\circ y'\le(1+u)(1+v)(c_2+c_3)$ by Lemma~\ref{cones:lem:strang} and
$\frac12(1-u^2-v^2)<1$, we get $\lambda_{\max}=(1+u)(1+v)$ and
$\lambda_{\min}=\frac12(1-u^2-v^2)$, so $\kappa_*$ is attained.
\end{proof}

\begin{corollary}\label{cones:cor:jordan-cases}
On $\mathrm{Herm}_2(\mathbb F)$, $\mathbb F=\R,\C,\mathbb H,\mathbb O$ (spin factors), the
contrast of the Jordan product is subadditive exactly on the ellipse
$u^2+uv+v^2\le1$; on $\mathrm{Herm}_n(\mathbb F)$, $\mathbb F=\R,\C,\mathbb H$, $n\ge3$, on
the Albert algebra $\mathrm{Herm}_3(\mathbb O)$, and on $\R\oplus V_n$, $V_m\oplus V_n$
($n\ge2$) exactly on the triangle $u+v\le1$. In particular, for positive complex
$n\times n$ matrices the threshold changes between $n=2$ and $n=3$.
\end{corollary}

The mechanism is the same as in Theorem~\ref{cones:thm:rank}. Strang's two bounds are attained by
different configurations: the upper one by an idempotent under both $E_{\max}$, the lower
one, for $u+v>1$, by a non-commuting pair in a rank-two Peirce algebra. Fitting both needs a
third orthogonal piece, i.e.\ room beyond a spin factor; in a spin factor the
two extremes compete, and the region of subadditivity is larger. The quantity in the proof
of (a) is the speed of $x\circ y$: for four-velocities with speeds $u,v$ and angle $\theta$
between them, $\OMC(x\circ y)=|u\bar a+v\bar b|/(1+uv\cos\theta)$ is the Galilean sum of the
velocities divided by the relativistic denominator.

The positivity criterion $u^2+v^2<1$ is \cite[Cor.~4.4]{Gus68} and \cite{Strang62}, and the
bounds of Lemma~\ref{cones:lem:strang} for matrices are \cite{Strang62,AB85}; de Bruijn
\cite[\S4]{dB80} uses a parametrisation close to the one in the proof of (a). The
joint attainability of the two extremes is not discussed in
\cite{Strang62,Gus68,Nich79,dB80,AB84,AB85}.

\begin{remark}[Multiplicative commutators of spinor boosts]\label{cones:rem:groupcomm}
For positive matrices the additive commutator is controlled by Strang's
$\|XY-YX\|\le\frac12(M-m)(N-n)$ \cite[Thm.~II]{Strang62}, attained in dimension two; in
terms of contrasts, $\|XY-YX\|\le2\,\OMC(X)\OMC(Y)\mu(X)\mu(Y)$. The multiplicative
commutator $K=XYX^{-1}Y^{-1}$ of two positive elements of $\mathrm{SL}(2,\C)$, i.e.\ of two
spinor boosts $X=\exp(\alpha\,\bar a\cdot\vec\sigma)$, $Y=\exp(\beta\,\bar b\cdot\vec\sigma)$
with $\OLC(X)=\alpha$, $\OLC(Y)=\beta$ (boost rapidities $2\alpha$, $2\beta$, Remark~\ref{cones:rem:spinor})
and angle $\theta$ between the axes, satisfies
\begin{gather*}
\tr K=2-4\sinh^2\!\alpha\,\sinh^2\!\beta\,\sin^2\theta,\\
\sinh\OLC(K)=\tfrac12\sinh2\alpha\,\sinh2\beta\,|\sin\theta|\,\bigl|\tanh\alpha+e^{i\theta}\tanh\beta\bigr| ;
\end{gather*}
both follow from the expansion of $K$ in the Pauli basis. Indeed, with
$(s_0+\bar s\cdot\vec\sigma)(t_0+\bar t\cdot\vec\sigma)=(s_0t_0+\bar s\cdot\bar t)+(s_0\bar t+t_0\bar s+i\,\bar s\times\bar t)\cdot\vec\sigma$
for $s_0,t_0\in\C$, $\bar s,\bar t\in\C^3$ (bilinear products), and
$X^{\pm1}=\cosh\alpha\pm\sinh\alpha\,\bar a\cdot\vec\sigma$,
$Y^{\pm1}=\cosh\beta\pm\sinh\beta\,\bar b\cdot\vec\sigma$, $\bar a\cdot\bar b=\cos\theta$, one finds
$XY=p_0+(\bar r+i\bar n)\cdot\vec\sigma$ and $X^{-1}Y^{-1}=p_0+(-\bar r+i\bar n)\cdot\vec\sigma$,
where $p_0=\cosh\alpha\cosh\beta+\sinh\alpha\sinh\beta\cos\theta$,
$\bar r=\sinh\alpha\cosh\beta\,\bar a+\cosh\alpha\sinh\beta\,\bar b$ and
$\bar n=\sinh\alpha\sinh\beta\,\bar a\times\bar b$. Since $\bar r\perp\bar n$ and $\det(XY)=1$,
$p_0^2=1+|\bar r|^2-|\bar n|^2$, and multiplying out gives
$K=(1-2|\bar n|^2)+(2ip_0\bar n-2\,\bar r\times\bar n)\cdot\vec\sigma$. Hence
$\tr K=2-4|\bar n|^2$ with $|\bar n|=\sinh\alpha\sinh\beta|\sin\theta|$; and since the singular
values of $K\in\mathrm{SL}(2,\C)$ are $e^{\pm\OLC(K)}$,
$4\sinh^2\OLC(K)=\tr(K^*K)-2=2(1-2|\bar n|^2)^2+8|\bar n|^2(p_0^2+|\bar r|^2)-2=16|\bar n|^2|\bar r|^2$,
where $|\bar r|=\cosh\alpha\cosh\beta\,|\tanh\alpha+e^{i\theta}\tanh\beta|$. The trace formula is the classical trace identity for commutators
in $\mathrm{SL}(2)$ (Fricke; see e.g.\ \cite{Bea83}), specialised to boosts. So for $\alpha,\beta>0$ and
$\theta\in(0,\pi)$ the element $K$ is elliptic, i.e.\ conjugate to a non-trivial rotation,
iff $\sinh\alpha\sinh\beta|\sin\theta|<1$, and $K$ has its spectrum on the unit circle for
all angles iff $\sinh\alpha\sinh\beta\le1$; the strict inequality is exactly $u^2+v^2<1$ with $u=\OMC(X)=\tanh\alpha$, $v=\OMC(Y)=\tanh\beta$,
the positivity criterion for $X\circ Y$. This is not a coincidence, and the sufficiency
holds for positive invertible $X,Y$ in every unital $C^*$-algebra: with
$G=X^{-1}Y^{-1}+Y^{-1}X^{-1}$ and $K^*=Y^{-1}X^{-1}YX$ one has
$GK=X^{-1}Y^{-1}XYX^{-1}Y^{-1}+X^{-1}Y^{-1}$ and hence $K^*GK=X^{-1}Y^{-1}+Y^{-1}X^{-1}=G$. If
$u^2+v^2<1$, then $G=2X^{-1}\circ Y^{-1}$ is positive and invertible by Lemma~\ref{cones:lem:strang},
applied to $X^{-1},Y^{-1}$, which have contrasts $u,v$; so $G^{1/2}KG^{-1/2}$ is an invertible
isometry, hence unitary, and the spectrum of $K$ lies on the unit circle. For positive $n\times n$
matrices with given contrasts the suprema of $\OLC(K)$ and of the spectral radius of $K$ are
non-decreasing in $n$ ($K$ does not change under scaling of $X$ and $Y$, so one may assume
that $1$ lies in both spectral intervals and pass to $X\oplus1$, $Y\oplus1$); numerically
they are the same for $n=3,4$ as for $n=2$, so unlike the Jordan product these quantities
seem to show no dimension threshold, but we have no proof.
\end{remark}

\section{Boosts in representations of \texorpdfstring{$\mathrm{SL}(2,\C)$}{SL(2,C)}}\label{cones:sec:reps}

The non-commuting equality case for $4\times4$ positive matrices Theorem~\ref{core:dimension}
raises the question whether it has a kinematic meaning, for instance in the Dirac
representation. It has none, in the following sense.

Let $\rho$ be a finite-dimensional continuous representation of $\mathrm{SL}(2,\C)$ on a
complex vector space $W$ (a real representation, such as the vector representation
on $\R^4$, is complexified; this does not change $\OLC$). Continuous representations
are smooth, and the real Lie algebra $\mathfrak{sl}(2,\C)$ is semisimple, so $\rho$ is
completely reducible, with irreducible components of type $(j_1,j_2)$,
$\mathcal D^{j_1}\otimes\overline{\mathcal D^{j_2}}$, $j_1,j_2\in\frac12\mathbb N_0$
\cite[Ch.~II]{Knapp86}; here $\overline{\mathcal D^{j}}(g)$ is the entrywise complex
conjugate of $\mathcal D^j(g)$ in an orthonormal basis (unitarily equivalent to the
physicists' $g\mapsto\mathcal D^j(g^*)^{-1}$). We call an inner product on $W$
\emph{admissible} if
\[
\rho(g^*)=\rho(g)^*\qquad(g\in\mathrm{SL}(2,\C)).
\]
Admissible inner products exist: on $\mathcal D^{j_1}\otimes\overline{\mathcal D^{j_2}}$
the tensor product of the $\mathrm{SU}(2)$-invariant inner products is admissible, and
an orthogonal sum of admissible inner products is admissible. On each irreducible
block an admissible inner product is unique up to a positive factor (Schur), and
for an admissible inner product invariant subspaces have invariant orthogonal
complements. (Mere $\mathrm{SU}(2)$-invariance is not enough: on the vector
representation $\langle x,y\rangle=x_0y_0+4\,\bar x\cdot\bar y$ is rotation-invariant,
but boosts are not self-adjoint for it.) For admissible inner products,
$\rho(b)=\rho(b^{1/2})^*\rho(b^{1/2})$ is positive definite for every positive
$b\in\mathrm{SL}(2,\C)$, and $\rho(b)^{1/2}=\rho(b^{1/2})$.

Let $s(\rho)=\max_k\bigl(j_1^{(k)}+j_2^{(k)}\bigr)$ over the irreducible components, and,
for a unit vector $\bar n\in\R^3$ and $\eta\ge0$, let
$b_{\bar n}(\eta)=\exp(\frac\eta2\,\bar n\cdot\vec\sigma)$, the pure boost of rapidity
$\eta$ along $\bar n$ (Remark~\ref{cones:rem:spinor}).

\begin{lemma}\label{cones:lem:kak}
For an admissible inner product, $\OLC(\rho(g))=2s(\rho)\,\OLC(g)$ for every
$g\in\mathrm{SL}(2,\C)$.
\end{lemma}

\begin{proof}
Write $g=k_1ak_2$ with $k_1,k_2\in\mathrm{SU}(2)$ and $a=\operatorname{diag}(e^{t/2},e^{-t/2})$,
$t=2\,\OLC(g)\ge0$ (the $KAK$, i.e.\ singular value, decomposition). The operators
$\rho(k_i)$ are unitary, since $\rho(k^*)=\rho(k^{-1})$, so $\OLC(\rho(g))=\OLC(\rho(a))$.
On a component of type $(j_1,j_2)$, $\rho(a)=\exp(tK)$ with
$K=J_3^{(1)}\otimes1+1\otimes\overline{J_3^{(2)}}$, whose eigenvalues are $m_1+m_2$,
$|m_i|\le j_i$, so the spectrum of $\rho(a)$ is $\{e^{tm}\}$ with $m$ ranging over a set
symmetric about $0$ whose largest element is $s(\rho)$. Hence
$\OLC(\rho(a))=\frac12\ln(e^{ts}/e^{-ts})=ts=2s\,\OLC(g)$.
\end{proof}

\begin{theorem}\label{cones:thm:reps}
Let $\rho$ be as above with $s(\rho)>0$ and an admissible inner product, let
$\eta_1,\eta_2>0$, and put $X=\rho(b_{\bar n_1}(\eta_1))$, $Y=\rho(b_{\bar n_2}(\eta_2))$.
Then $X,Y$ are positive definite, $\OLC(X)=s(\rho)\eta_1$, $\OLC(Y)=s(\rho)\eta_2$, and
the following are equivalent:
\begin{enumerate}
\item $\OLC(XY)=\OLC(X)+\OLC(Y)$;
\item $\OLC(X^{1/2}YX^{1/2})=\OLC(X)+\OLC(Y)$;
\item $\bar n_1=\bar n_2$.
\end{enumerate}
\end{theorem}

\begin{proof}
Positivity was noted above, and $\OLC(X)=2s\,\OLC(b_{\bar n_1}(\eta_1))=s\eta_1$ by
Lemma~\ref{cones:lem:kak}. Since $XY=\rho(b_1b_2)$ and $X^{1/2}YX^{1/2}=\rho(b_1^{1/2}b_2b_1^{1/2})$,
with $b_i=b_{\bar n_i}(\eta_i)$, Lemma~\ref{cones:lem:kak} reduces (a) and (b) to the same
statements for the $2\times2$ matrices $b_1,b_2$, i.e.\ to $\rho=\mathrm{id}$. For positive
$2\times2$ matrices the equality cases of $XY$ and of $X^{1/2}YX^{1/2}$ coincide (Remark
after Theorem~\ref{cones:thm:eq}), and in $\mathrm{Herm}_2(\C)\cong V_3$ the latter is
Theorem~\ref{cones:thm:hyp}(c): equality holds iff the space parts of $b_1,b_2$, which are
$\sinh(\eta_i/2)\,\bar n_i$, are non-negatively proportional, i.e.\ iff
$\bar n_1=\bar n_2$.
\end{proof}

\begin{corollary}\label{cones:cor:dirac}
In no non-trivial finite-dimensional representation of $\mathrm{SL}(2,\C)$ with an
admissible inner product does a pair of pure boosts give a non-commuting equality
case: equality forces equal axes, hence commutation. In particular no
non-commuting equality pair of the kind of Theorem~\ref{core:dimension} arises from boosts of
a Dirac spinor $(\frac12,0)\oplus(0,\frac12)$ (for the specific example of
\cite{NN26} this is already clear from its spectrum $\{4,3,1,1\}$, which is not of the
form $\{a,a,a^{-1},a^{-1}\}$ up to scaling).
\end{corollary}

What non-commuting equality pairs do realise is orthogonal blocks, i.e.\ direct sums
(not tensor products, cf.\ Proposition~\ref{core:tensor-equality}). On
two spinors, i.e.\ for $\mathrm{SL}(2,\C)\times\mathrm{SL}(2,\C)$ acting on
$\C^2\oplus\C^2$, let
\[
X=b_z(2\ln2)\oplus b_z(\ln2),\qquad Y=b_z(2\ln2)\oplus b_x(\ln2).
\]
Then $\lambda_{\max}=2$, $\lambda_{\min}=\frac12$ for both, the first blocks give
$\|XY\|=4$ and $\|(XY)^{-1}\|=4$ (the second blocks have norms at most $2$), so
$\OLC(XY)=\OLC(X)+\OLC(Y)$, and likewise for $X^{1/2}YX^{1/2}$, while $XY\ne YX$: the
first spinor is boosted collinearly by both, the second along different axes by a
smaller rapidity, which does not reach the extremes. The example of Theorem~\ref{core:dimension} has the same structure. This is the
structure behind Theorem~\ref{cones:thm:rank}: the extremes are fixed by idempotents
$c\le E_{\max}$, $c'\le E_{\min}$, and the non-commutation sits in an orthogonal
rank-two piece. For non-positive elements equality without commutation is easy even
for one spinor: $g_1=BU$, $g_2=U^*B$ with $B>0$, $U$ unitary and $[U,B^2]\ne0$ give
$g_1g_2=B^2$ and $\OLC(g_1g_2)=\OLC(g_1)+\OLC(g_2)$; these lie outside the positive
setting of Theorem~\ref{core:dimension}.

\section{Infinite dimension: JB- and JBW-algebras}\label{cones:sec:jbw}

Let $A$ be a unital JB-algebra \cite{HOS84}, $A_+$ its cone of squares, and $A_{++}$
the set of invertible elements of $A_+$. For $x\in A_{++}$ put
$\lambda_{\max}(x)=\max\sigma(x)$, $\lambda_{\min}(x)=\min\sigma(x)>0$ and
$\OLC(x)=\frac12\ln(\lambda_{\max}(x)/\lambda_{\min}(x))$.

\begin{proposition}\label{cones:prop:jb}
For $x,y\in A_{++}$, $\OLC(P(x^{1/2})y)\le\OLC(x)+\OLC(y)$.
\end{proposition}

\begin{proof}
The proof of Theorem~\ref{cones:thm:sub} uses only that $P(a)$ is a positive map
\cite[Prop.~3.3.6]{HOS84} and that $\lambda_{\min}(z)e\le z\le\lambda_{\max}(z)e$ for
$z\in A_{++}$.
\end{proof}

Let now $M$ be a JBW-algebra and define its \emph{rank} $\rk M\in\{1,2,\dots,\infty\}$
as the supremum of the cardinalities of families of pairwise orthogonal nonzero
projections; for Euclidean Jordan algebras this is the usual rank. Elements $x,y$
\emph{operator commute} if $[T_x,T_y]=0$, where $T_x=x\circ\cdot$.

\begin{theorem}\label{cones:thm:jbw}
For a JBW-algebra $M$ the following are equivalent:
\begin{enumerate}
\item for all $x,y\in M_{++}$, $\OLC(P(x^{1/2})y)=\OLC(x)+\OLC(y)$ implies that $x$ and $y$
operator commute;
\item $\rk M\le3$, or $M$ is associative.
\end{enumerate}
\end{theorem}

\begin{proof}
(b)$\Rightarrow$(a). If $M$ is associative there is nothing to prove. Let $\rk M\le3$.
Then the centre of $M$ has at most three minimal projections, and every factor summand
has rank at most $3$; factors of type $\mathrm{II}$, $\mathrm{III}$ and $\mathrm I_\infty$
contain infinite orthogonal families of projections, so every summand is of type
$\mathrm I_n$, $n\le3$. The JBW-factors of types $\mathrm I_1$, $\mathrm I_2$, $\mathrm I_3$
are $\R$, the spin factors $V_\cH=\R\oplus\cH$ over real Hilbert spaces $\cH$ of arbitrary
dimension $\ge2$, and $\mathrm{Herm}_3(\mathbb F)$, $\mathbb F=\R,\C,\mathbb H,\mathbb O$
\cite[6.1.8, 6.3.1, Thm.~6.4.1]{HOS84}, \cite{Stormer66,Shultz79}; that a JBW-algebra of
finite rank is a finite direct sum of algebras of type $\mathrm I_n$, $n$ bounded by the
rank, follows from \cite[Lemma~5.3.4, Thm.~5.3.5]{HOS84}. So $M$ is a direct sum, of total rank
at most $3$, of such algebras. Every summand carries a faithful associative trace
($\tr(\alpha+v)=2\alpha$ on $V_\cH$), every element has a spectral decomposition over a
Jordan frame of at most three minimal projections, the Peirce spaces $V_{ii}$ of such a
frame are $\R c_i$, and the Peirce rules hold. These are exactly the facts used in the
proofs of Lemmas~\ref{cones:lem:face}, \ref{cones:lem:Pc}, Theorem~\ref{cones:thm:eq} and of
Theorem~\ref{cones:thm:rank} in the case $\rk V\le3$, together with (F1)--(F8); hence these
proofs apply verbatim.

(a)$\Rightarrow$(b). Let $M$ be non-associative with $\rk M\ge4$. We construct pairwise
orthogonal nonzero projections $c_1,c_2,c_3,c_4$ with sum $e$ and a symmetry $s$ with
$P(s)c_2=c_3$. Take nonzero orthogonal projections $q_1,\dots,q_4$ with sum $e$.
\emph{Case A: some $q_i$, say $q_1$, is not central.} Then
$M_{1/2}(q_1)=\bigoplus_{j\ne1}M_{1/2}(q_1)\cap M_{1/2}(q_j)\ne0$, say the term with $j=2$.
In the JBW-algebra $N=U_{q_1+q_2}M$ with unit $q_1+q_2$, Theorem~5.2.11 of \cite{HOS84}
writes $q_1=p_1+p_2$, $q_2=p_1'+p_2'$ with $p_1,p_1'$ exchanged by a symmetry of $N$ and
$c_N(p_2)\perp c_N(p_2')$; if $p_1=0$, then $c_N(q_1)\perp c_N(q_2)$, so $q_1$ would be central
in $N$, contradicting $M_{1/2}(q_1)\cap N\ne0$. Adding $e-q_1-q_2$ to the symmetry, we get
nonzero $c_2=p_1\le q_1$, $c_3=p_1'\le q_2$ exchanged by a symmetry $s$ of $M$, and we take $c_1=q_3$,
$c_4=q_4+(q_1-c_2)+(q_2-c_3)$. \emph{Case B: all $q_i$ are central.} Some summand
$U_{q_1}M$ is non-associative, hence contains a projection $p$ that is not central in it
(a JBW-algebra is generated by its projections \cite[4.2.3]{HOS84}, and the centre is
associative \cite[2.5.3]{HOS84}); \cite[Thm.~5.2.11]{HOS84} applied in $U_{q_1}M$ to $p$ and
$q_1-p$ gives, as in Case A, nonzero $c_2\le p$, $c_3\le q_1-p$ exchanged by a symmetry, and
we take $c_1=q_2$,
$c_4=q_3+q_4+(q_1-c_2-c_3)$. Put $u=2\{c_2\,s\,c_3\}$ (Jordan triple product). The
unital Jordan subalgebra generated by $s$ and $c_2$ contains $c_3=P(s)c_2$ and $u$, and is
special by the Shirshov--Cohn theorem \cite[2.4.14]{HOS84}; computing there, with
$c_2c_3=0$ and $sc_3s=c_2$, gives $u=c_2sc_3+c_3sc_2$, $c_2\circ u=c_3\circ u=\frac12u$
and $u^2=c_2+c_3$. Now let
\[
x=4c_1+3c_2+c_3+c_4,\qquad y=4c_1+2(c_2+c_3)+u+c_4 .
\]
As in the proof of Theorem~\ref{cones:thm:rank}, $y\in M_{++}$ (the element $2(c_2+c_3)+u$ has
spectrum $\{1,3\}$ in $U_{c_2+c_3}M$), the Peirce rules give
$P(x^{1/2})y=16c_1+P(x'^{1/2})y'+c_4$ with $x'=3c_2+c_3$, $y'=2(c_2+c_3)+u$, and by
Proposition~\ref{cones:prop:jb} the spectrum of $P(x'^{1/2})y'$ in $U_{c_2+c_3}M$ lies in $[1,9]$;
so $\lambda_{\max}=16=4\cdot4$, $\lambda_{\min}=1$, and equality holds. Finally, a direct
computation with the Peirce rules gives $T_xT_yc_2=6c_2+u$ and $T_yT_xc_2=6c_2+\frac32u$,
so $[T_x,T_y]c_2=-\frac12u\ne0$: $x$ and $y$ do not operator commute.
\end{proof}

\begin{remark}
In $V_\cH$ over an infinite-dimensional Hilbert space everything in
Section~\ref{cones:sec:lorentz} holds verbatim: the boosts are those of $\mathrm{SO}^+(1,\cH)$,
and equality holds iff the space parts are non-negatively proportional. For algebras of
infinite rank, such as $B(\cH)_{\rm sa}$, the extreme points of the spectrum need not be
eigenvalues. Thus the finite-dimensional criterion in Theorem~\ref{cones:thm:eq}, expressed using common extreme eigenvectors, does not extend to arbitrary infinite-rank algebras without additional hypotheses.
\end{remark}

\section{A characterisation of traces}

A positive linear functional $\varphi$ on $V$ is a \emph{trace} if
$\varphi(x\circ(y\circ z))=\varphi((x\circ y)\circ z)$; for
$V=\bigoplus_kV_k$ (simple components) the traces are exactly
$\varphi=\sum_k\alpha_k\tr_k$, $\alpha_k\ge0$ (for simple $V$ this is
\cite[Prop.~III.4.1]{FK94}; the direct-sum case follows by restricting to components). A symmetry is
$s\in V$ with $s^2=e$; then $P(s)$ is an automorphism of $V$ and $P(s)e=e$.

\begin{theorem}\label{cones:thm:trace}
Let $\varphi\ne0$ be a positive linear functional on $V$ and $p>0$. The following
are equivalent:
\begin{enumerate}
\item $\varphi$ is a trace;
\item $\varphi(P(s)x)=\varphi(x)$ for all symmetries $s$ and all $x\in V$;
\item $\varphi(a^p)\,\varphi(P(s)a^{-p})\ge\varphi(e)^2$ for all $a\in\Omega$ and all
symmetries $s$.
\end{enumerate}
\end{theorem}

\begin{proof}
(a)$\Rightarrow$(b): $P(s)x=2s\circ(s\circ x)-s^2\circ x=2s\circ(s\circ x)-x$, and by
associativity of $\varphi$, $\varphi(s\circ(s\circ x))=\varphi(s^2\circ x)=\varphi(x)$.
(b)$\Rightarrow$(c): $\varphi(P(s)a^{-p})=\varphi(a^{-p})$, and
$(u,v)\mapsto\varphi(u\circ v)$ is positive semidefinite, so by Cauchy--Schwarz with
$u=a^{p/2}$, $v=a^{-p/2}$, $u\circ v=e$: $\varphi(e)^2\le\varphi(a^p)\varphi(a^{-p})$.
(c)$\Rightarrow$(b): for $h\in V$ and a symmetry $s$ let
$f(\varepsilon)=\varphi(e^{p\varepsilon h})\varphi(P(s)e^{-p\varepsilon h})-\varphi(e)^2\ge0$.
Then $f(0)=0$, so $0=f'(0)=p\varphi(e)\bigl(\varphi(h)-\varphi(P(s)h)\bigr)$, and
$\varphi(e)>0$ (if $\varphi(e)=0$, then $0\le x\le\lambda_{\max}(x)e$ gives $\varphi=0$
on $\Omega$, which spans $V$).
(b)$\Rightarrow$(a): let $V_k$ be a simple component of rank $r_k\ge2$ and $c_i\perp c_j$
primitive in $V_k$. With $u$ from (F6), $s=(e-c_i-c_j)+u$ is a symmetry
($(e-c_i-c_j)\circ u=0$ by (F4), $u^2=c_i+c_j$) and $P(s)c_i=P(u)c_i=c_j$, the mixed
terms vanishing by (F4). Hence $\varphi(c_i)=\varphi(c_j)$; every primitive idempotent
of $V_k$ belongs to a Jordan frame of $V_k$, so $\varphi(c)=\varphi(e_k)/r_k$ for all
primitive $c\in V_k$, and by the spectral theorem $\varphi|_{V_k}=\frac{\varphi(e_k)}{r_k}\tr_k$.
Components of rank $1$ are $\R$, where every functional is a multiple of $\tr$.
\end{proof}

\begin{remark}
With $s=e$ only, (c) holds for every positive functional (Cauchy--Schwarz).
\end{remark}

\section{Open questions}

\begin{enumerate}[label=\textup{(\arabic*)}]
\item The equality case of Theorem~\ref{cones:thm:eq} in JBW-algebras of infinite rank, such
as $B(\cH)_{\rm sa}$, where the extreme points of the spectrum need not be eigenvalues.
\item Schatten contrasts in JBW-algebras with a faithful normal tracial state: is the
defect at most $\frac1p\ln(1/\varepsilon)$, with $\varepsilon$ the infimum of the traces of
nonzero projections?
\item For which other symmetric cones, besides the Lorentz cones, does the Wigner rotation
$k$ of Lemma~\ref{cones:lem:polar} have a geometric meaning, e.g.\ as a holonomy of the
Riemannian symmetric space $G(\Omega)/\mathrm{Aut}(V)$?
\item The analogue of Theorem~\ref{cones:thm:jordan} for the spectral contrasts $C_\Phi$, which
needs a defect (already $C_p(x^2)>2C_p(x)$); the region for fixed rank $\ge3$ is not known.
\end{enumerate}

\cleardoublepage

\end{document}